\documentclass[journal]{IEEEtran}
\usepackage{amsmath}
\usepackage{bbm}
\usepackage{amssymb}
\usepackage{mathtools}
\usepackage{amsthm}
\usepackage{booktabs}
\usepackage{amsmath}
\usepackage{fontawesome5}
\usepackage{algorithm,algorithmic}
\usepackage{hyperref}
\usepackage{bbm}
\usepackage{dsfont} 
\usepackage{subfigure}
\usepackage{xcolor}

\newtheorem{theorem}{Theorem}

\newtheorem{Lemma}{Lemma}
\newtheorem{corollary}{Corollary}
\newtheorem{definition}{Definition}

\newtheorem{assumption}[theorem]{Assumption}

\newcommand{\bd}{\boldsymbol}
\newcommand{\mc}{\mathcal}
\newcommand{\mb}{\mathbb}
\allowdisplaybreaks
\def\bx{\boldsymbol{x}}
\def\by{\boldsymbol{y}}
\def\bz{\boldsymbol{z}}
\def\bg{\boldsymbol{g}}
\def\bm{\boldsymbol{m}}

\def\bh{\boldsymbol{h}}

\def\bg{\boldsymbol{g}}
\def\bm{\boldsymbol{m}}
\def\bu{\boldsymbol{u}}

\DeclareMathOperator{\mV}{{\scriptstyle\boldsymbol{\mathcal{V}}}}

\DeclareMathOperator{\mM}{{\scriptstyle\boldsymbol{\mathcal{M}}}}

\DeclareMathOperator{\mU}{{\scriptstyle\boldsymbol{\mathcal{U}}}}
\DeclareMathOperator{\mG}{{\scriptstyle\boldsymbol{\mathcal{G}}}}
\DeclareMathOperator{\mX}{{\scriptstyle\boldsymbol{\mathcal{X}}}}
\DeclareMathOperator{\mY}{{\scriptstyle\boldsymbol{\mathcal{Y}}}}
\DeclareMathOperator{\mH}{{\scriptstyle\boldsymbol{\mathcal{H}}}}
\DeclareMathOperator{\mnH}{{\scriptstyle\mathcal{H}}}

\def\mc{\mathcal}

\def\bm{\boldsymbol{m}}
\def\bxi{\boldsymbol{\xi}}

\def\mE{\mathbb{E}}
\def\mR{\mathbb{R}}

\ifCLASSINFOpdf
\else
\fi

\begin{document}
\title{Communication-Efficient Decentralized Stochastic Minimax Optimization}

\author{Haoyuan Cai,
        Sulaiman A.~Alghunaim,
        and~Ali H.~Sayed,
        \thanks{Haoyuan Cai and Ali H.~Sayed are with the School of Engineering, École Polytechnique Fédérale de Lausanne, Switzerland. 
Emails: \{haoyuan.cai, ali.sayed\}@epfl.ch}
\thanks{Sulaiman A.~Alghunaim is with the Electrical Engineering Department, Kuwait University, Kuwait. 
Email: sulaiman.alghunaim@ku.edu.kw.}
}

\markboth{Journal of \LaTeX\ Class Files,~Vol.~14, No.~8, August~2015}%
{Shell \MakeLowercase{\textit{et al.}}: Bare Demo of IEEEtran.cls for IEEE Journals}

\maketitle

\begin{abstract}
\color{black}
We study decentralized stochastic nonconvex Polyak–Łojasiewicz minimax problems.
Motivated by local SGD in federated learning, we improve communication efficiency by allowing agents to perform several local updates before the communication stage. 
Existing work gives a communication bound of
$\mathcal O\Big(\frac{\kappa^3\varepsilon^{-2}}{(1-\lambda)^2}\Big)$ for local
decentralized gradient descent--ascent using
$\mathcal O\Big(\varepsilon^{-2}+\frac{\kappa^2\varepsilon^{-2}}{K(1-\lambda)}\Big)$
local updates per round. Here $\varepsilon$ is the target accuracy,
$1-\lambda$ is the spectral gap and $\kappa\ge1$ is the condition number. However,
too many local updates can leave available communication resources underused and increase local drift, as agents follow their own objectives for longer before exchanging information. We therefore pursue improved communication efficiency while keeping the number of local steps small. To this end,
we propose DiMA,  a structured local–global learning strategy that uses local updates to
build an accelerated momentum and perform decentralized learning through either standard gossip or accelerated mixing.

By exploiting the positive semidefiniteness of the mixing matrix, we derive complexity bounds for DiMA under different updating strategies: i) with standard gossip and $N=\mathcal{O}\Big(
\frac{\kappa(1-\lambda)}{K\varepsilon}
\Big)$ local steps, we achieve a communication complexity of $\mathcal{O}\Big(
\frac{\kappa^2\varepsilon^{-2}}{(1-\lambda)}
\Big)$, improving upon local
decentralized gradient descent--ascent; ii)
accelerated mixing further improves the spectral gap dependence under a suitable choice of $N$. Both strategies attain a total oracle complexity of \(\mathcal{O}(\kappa^3\varepsilon^{-3})\) under the corresponding accuracy conditions; and iii) with $N=1$, standard gossip yields a communication complexity of \(\mathcal{O}\Big(\frac{\kappa^3\varepsilon^{-3}}{K}\Big)\) with no explicit spectral gap factor in the stated accuracy regime. This bound exhibits linear speedup in \(K\), at the cost of worse dependence on \(\kappa\) and \(\varepsilon\).
Experiments on robust logistic regression demonstrate improved performance over several existing baselines.
\end{abstract}


\IEEEpeerreviewmaketitle

\vspace{-1em}
\section{Introduction}
\label{sec:intro}
\IEEEPARstart{M}{inimax} optimization has attracted significant attention in machine learning due to its broad applications, such as adversarial learning \cite{tu2019theoretical}, fair machine learning \cite{celis2019improved}, and robust image classification \cite{cheng2014minimax}.
In this work, we consider the following stochastic minimax optimization problem:
\begin{align} \label{distributed_min_max}
\min_{x\in \mathbb{R}^{M_1}}
&\max_{y \in \mathbb{R}^{M_2}} J(x,y)
 \triangleq \frac{1}{K} \sum_{k=1}^{K}
J_k(x,y), \notag
\\
\text{ where }
&
\ J_{k}(x,y) = \mathbb{E}_{\bd{\xi}_k \sim \mathcal{D}_{k}}[Q_{k}(x,y;\bd{\xi}_{k})].
\end{align}
Here, 
$x$ and $y$ are the optimization variables,
$\mE[\cdot]$ denotes the expectation operator,
\( K \) represents the number of agents (nodes, machines, workers, or clients), and \( Q_{k}(x,y;\bd{\xi}_k) \), \( \bd{\xi}_{k} \), and \( \mathcal{D}_{k} \) denote the loss function, sample, and data distribution specific to agent \( k \), respectively.
The global cost $J(x, y)$ is assumed to be $L_f$-smooth, and  $-J(x, y)$ is possibly nonconvex in $y$ 
but satisfies a 
$\nu$-Polyak Łojasiewicz (PL) condition in $y$ -- see eq.~\eqref{pl_eq}.
In multi-agent scenarios, when $K\ge 2$, there are two primary configurations: (i) federated/centralized setups, where all agents communicate with a central server \cite{mcmahan2017communication}, and (ii) decentralized setups, where agents are interconnected via a network and can only communicate with their immediate neighbors \cite{sayed2014adaptation}.
This work focuses on the decentralized case due to its scalability and ability to handle large-scale problems.

Solving a general minimax optimization problem is computationally intractable \cite{daskalakis2021complexity}.
This motivates the study of minimax problems with additional structure. One line of work formulates these problems as variational inequalities and establishes convergence under monotonicity or other structural conditions \cite{mertikopoulos2018optimistic,diakonikolas2021efficient}.
Another widely studied setting is nonconvex–strongly concave minimax optimization, which arises in applications such as robust training of neural networks \cite{nouiehed2019solving,lin2020gradient}.  The strong concavity condition ensures that the inner maximizer is unique for each $x$ whenever it exists. In the centralized setting, two-time-scale stochastic gradient descent-ascent (GDA) finds an $\varepsilon$-stationary point of the value function in $\mathcal{O}(\kappa^{3}\varepsilon^{-4})$ stochastic gradient evaluations under the smooth nonconvex-strongly concave assumptions \cite{lin2020gradient}, where $\kappa \ge 1$ is the condition number.
In this work, we study a decentralized nonconvex Polyak–Łojasiewicz (PL) minimax formulation, which relaxes the commonly used strong concavity assumption on the $y$-variable to a PL condition. 
\textcolor{black}{
When communication bandwidth is limited, exchanging information with neighbors after every update can be costly. Local updates reduce the frequency of communication by allowing each agent to take several steps using its own data between communication exchanges \cite{stich2018local}. Momentum offers a complementary approach of acceleration by reducing the number of stochastic oracle evaluations needed to reach an approximate stationary point \cite{cai2024accelerated,xian2021faster}. This leads to the following question. 
 }

\textbf{Q1:}
 \textcolor{black}{\em 
Can momentum further reduce the communication needed to reach 
$\varepsilon$-stationarity in decentralized stochastic minimax optimization with local updates?
}

{\color{black} This work answers the question affirmatively. Achieving this improvement, however, introduces new challenges:} algorithmically, it requires adapting communication-efficient strategies to minimax acceleration methods; theoretically, it requires controlling nontrivial errors arising from decentralization, communication protocols, and acceleration design.
We next review the related works.
\begin{table*}[t]
\color{black}
\centering
\renewcommand{\arraystretch}{1.9}
\setlength{\tabcolsep}{3pt}
\caption{\normalfont
Comparison of decentralized algorithms for online
nonconvex--strongly concave (PL) minimax optimization.
OS denotes \(\mE\|\nabla P(\bx)\|\leq\varepsilon\);
GS denotes
\(\mE\|\nabla_xJ(\bx,\by)\|\leq\varepsilon\) and
\(\mE\|\nabla_yJ(\bx,\by)\|\leq\varepsilon\).
Communication complexity counts rounds of exchange between neighbors. Oracle complexity counts oracle calls across all agents. For \textbf{DiMA}, the reported oracle complexity applies separately to stochastic gradient and Hessian–vector product (HVP) evaluations. For the first-order baselines, we count only stochastic gradient evaluations.
Here, \(1-\lambda\) is the network spectral gap and
\(\kappa\) is the condition number.
}
\label{table:comparison}
\begin{tabular}{cccccccc}
\toprule
\textbf{Algorithm}
& \shortstack{\textbf{Batch}\\\textbf{flexible?}}$^{\triangleleft}$
& \textbf{Acceleration?}
& \shortstack{\textbf{Local}\\\textbf{steps}}
& \shortstack{\textbf{Mixing}\\\textbf{operations}}$^{\square}$
& \textbf{Metric}
& \textbf{Comm. complexity}
& \textbf{Oracle complexity}
\\
\midrule

DM-HSGD \cite{xian2021faster}
& \faCheck & \faCheck
& $1$ & $4$ & OS
& $\mc{O}\Big(
\frac{\kappa^3\varepsilon^{-2}}{(1-\lambda)^2\min\{1,K\varepsilon\}}
\Big)$
& $\mc{O}\Big(
\frac{K\kappa^3\varepsilon^{-2}}{(1-\lambda)^2\min\{1,K\varepsilon\}}
\Big)$
\\
\hline

DM-GDA \cite{huang2023near}
& \faCheck & \faCheck
& $1$ & $4$ & OS
& $\mc{O}(\varepsilon^{-3})^{\dagger}$
& $\mc{O}(K\varepsilon^{-3})^{\dagger}$
\\
\hline

Dec-FedTrack \cite{ghiasvand2025robust}
& \faCheck & \faTimes
& $\mc{O}\Big( \varepsilon^{-2}+
\frac{\kappa^2\varepsilon^{-2}}{K(1-\lambda)}
\Big)$ & $4$ & OS
& $\mc{O}\Big(
\frac{\kappa^3\varepsilon^{-2}}{(1-\lambda)^2}
\Big)$
& $\mc{O}\Big(
\frac{K\kappa^3\varepsilon^{-4}}{(1-\lambda)^2}
+\frac{\kappa^5\varepsilon^{-4}}{(1-\lambda)^3}
\Big)$
\\
\hline

DGDA-VR \cite{zhang2024jointly}
& \faTimes & \faCheck
& $1$ & $4$ & OS/GS$^{\circ}$
& $\mc{O}\Big(
\frac{\kappa^2\varepsilon^{-2}}
{\min\{1/\kappa,(1-\lambda)^2\}}
\Big)$
& $\mc{O}\Big(
\frac{K\kappa^3\varepsilon^{-3}}
{\min\{1/\kappa,(1-\lambda)^2\}}
\Big)$
\\
\hline

DREAM \cite{chen2024efficient}
& \faTimes & \faCheck
& $1$
& $\mc{O}\Big(
\frac{\log K}{\sqrt{1-\lambda}}
\Big)$
& OS
& $\mc{O}\Big(
\frac{\kappa^2\varepsilon^{-2}\log K}
{\sqrt{1-\lambda}}
\Big)^{\ddagger}$
& $\mc{O}(\kappa^3\varepsilon^{-3})$
\\
\midrule

\multicolumn{8}{l}{\textbf{DiMA (Ours)}$^{\Diamond}$}
\\
\hline

DiMA $(\tau=0)$
& \faCheck & \faCheck
& $\mathcal{O}\Big(
\frac{\kappa(1-\lambda)\varepsilon^{-1}}{K}
\Big)$
& $4$ & OS/GS
& $\mc{O}\Big(
\frac{\kappa^2\varepsilon^{-2}}{1-\lambda}
\Big)$
& $\mc{O}\Big(
\kappa^3\varepsilon^{-3}
\Big)$
\\
\hline
DiMA $(\tau=0)$
& \faCheck & \faCheck
& $1$ & $4$ & OS/GS
& $\mc{O}\Big(
\frac{\kappa^3\varepsilon^{-3}}
{K}
\Big)$
& $\mc{O}\Big(
\kappa^3\varepsilon^{-3}
\Big)$
\\
\hline
DiMA $(\tau=1)$
& \faCheck & \faCheck
& $\mathcal{O}\Big(\frac{\kappa\varepsilon^{-1}}{K}\Big)$
& $\mathcal{O}\Big(\frac{1}{\sqrt{1-\lambda}}\Big)$ & OS/GS
& $\mc{O}\Big(
\frac{\kappa^2\varepsilon^{-2}}{\sqrt{1-\lambda}}
\Big)$
& $\mc{O}(\kappa^3\varepsilon^{-3})$
\\
\hline
DiMA $(\tau=1)$
& \faCheck & \faCheck
& $1$ & $\mathcal{O}\Big(\frac{1}{\sqrt{1-\lambda}}\Big)$ & OS/GS
& $\mc{O}\Big(
\frac{\kappa^3\varepsilon^{-3}}
{K\sqrt{1-\lambda}}
\Big)$
& $\mc{O}(\kappa^3\varepsilon^{-3})$
\\
\bottomrule
\end{tabular}

\vspace{0.5em}
\begin{minipage}{\textwidth}
\footnotesize
$^{\triangleleft}$
A check mark indicates that the algorithm does not require large-batch gradient refreshes for variance reduction.

$^{\square}$
The number of mixing operations per outer iteration,
counting each network variable separately.

$^{\dagger}$
The cited complexity does not explicitly report the
dependence on $1-\lambda$, which cannot be easily derived.

$^{\circ}$ For DGDA-VR, GS follows by smoothness from the bounds on
the value gradient and the dual response error
\cite[Thms.~1--2]{zhang2024jointly}. The oracle complexity multiplies the
per-agent bound by $K$.

$^{\ddagger}$ DREAM reports a $\log K$ factor in its main communication bound. Their Appendix~G \cite{chen2024efficient} removes this factor when the gradient refresh batch size is at least $K$.

$^{\Diamond}$ 
The table reports the complexity bounds for $\varepsilon\le\kappa(1-\lambda)/K$ with standard gossip ($\tau=0$) and for $\varepsilon\le\kappa/K$ with accelerated mixing ($\tau=1$). See Corollaries~\ref{corollary:gossip} and~\ref{corollary:fastmixing} for the detailed results.
\end{minipage}
\vspace{-2em}
\end{table*}

\vspace{-1em}
\subsection{Related works}
\textbf{Nonconvex minimax optimization.}
Existing works 
have explored both 
single- and double-loop algorithms 
for solving stochastic nonconvex strongly-concave 
minimax problems; see, e.g., references \cite{luo2020stochastic,lin2020gradient, yang2022faster}. A single-loop structure is often preferred for its simplicity but typically requires a two-time-scale step-size policy to achieve stationarity of the value function $P(x)=\max_y J(x,y)$. The works 
\cite{lin2020gradient, luo2020stochastic,zhang2020single, ostrovskii2021efficient, kong2021accelerated, zhang2022sapd, xu2023unified} proposed algorithms for solving nonconvex-concave 
minimax optimization. 
They explored the nonsmooth optimization techniques of the Moreau envelope to approximate a stationary point of the smooth surrogate cost. 
The works
\cite{nouiehed2019solving, yang2022faster, huang2023enhanced, cai2024accelerated, yang2020global} studied nonconvex-PL minimax optimization problems.
The PL condition relaxed the strong concavity assumption on the $y$-variable, thereby accommodating certain nonconcave structures \cite{karimi2016linear}. Existing works mainly focus on finding an optimization stationary point or game stationary point, whereas we adopt a unified metric and establish convergence under both metrics; see Definitions \ref{main:definition:OS}, \ref{main:definition:GS} and Lemma \ref{main:lemma:unified_metrics}. Moreover, we consider the more challenging decentralized setting with an emphasis on communication efficiency, which requires a novel algorithmic design and a more involved analysis

\textbf{Decentralized nonconvex minimax optimization.}
Decentralized minimax optimization has been studied in 
\cite{deng2021local, wu2024solving, xian2021faster, huang2023near,sharma2022federated, cai2024diffusion, chen2024efficient, shen2024stochastic,zhang2025federated, zhang2023communication, cai2025dama}.
A common framework involves a central coordinating entity, often termed as a federated minimax problem, studied 
by \cite{deng2021local, wu2024solving, sharma2022federated, shen2024stochastic,zhang2025federated,zhang2025federated, zhang2023communication}.
 In contrast, decentralized minimax problems provide a more flexible framework and are more robust against server-agent failures; they are studied in \cite{xian2021faster, huang2023near, chen2024efficient,cai2024diffusion, cai2025dama}.  Table~\ref{table:comparison} gives the communication bounds for the methods compared here.

\textbf{Acceleration technique.}
Acceleration techniques 
such as variance reduction \cite{nguyen2017sarah, johnson2013accelerating, reddi2016stochastic}, momentum \cite{cutkosky2019momentum, levy2021storm, tran2022better, cai2024accelerated}, and implicit gradient transport \cite{cutkosky2020momentum, arnold2019reducing}, have been widely studied for single-variable minimization problems.
Noticeably, variance reduction approaches such as 
SARAH, SVRG have been adapted to minimax optimization problems \cite{luo2020stochastic, gao2022decentralized, yang2020global, gao2022decentralized,chen2024efficient}.
STORM momentum has also been adapted to solve the minimax problem \cite{xian2021faster, huang2023enhanced, shen2024stochastic}.
In the context of nonconvex-strongly concave minimax optimization, 
acceleration techniques have been proposed to achieve improved convergence performance on the condition number $\kappa$, including the inexact proximal-point-type method \cite{lin2020near,yang2020catalyst, yang2022faster}. Unlike these works, we devise decentralized accelerated methods that utilize second-order Lipschitz information to improve the sample complexity on $\varepsilon$.

\textbf{Communication-efficient learning.}
Local SGD is a communication-efficient learning paradigm widely used in federated learning. 
It reduces communication overhead by allowing agents to perform multiple local updates without synchronizing after every gradient step \cite{stich2018local}.
This local SGD method has inspired fully decentralized optimization \cite{liu2024decentralized, nguyen2023performance, alghunaim2024local} and has been extended to local GDA methods \cite{ghiasvand2025robust, li2024fast}. However, without acceleration, local-SGD-type methods may require an excessive number of local steps to achieve optimal communication complexity, potentially underusing available communication resources. Their sample complexity also remains to be improved. Another popular approach is model compression \cite{zhang2023communication}; however, compressed communication may incur information loss. We therefore focus on local updates with momentum acceleration, which reduces the communication burden while preserving lossless information exchange.

\vspace{-1em}

\subsection{Contributions}
We develop accelerated communication-efficient minimax algorithms by leveraging momentum acceleration with local updates. We use a momentum acceleration scheme based on partial second-order information \cite{cai2024accelerated}. While this information has been widely used in minimization \cite{cutkosky2020momentum, tran2022better, carmon2018accelerated}, minimax optimization \cite{li2022tiada}, and bilevel programming \cite{dagreou2022framework}, its potential in decentralized minimax optimization remains largely unexplored.
Below, we summarize our main contributions.

\begin{itemize}
\item[$\bullet$]
\textcolor{black}{
We develop \textbf{DiMA} (Local \textbf{Di}ffusion \textbf{M}inimax \textbf{A}lgorithm) by combining local updates with momentum acceleration. Between communication stages, each agent takes several normalized gradient descent–ascent (GDA) steps to limit local drift, while accumulating stochastic gradients and Hessian–vector-product (HVP) along the local trajectory. These accumulated quantities are used to update its bias-corrected momentum. Agents then track the average corrected momentum across the network and combine their updated models with those of their neighbors. Both operations can use standard gossip, following the adapt-then-combine (ATC) diffusion strategy \cite{sayed2014adaptation,sayed2022inference}, or accelerated mixing \cite{ye2023multi}.
}
 
\color{black}
\item[$\bullet$]
We establish convergence guarantees for \textbf{DiMA} with different numbers of local steps and two mixing strategies.  Table~\ref{table:comparison} compares our guarantees with existing results. Main theoretical contributions are listed below.
\begin{enumerate}
\item[(i)] \textbf{Multiple local updates with standard gossip.} By exploiting the positive semidefiniteness (PSD) of the mixing matrix, we show that, when $K\varepsilon\le\kappa(1-\lambda)$, choosing
$
N=\left\lceil\frac{\kappa(1-\lambda)\varepsilon^{-1}}{K}\right\rceil
$
gives a communication complexity of
$\mathcal{O}\big(\frac{\kappa^2\varepsilon^{-2}}{1-\lambda}\big)$
and a total oracle complexity of
\(\mathcal{O}(\kappa^3\varepsilon^{-3})\), counting stochastic gradient and HVP evaluations across all agents. An HVP can be computed without forming the Hessian, at a computational cost comparable to that of a gradient evaluation \cite{pearlmutter1994fast,tran2022better}.
Compared with Dec-FedTrack \cite{ghiasvand2025robust}, \textbf{DiMA} requires fewer local updates and improves the dependence of communication complexity on both \(\kappa\) and \(1-\lambda\), thereby addressing question \textbf{Q1}. \textbf{DiMA} also improves the gradient oracle complexity of Dec-FedTrack, with no spectral gap factor dependency under the regime $K\varepsilon\le\kappa(1-\lambda)$. Among other methods using standard gossip, \textbf{DiMA} improves the $\varepsilon$ dependence of the communication bound over DM-HSGD \cite{xian2021faster} and the spectral gap dependence over DGDA-VR \cite{zhang2024jointly}.

\item[(ii)] \textbf{Multiple local updates with accelerated mixing.}
When $K\varepsilon\le\kappa$, choosing $N=\mathcal{O}\Big(\frac{\kappa\varepsilon^{-1}}{K}\Big)$ gives a communication complexity of \(\mathcal{O}\big(\frac{\kappa^2\varepsilon^{-2}}{\sqrt{1-\lambda}}\big)\) and a total oracle complexity of \(\mathcal{O}(\kappa^3\varepsilon^{-3})\). 
The communication result includes \(\mathcal{O}\big(\frac{1}{\sqrt{1-\lambda}}\big)\) gossip operations per outer iteration for momentum tracking and model mixing.
DREAM’s general  communication bound includes a $\log K$ factor, which its refined bound removes when the gradient refresh batch size is at least $K$ per agent \cite{chen2024efficient}. \textbf{DiMA} with accelerated mixing matches the order of that refined communication bounds through local updates, without recomputing the gradient estimate on a large batch of samples.

\item[(iii)] \textbf{Linear speedup with single local updates.}
Under the same accuracy conditions as in (i) and (ii), respectively, setting $N=1$ gives a communication complexity of \(\mathcal{O}\big(\frac{\kappa^3\varepsilon^{-3}}{K}\big)\) with standard gossip and \(\mathcal{O}\big(\frac{\kappa^3\varepsilon^{-3}}{K\sqrt{1-\lambda}}\big)\) with accelerated mixing. Both results exhibit linear speedup in the number of agents $K$, while the total oracle complexity remains \(\mathcal{O}(\kappa^3\varepsilon^{-3})\). Compared with the multiple local step choices above, the single-step bounds of standard gossip exhibit linear speedup in \(K\) and improved spectral gap dependence in communication, at the cost of worse dependence on \(\kappa\) and \(\varepsilon\).
\end{enumerate}

\color{black}

\item[$\bullet$]
\color{black}
 Our theoretical analysis is nontrivial and relies on a carefully constructed potential function.
A major challenge is that agents normalize their directions locally, so the average of the normalized directions need not equal the normalized network average. 
We exploit the PSD mixing matrix to bound iterate differences between successive communication stages. This produces a tighter inverse spectral gap factor
and sharpens the bound on the momentum estimation error, further improving convergence guarantees for game stationarity and optimization stationarity.
Normalization permits a step-size ratio of \(\frac{\eta_x}{\eta_y}=\mathcal{O}\big(\frac{1}{\kappa}\big)\), compared with \(\mathcal{O}\big(\frac{1}{\kappa^2}\big)\) in the two-timescale GDA analysis of \cite{lin2020gradient}. This allows a larger primal step relative to the dual step when \(\kappa\) is large. The relaxed ratio alone does not imply a better oracle complexity.
\color{black}

\end{itemize}
\textbf{Notation.}
We use normal, bold, and calligraphic fonts, e.g., $x, \bd{x}$, and $\bd{\scriptstyle\mathcal{X}}$,
to denote deterministic, stochastic, and network-extended quantities, respectively.
$\mc{O}(\cdot)$ is standard complexity notation.
We use  $z=\mbox{\rm col}\{x,y\}$ to denote concatenated vector of $x$ and $y$. 
$I_{M}$ denotes $M$-dimensional identity matrix, where $M=M_1+M_2$.
The symbol $\otimes$  denotes the Kronecker product. The symbol $\|\cdot\|$ denotes the Euclidean norm for vectors and the induced spectral norm for matrices. 
\textcolor{black}{$\lceil a\rceil$ denotes $a$ rounded up to the nearest integer.}
The symbol
${\rm col}\{a,b,c\}$ or 
$[a;b;c]$ means concatenating the vectors $a, b, c$ along the column dimension.
\textcolor{black}{Let $\mathcal N(v)=v/\|v\|$ for $v\ne0$ and $\mathcal N(0)=0$.}
We often use the notation, e.g.,  $\bx_{k,i,n}$, to denote a local quantity, where the subscripts $k,i$, and $n$  represent the agent,  communication round, and local update index, respectively.

\section{Technical Preliminaries}

In this section,
we present technical preliminaries.
We first introduce some technical assumptions for smooth optimization based on the following notions.
\begin{definition}[\textbf{$L_f$-smooth}]
\label{definition:lipchitzgradient}
The function $f(z)$ is $L_f$-smooth if
    \begin{align}
        \|\nabla f(z_1) - \nabla f(z_2)\| \le L_f \|z_1 - z_2\|, \ \forall z_1, z_2 .
    \end{align}
\end{definition}
\begin{definition}[\textbf{$L_h$-second-order Lipschitz}]
 The function $f(z)$  is $L_h$-second-order Lipschitz if
    \begin{align}
        \|\nabla^2 f(z_1) - \nabla^2 f(z_2)\| \le L_h \|z_1 - z_2\|, \ \forall z_1, z_2. 
    \end{align}
\end{definition}
\begin{assumption}[\textbf{Smooth function}]
\label{assumption:lipchitzgradient}
We assume each local cost 
    $J_k(x,y)$ is $L_f$-smooth and $L_h$-second-order Lipschitz  with respect to  $z=\mbox{\rm col}\{x,y\}$.
\end{assumption}
The $L_f$-smooth assumption is standard in the analysis of smooth minimax optimization \cite{lin2020gradient, luo2020stochastic, yang2022faster}. 
\textcolor{black}{The second condition, requiring an \(L_h\)-Lipschitz Hessian, controls the approximation error in our accelerated momentum estimator \cite{arnold2019reducing,tran2022better,cutkosky2020momentum,cai2024accelerated}.}
Given that our algorithms access stochastic oracles only, we adopt the following assumption.
\begin{assumption}[\textbf{Bounded variance}]
\label{assumption:boundvariance}
We assume
the local stochastic gradient and stochastic Hessian are unbiased and have bounded variance
    \begin{equation}
        \begin{aligned}
        &\mE_{\bxi_k \sim \mc{D}_k} [\nabla_z Q_k(x,y;\bxi_k)] = \nabla_z J_k(x, y),  \\
        & \mE_{\bxi_k \sim \mc{D}_k} [\nabla^2_z Q_k(x,y;\bxi_k)] = \nabla^2_z J_k(x, y),  \\
        &\mE_{\bxi \sim \mc{D}_k} \|\nabla_z Q_k(x,y;\bxi_k) -  \nabla_z J_k(x, y) \|^2 \le  \sigma^2,  \\
      &\mE_{\bxi \sim \mc{D}_k} \|\nabla^2_z Q_k(x,y;\bxi_k) -  \nabla^2_z J_k(x, y) \|^2 \le  \sigma^2_h.
    \end{aligned}
    \end{equation}
{\color{black}Here $z=\operatorname{col}\{x,y\}$. All initialization samples $\bxi_{k,-1}$ and all samples $\bxi_{k,i,n}$ are mutually independent, and each sample drawn by agent \(k\) follows distribution \(\mathcal D_k\). The gradient and HVP at a given local step use the same sample.}
\end{assumption}
Note that an assumption about the stochastic Hessian is made; we only limit our queries to {\em partial} second-order information, which can be obtained efficiently, as will be elaborated soon.
\begin{assumption}[\textbf{Cost function}]
\label{assumption:riskfunction}
 The cost function $J(x, y)$ is nonconvex in $x$ and $-J(x, y)$ is $\nu$-PL in $y$, i.e., 
\begin{align} \label{pl_eq}
\|\nabla_y J(x, y)\|^2 \ge 2\nu [\max_{y} J(x, y) - J(x, y)],
\end{align}
for a constant $0<\nu\le L_f$. Write $\kappa=L_f/\nu\ge1$.
\color{black}
For every $x$, the set
$Y^{o}(x)=\operatorname*{argmax}_yJ(x,y)$ is nonempty
and $P(x)=\max_yJ(x,y)$ is finite. Moreover,
$P^\star=\inf_xP(x)>-\infty$. 
\end{assumption}
We remark that the $\nu$-PL condition is critical for ensuring the smoothness of the value function $P(x)$ implied by the Danskin-type lemma \cite{nouiehed2019solving}. This enables a smooth analysis without relying on nonsmooth optimization tools, such as the Moreau envelope, while allowing us to establish a sharper convergence rate. Moreover, we note that the PL condition is weaker than strong concavity \cite{karimi2016linear, nouiehed2019solving}.

We next review two commonly used convergence metrics in minimax optimization, namely,  optimization stationarity (OS) and game stationarity (GS) \cite{razaviyayn2020nonconvex,li2026smoothing}.
\begin{definition}
[\textbf{$\varepsilon$-Optimization stationarity}]
\label{main:definition:OS}
A point $\bx$ is an $\varepsilon$-optimization stationary point of the value function $P(x) = \max_y J(x,y)$ if $\mE\|\nabla P(\bx)\| \le \varepsilon$.
\end{definition}
The OS notion has been adopted in prior works, e.g.,  \cite{lin2020gradient, luo2020stochastic, xian2021faster}, and is well suited to settings such as adversarial learning \cite{tu2019theoretical}, where minimizing the value function $P(x)$ corresponds to finding a model robust to the worst-case cost. 

\begin{definition}[\textbf{$\varepsilon$-Game stationarity}]
\label{main:definition:GS}
The point $(\bx, \by)$ is an $\varepsilon$-game stationary point if $
\mE\|\nabla_x J(\bx, \by)\| \le \varepsilon$ and $\mE\|\nabla_y J(\bx, \by)\| \le \varepsilon$.
\end{definition}
The GS notion has been used 
in centralized works \cite{shen2024stochastic,yang2022faster, li2022tiada}
and is
suitable when the minimax optimization is explicitly viewed as a two-player game, e.g., a generative adversarial network.
 In fact, under the structured nonconvex-PL formulation, these two convergence notions can be established simultaneously by bounding the following metrics.
\begin{Lemma}
\label{main:lemma:unified_metrics}
If the function $J(x,y)$ is $L_f$-smooth,  nonconvex in $x$,  and $\nu$-PL in $y$, then
the point $(\bx, \by)$ 
satisfies both $\varepsilon$-OS/GS notions if 
\begin{align}
\mE\|\nabla P(\bx)\| +  2 L_f \mE\|\by-\by^o(\bx)\| \le \varepsilon,
\label{main:convergence_criterion}
\end{align}
\color{black}
where $\by^o(\bx)\in
\operatorname*{argmin}_{\boldsymbol{v} \in Y^o(\bx)}\|\by-\boldsymbol{v}\|$.
\end{Lemma}
\begin{proof}
It is straightforward to see that $\mE\|\nabla P(\bx)\| \le \varepsilon$ under \eqref{main:convergence_criterion}. On the other hand, we have
\begin{align}
&\mE\|\nabla_x J(\bx,\by)\|
+\mE\|\nabla_y J(\bx,\by)\| \notag\\
&\le \mE\|\nabla_x J(\bx,\by)-\nabla_x J(\bx,\by^o(\bx))+\nabla_x J(\bx,\by^o(\bx))\|
\notag\\
&\quad 
+\mE\|\nabla_y J(\bx,\by) - \nabla_y J(\bx,\by^o(\bx))\| \notag \\
&\le \mE\|\nabla P(\bx)\|
+2L_f\mE\|\by - \by^o(\bx)\| \le \varepsilon,
\end{align}
where the last inequality follows from condition \eqref{main:convergence_criterion}, $L_f$-smooth assumption and triangle inequality.
\end{proof}

\vspace{-1em}
\section{Algorithm Development}
In this section, we first introduce the core ideas behind our algorithm design and then present our decentralized communication-efficient algorithms.
\vspace{-1em}
\subsection{Local update and momentum acceleration}
A key bottleneck in decentralized optimization is limited communication bandwidth. We therefore employ two complementary strategies to alleviate this issue, i.e.,
(i) acceleration techniques that enable convergence in fewer communication rounds \cite{cutkosky2019momentum,tran2022better,cai2024accelerated}, and
(ii) local updates that directly reduce communication overhead by performing multiple updates between gossiping steps \cite{nguyen2023performance, liu2024decentralized, alghunaim2024local}.
We  use them jointly to boost communication efficiency.
We introduce these techniques below.

$\bullet$ \textbf{Momentum acceleration}.
Momentum-based acceleration strategies have shown empirical success across various applications \cite{ xian2021faster, tran2022better}. 
In this work, we adopt a momentum scheme similar to that in \cite{cai2024accelerated}. Notably, extending this momentum scheme to the decentralized setting is nontrivial due to new challenges introduced by decentralization, communication protocols, and nonlinear learning steps.
Momentum acceleration scheme performs updates using the moving average of the current gradient direction and a {\em corrected} past momentum direction \cite{cai2024accelerated}.
Let $\bm_{x,i}, \bm_{y,i}$ be the momentum direction at iteration $i$, it runs $\forall i \ge 0$:

{\footnotesize
\begin{align}
&\begin{bmatrix}
    \bm_{x,i}\\
    \bm_{y,i}
\end{bmatrix}
=
(1-\beta)
\Bigg(\begin{bmatrix}
 \bm_{x,i-1} \\
\bm_{y,i-1}
\end{bmatrix} +
\notag \\
& 
\underbrace{\begin{bmatrix}
\nabla^2_x Q(\bd{x}_{i}, \bd{y}_{i};\bd{\xi}_i) &\nabla^2_{xy}Q(\bd{x}_{i}, \bd{y}_{i};\bd{\xi}_i) \\
\nabla^2_{yx} Q(\bd{x}_{i}, \bd{y}_{i};\bd{\xi}_i) &
\nabla^2_{y} Q(\bd{x}_{i}, \bd{y}_{i};\bd{\xi}_i) 
\end{bmatrix}
\begin{bmatrix}
\bx_{i} - \bx_{i-1} \\
\by_{i} - \by_{i-1}
\end{bmatrix}}_{\text{Momentum Correction}}
\Bigg) \notag \\
& +
\beta
\begin{bmatrix}
    \nabla_x Q(\bx_i,\by_i;\bxi_i)
    \\
    \nabla_y Q(\bx_i,\by_i;\bxi_i)
\end{bmatrix},
\label{main:corrected_momentum}
\end{align}
}
where $\beta$ is the smoothing factor.
\textcolor{black}{The Hessian–vector product (HVP) in \eqref{main:corrected_momentum} approximates the change in the stochastic gradient between successive iterates. Both the gradient and the HVP are evaluated at the current iterate; the previous iterate is used only to form the displacement. This correction therefore requires no derivative evaluation at the previous iterate on the fresh sample, a distinction that motivates HVP-based variance reduction \cite{arjevani2020second}.
When the same fresh sample can be used at both iterates, STORM \cite{cutkosky2019momentum} provides a first-order alternative. It computes the gradient difference directly, avoiding the Taylor remainder of the HVP approximation. Here, we establish communication and oracle complexity guarantees for the HVP formulation of \textbf{DiMA}. We leave the analysis of whether a STORM variant attains the same guarantees to future work.}

We define a stochastic oracle as follows.
\begin{definition}[\textbf{Stochastic Hessian-vector-product}]
Given a loss function $Q_k(x,y;\bxi_k)$, the stochastic HVP at two points $(x_1, y_1), (x_2, y_2)$ and random sample $\bxi_k$
is defined as follows: 
\begin{align}
\textbf{hvp}&\Big((x_1, y_1), (x_2, y_2); \bxi_k\Big)  \\
\triangleq 
&\begin{bmatrix}
\nabla^2_x Q_k(x_1, y_1; \bxi_k)&\nabla^2_{xy} Q_k(x_1, y_1; \bxi_k)\\
\nabla^2_{yx} Q_k(x_1, y_1; \bxi_k)&
\nabla^2_y Q_k(x_1, y_1; \bxi_k)
\end{bmatrix}
\begin{bmatrix}
x_1 - x_2\\
y_1 - y_2
\end{bmatrix}.\notag
\end{align} 
\end{definition}
We note that, in many settings, querying an HVP oracle has a cost comparable to that of querying a stochastic gradient oracle \cite{pearlmutter1994fast}. For instance, PyTorch provides automatic differentiation tools for efficiently computing exact HVPs and Jacobian-vector products \cite{paszke2017automatic}. In cases with cost functions with explicit structure, such as logistic regression, the Hessian and Jacobian often admit simple outer-product forms and can be computed efficiently.

$\bullet$ \textbf{Local updates with normalization}.
Traditional distributed optimization methods synchronize agents after each update \cite{sayed2014adaptation}, leading to high communication overhead.
Local update strategies reduce this cost by allowing agents to perform multiple local updates between two gossip rounds \cite{nguyen2023performance, liu2024decentralized, alghunaim2024local}.
As such, a simple strategy is to aggregate local models periodically. However, naively skipping communication and periodically aggregating models may lead to suboptimal performance \cite{liu2024decentralized}. As shown in \cite{liu2024decentralized}, leveraging accumulated local gradient information to update the model can improve stochastic error bounds; in particular, compared with periodic aggregation methods, the second dominant term can be reduced by a factor proportional to the number of local steps.
Nevertheless, excessive or improperly designed local updates can cause local models to drift toward their own local optima, making consensus among agents difficult to achieve. To address this issue, we incorporate normalization into the local learning stage, which tightly controls the drift of each local model from its initial state.

\begin{algorithm}[t]
\color{black}
\caption{Mix(${\rm col} \{\bd{v}_k\}_{k=1}^K;A,r,\tau$)}
\label{alg:mix}
\begin{algorithmic}[1]
\footnotesize
\STATE Choose the mixing mode:
\[
(R,\eta_a)=
\begin{cases}
(1,0), & \tau=0 \quad \text{(Gossip)},\\[1mm]
\bigl(r,(1+\sqrt{1-\lambda^2})^{-1}\bigr),
& \tau=1 \quad \text{(FastMix)}.
\end{cases}
\]
where $\lambda=\|A-\mathds{1}_K \mathds{1}_K^\top/K\|_2$.
\STATE Initialize
$\bd{v}_k^{(-1)}=\bd{v}_k^{(0)}=\bd{v}_k$,
$\forall k\in[K]$.
\FOR{\textbf{gossip round} $s=0,\ldots,R-1$}
\FOR{\textbf{each agent} $k$ in parallel}
\STATE Send $\bd{v}_k^{(s)}$ to neighbors and receive
${\rm col}\{\bd{v}_{\ell}^{(s)}\}_{\ell\in\mathcal{N}_k}$.
\STATE Update
\[
\bd{v}_k^{(s+1)}
=(1+\eta_a)\sum_{\ell\in\mathcal{N}_k}
a_{k\ell}\bd{v}_{\ell}^{(s)}
-\eta_a\bd{v}_k^{(s-1)}.
\]
\ENDFOR
\ENDFOR
\STATE \textbf{Return} ${\rm col}\{\bd{v}_k^{(R)}\}_{k=1}^K$.
\end{algorithmic}
\end{algorithm}


\begin{algorithm}[t]
\caption{Local \textbf{Di}ffusion \textbf{M}inimax
\textbf{A}lgorithm (\textbf{DiMA})}
\label{DiMA}
\begin{algorithmic}[1]
\footnotesize
\renewcommand{\REQUIRE}{\item[\textbf{Initialize:}]}
\REQUIRE
Choose shared deterministic models $x_0,y_0$. For every
$k\in[K]$ and $n=0,\ldots,N-1$, set
$\bx_{k,0}=\bx_{k,-1}=\bx_{k,-1,n}=x_0$ and
$\by_{k,0}=\by_{k,-1}=\by_{k,-1,n}=y_0$.
Set initial direction $\bu^x_{k,-1}={\bd{m}}^x_{k, -1}=\nabla_x Q_k(\bz_{k,-1};\bxi_{k,-1}), \bu^y_{k,-1}={\bd{m}}^y_{k,-1} =\nabla_y Q_k(\bz_{k,-1};\bxi_{k,-1}), \forall k \in [K]$.
Step sizes $\eta_x,\eta_y,\bar{\eta}_x,\bar{\eta}_y>0$,
momentum parameter $\beta\in(0,1]$, 
\textcolor{black}{local step $N$, and gossip steps $r\geq1$. Mixing subroutine $\tau\in\{0,1\}$: $\tau=0$ for standard gossip and $\tau=1$ accelerated mixing.
}

\FOR{\textbf{iteration} $i=0,\ldots,T-1$}
\FOR{\textbf{each agent} $k$ in parallel}
\STATE \underline{\texttt{Initialize local values}}
\[
\bx_{k,i,0}=\bx_{k,i},\qquad
\by_{k,i,0}=\by_{k,i},\qquad
[\bg^x_{k,i};\bh^x_{k,i};\bg^y_{k,i};\bh^y_{k,i}]=0.
\]

\FOR{\textbf{local step} $n=0,\ldots,N-1$}
\STATE \underline{\texttt{Update local model}}
\begin{align*}
\bg^x_{k,i,n}
&=\nabla_xQ_k(\bx_{k,i,n},\by_{k,i,n};\bxi_{k,i,n}),\\
\bg^y_{k,i,n}
&=\nabla_yQ_k(\bx_{k,i,n},\by_{k,i,n};\bxi_{k,i,n}),\\
\bx_{k,i,n+1}
&=\bx_{k,i,n}
-\bar{\eta}_x\mathcal{N}(\bg^x_{k,i,n}),\\
\by_{k,i,n+1}
&=\by_{k,i,n}
+\bar{\eta}_y\mathcal{N}(\bg^y_{k,i,n}).
\end{align*}

\STATE \underline{\texttt{Accumulate local stochastic gradient}}
\[
[\bg^x_{k,i};\bg^y_{k,i}]
\leftarrow
[\bg^x_{k,i};\bg^y_{k,i}]
+\frac{1}{N}[\bg^x_{k,i,n};\bg^y_{k,i,n}].
\]

\STATE \underline{\texttt{Accumulate momentum correction}}
\begin{align*}
&[\bh^x_{k,i,n};\bh^y_{k,i,n}] \notag\\
&=\operatorname{\textbf{hvp}}\Big(
(\bx_{k,i,n},\by_{k,i,n}),
(\bx_{k,i-1,n},\by_{k,i-1,n});
\bxi_{k,i,n}\Big),\\
&[\bh^x_{k,i};\bh^y_{k,i}]
\leftarrow
[\bh^x_{k,i};\bh^y_{k,i}]
+\frac{1}{N}[\bh^x_{k,i,n};\bh^y_{k,i,n}].
\end{align*}
\ENDFOR

\color{black}
\STATE \underline{\texttt{Momentum update}}
\begin{align*}
\bm^x_{k,i}
&=(1-\beta)(\bm^x_{k,i-1}+\bh^x_{k,i})
+\beta\bg^x_{k,i},\\
\bm^y_{k,i}
&=(1-\beta)(\bm^y_{k,i-1}+\bh^y_{k,i})
+\beta\bg^y_{k,i}.
\end{align*}
\ENDFOR
\color{black}
\STATE \underline{\texttt{Momentum tracking}}
\begin{align*}
\{\bu^x_{k,i}\}_{k=1}^K
&=\operatorname{Mix}\left({\rm col}
\{\bu^x_{k,i-1}+\bm^x_{k,i}-\bm^x_{k,i-1}\}_{k=1}^K;
A,r,\tau\right),\\
\{\bu^y_{k,i}\}_{k=1}^K
&=\operatorname{Mix}\left({\rm col}
\{\bu^y_{k,i-1}+\bm^y_{k,i}-\bm^y_{k,i-1}\}_{k=1}^K;
A,r,\tau\right).
\end{align*}

\STATE \underline{\texttt{Model update}}
\begin{align*}
\{\bx_{k,i+1}\}_{k=1}^K
&=\operatorname{Mix}\left({\rm col}
\left\{\bx_{k,i}
-\eta_x\mathcal{N}(\bu^x_{k,i})\right\}_{k=1}^K;
A,r,\tau\right),\\
\{\by_{k,i+1}\}_{k=1}^K
&=\operatorname{Mix}\left({\rm col}
\left\{\by_{k,i}
+\eta_y\mathcal{N}(\bu^y_{k,i})\right\}_{k=1}^K;
A,r,\tau\right).
\end{align*}
\ENDFOR
\end{algorithmic}
\end{algorithm}

\vspace{-1em}
\subsection{Local Diffusion Minimax Algorithm}

We now introduce the communication-efficient local-global learning strategy 
for solving problem \eqref{distributed_min_max}.
In the decentralized setting, agents are only partially connected and communicate with their neighbors through a graph induced by a gossiping matrix $A=[a_{k \ell}]$, where $a_{k \ell}$ is the 
scaling weight
for the 
information flowing from 
neighboring agents $\ell \in \mathcal{N}_k$
to agent $k$. 
Our algorithm is termed \textbf{DiMA}, whose details are given in \textbf{Algorithm} \ref{DiMA}. At a high level, our algorithm consists of the local-global learning phase with three key stages:
\begin{enumerate}
\item \textbf{Local learning stage} (Steps 3--8): Agents perform multiple normalized GDA and accumulate local information through $N$ local steps.
\item \textbf{Momentum update stage} (Steps 9--10): Agents use the accumulated information to perform enhanced momentum updates.
\item \textbf{Global learning stage} (Steps 11--12): Agents track the network-average momentum, update their models using the normalized trackers, and combine the updated models with those of their neighbors.
\end{enumerate}

At the local learning stage, each agent first initializes its models using the information from the most recent communication round. Specifically, the variables $\bg^x_{k,i}$ and $\bg^y_{k,i}$ are devised to accumulate local gradient information, while $\bh^x_{k,i}$ and $\bh^y_{k,i}$ accumulate correction terms based on local samples after communication round $i$.
In step 5, each agent performs $N$ steps of normalized GDA using local step sizes $\bar{\eta}_x$ and $\bar{\eta}_y$. 
\textcolor{black}{Throughout, we use $\mathcal{N}(\bg) =\bg/\|\bg\|$  if $\bg\not= 0$ otherwise $\mathcal{N}(\bg)=0$.}
The normalization scheme ensures that local iterates do not drift significantly from their initial values, thereby allowing $\bg^x_{k,i}$ and $\bg^y_{k,i}$ to accumulate reliable gradient information.
\textcolor{black}{A small learning rate can limit local drift, but the step length also depends on the gradient norm. Agents using the same learning rate may therefore take steps of very different lengths. Normalization bounds each update direction’s norm to $1$, allowing the learning rate to directly control the step length.}
In step 6, the local stochastic gradients are incrementally accumulated into $\bg^x_{k,i}$ and $\bg^y_{k,i}$. 
\textcolor{black}{We accumulate unnormalized local gradients because the momentum vectors $\bm^x_{k,i}, \bm^y_{k,i}$ designed below estimate the unnormalized gradient direction.  Accumulating normalized gradients would not allow the variance reduction analysis because the accumulated normalized local gradient is a biased estimate of its true value.}
Since this accumulation occurs over $N$ steps, the final values correspond to the average of $N$ local stochastic gradients.
{\color{black}In steps 7--8, agent $k$ evaluates the sample Hessian at $\bz_{k,i,n}$ and applies it to $\bz_{k,i,n}-\bz_{k,i-1,n}$. Averaging these products over $n$ gives the correction $[\bh^x_{k,i};\bh^y_{k,i}]$ in the momentum update.}
\textcolor{black}{
We accumulate both stochastic gradients and HVP corrections along the local trajectory rather than at a fixed iterate. This allows each local update to move the model before subsequent derivatives are evaluated. A local ascent step may therefore improve the estimate of \(y\)-variable used in later \(x\)-gradients, motivating the use of a local trajectory instead of a mini-batch at a fixed iterate. The HVP corrections require storing the previous local trajectory, resulting in \(\mathcal{O}(NM)\) storage per agent, where \(M=M_1+M_2\). We aim to improve communication efficiency and do not claim lower overall memory usage.
}
In step 9, each agent updates its momentum using the stochastic gradients and HVP corrections accumulated during the local phase, following \eqref{main:corrected_momentum}.

\textcolor{black}{The global learning phase begins in step 11. Agents first use the mixing subroutine in Algorithm~\ref{alg:mix} to track the network-average momentum. In step 12, each agent updates its model using the normalized momentum tracker and combines the updated model with those of its neighbors through the same subroutine. As in the local phase, normalization controls the update length, which is essential for bounding the momentum correction errors in our analysis.}

\vspace{-1em}

\section{Theoretical Analysis}
\label{main:subsection:DiMA}
In this section, we establish the convergence of $\textbf{DiMA}$ towards an $\varepsilon$-GS/OS stationary point. As we focus on a fully decentralized manner, the following assumption on the communication matrix $A \in \mR^{K \times K}$ is imposed.

\begin{assumption}[\textbf{Mixing matrix}]
\label{assumption:DiMA}
\color{black}
We assume the mixing matrix $A = [a_{k\ell}]$ is symmetric, primitive, PSD, and
doubly stochastic satisfying
\begin{align}
\mathds{1}^\top_K A = \mathds{1}_K^\top, \quad A \mathds{1}_K = \mathds{1}_K,
\end{align}
\end{assumption}
\color{black}
Under this assumption, the second largest eigenvalue of $A$  satisfies
$\lambda= \|A-\mathds{1}_K\mathds{1}_K^\top/K\|<1$. Here, 
symmetry, primitivity, and double stochasticity are standard assumptions in decentralized optimization \cite{sayed2014adaptation,sayed2022inference}. 
We further use the PSD condition to bound the difference between the iterates at successive communication rounds directly, without introducing the network average as an intermediate term (see Lemma~\ref{lemma:temporal:increment}). The FastMix contraction in Lemma~\ref{lem:mix} also requires $A\succeq 0$ \cite{ye2023multi}. The PSD requirement can be easily satisfied by replacing a symmetric, primitive, and doubly stochastic matrix $A'$ with $A=(A'+I_K)/2$. 
Each agent can implement this change locally by halving its mixing weights and setting its own weight to $(1+a_{kk})/2$.

\setcounter{theorem}{0}

We now present the convergence results for the network centroid of two variables:
\begin{align}
\bx_{c,i}&=\frac1K\sum_{k=1}^{K}\bx_{k,i},
\quad 
\by_{c,i}=\frac1K\sum_{k=1}^{K}\by_{k,i},
\\
\by^o(\bx_{c,i})
&\in
\operatorname*{argmin}_{\boldsymbol v\in
Y^o(\bx_{c,i})}
\|\by_{c,i}-\boldsymbol v\|.
\end{align}
\begin{theorem}[\textbf{Convergence guarantee}]
\label{theorem:DiMA:gradientnorm:main}
\color{black}
Under Assumptions~\ref{assumption:lipchitzgradient},
\ref{assumption:boundvariance},
\ref{assumption:riskfunction}, and
\ref{assumption:DiMA}, initialize all quantities according to Algorithm \ref{DiMA}.
Choose the hyperparameters of \textbf{DiMA} as
\begin{align}
\beta=\eta_y=h,\quad
\eta_x=\frac{h}{2\kappa},\quad
\bar\eta_x=\frac{\eta_x}{NL_f},\quad
\bar\eta_y=\frac{\eta_y}{NL_f},
\end{align}
where $0<h\le1$ is a constant.
For standard gossip, set $\tau=0,r=1$.
For FastMix, choose
$\tau=1, r=\left\lceil
\frac{\log(2c_1)}{c_2\sqrt{1-\lambda}}
\right\rceil$, where $
c_1=\sqrt{14}, c_2=1-1/\sqrt{2}$.
Then, for every integer $T\ge1$, we have
\begin{align}
&\frac1T\sum_{i=0}^{T-1}
\mb{E}\Big[
\|\nabla P(\bx_{c,i})\|
+2L_f\|\by_{c,i}-\by^o(\bx_{c,i})\|
\Big] \notag\\
&\le
\frac{\kappa A_0}{hT}
+C'\kappa\left[\frac{1}{hT}+
\sqrt{\frac{h}{NK}}
+\frac{h}{1-\rho}
+\frac{\sqrt K h^2}{1-\rho}
\right]
+\frac{B_0}{T}.
\label{proof:potential_incremental44}
\end{align}
The network parameter $\rho$ satisfies $\rho=\lambda$ for standard gossip and $\rho\le 1/2$ for accelerated mixing. {\color{black}The constant $A_0$ contains the initial value objective and duality gaps.
The term $B_0/T$ includes the initial tracker disagreement and noise
with coefficient proportional to $\rho/(1-\rho)$; see
\eqref{proof:DiMA:new_initialization}. The constant $C'$ depends only
on the smoothness and noise parameters.}
See Appendix~\ref{appendix:DiMA} for the proof.
\end{theorem}

\color{black}

{\color{black}Let $I$ be uniform on $\{0,\ldots,T-1\}$ and independent of the algorithm.
If the right-hand side of Theorem~1 is at most $\varepsilon$, the centroid
pair $(\bx_{c,I},\by_{c,I})$ satisfies both stationarity criteria.}
Several terms in these bounds are independent of the network spectral gap. The factor $1-\rho$, which appears in network-related terms, is determined by the mixing subroutine: $1-\rho = 1-\lambda$ for standard gossip $(\tau=0)$ and $1-\rho \ge 1/2$ for accelerated mixing $(\tau=1)$.

The following corollaries further give the convergence rates and the complexities of reaching $\varepsilon$-stationarity for different local step counts $N$ and $\tau \in \{0,1\}$. Each outer iteration involves four mixing calls: two for momentum tracking and two for model mixing, one for each variable. Each call requires one gossip step for standard gossip and $r$ gossip steps for accelerated mixing. Counting these operations separately gives $CC=4T=\mathcal{O}(T)$ for standard gossip and $CC=4rT=\mathcal{O}(rT)$ for accelerated mixing. Each agent performs $N$ local steps per outer iteration, each requiring two stochastic gradient evaluations and two HVP evaluations. Oracle complexity (OC) counts all gradient and HVP oracle calls across $K$ agents, which gives $OC=4NKT+2K=\mathcal{O}(NKT)$ for either mixing strategy, where $2K$ accounts for momentum initialization.

{\color{black}Below, both corollaries assume $0<\varepsilon\le\kappa$ and treat the initial error terms as fixed.}
\begin{corollary}[\textbf{Standard gossip: $\tau=0$}]
\label{corollary:gossip}
Under the conditions of
Theorem~\ref{theorem:DiMA:gradientnorm:main},
we have $\rho=\lambda$. Choose
\begin{align}
h=\min\left\{
\frac{(NK)^{1/3}}{T^{2/3}},
\sqrt{\frac{1-\lambda}{T}},
\frac1{\sqrt K}
\right\}.
\label{equation:gossip:bound}
\end{align}
Then
\begin{align}
&\frac1T\sum_{i=0}^{T-1}
\mb{E}\Big[
\|\nabla P(\bx_{c,i})\|
+2L_f\delta^y_{c,i}
\Big] \label{equation:bound:gossip}\\
&\le
\mathcal{O}\left(
\frac{\kappa}{(NKT)^{1/3}}
+\frac{\kappa}{\sqrt{(1-\lambda)T}}
+\frac{\kappa\sqrt K}{T}
+\frac{\kappa}{(1-\lambda)T}
\right). \notag
\end{align}
\begin{enumerate}
\item
After transient time $T_{\rm tran} = \mathcal{O}(\frac{(NK)^2}{(1-\lambda)^3})$, the leading term
$\mathcal{O}\Big(\frac{\kappa}{(NKT)^{1/3}}\Big)$ dominates.
If $K\varepsilon\le\kappa(1-\lambda)$, choosing
$
N=\left\lceil
\frac{\kappa(1-\lambda)\varepsilon^{-1}}{K}
\right\rceil
$ gives the transient time $T_{tran}
=\mathcal{O}\left(\frac{\kappa^2\varepsilon^{-2}}{1-\lambda}\right)$. 
From the leading term, we get 
\begin{align}
CC&
=\mathcal{O}\left(
\frac{\kappa^2\varepsilon^{-2}}{1-\lambda}
\right), 
\quad 
OC=\mathcal{O}\left(
\kappa^3\varepsilon^{-3}
\right), \label{equation:CC;optimized_N:diffusion} 
\end{align}
where the $CC$ matches the order of the transient threshold, so the outer iteration includes the transient period.
\item
Choosing $N=1$ gives, for $\varepsilon\le\kappa$,
\begin{align}
CC
&=  \mathcal{O}\left(
\frac{\kappa^3\varepsilon^{-3}}{K}
+\frac{\kappa^2\varepsilon^{-2}}{1-\lambda}
+\kappa\sqrt K\,\varepsilon^{-1}
\right),\\
OC
&=\mathcal{O}\left(
\kappa^3\varepsilon^{-3}
+\frac{K\kappa^2\varepsilon^{-2}}{1-\lambda}
+\kappa K^{3/2}\varepsilon^{-1}
\right).
\end{align}
When $K\varepsilon\le\kappa(1-\lambda)$,
they reduce to
\begin{align}
CC=\mathcal{O}\left(
\frac{\kappa^3\varepsilon^{-3}}{K}
\right),
\qquad
OC=\mathcal{O}\left(
\kappa^3\varepsilon^{-3}
\right),
\end{align}
exhibiting linear communication speedup in $K$.
\end{enumerate}
\end{corollary}

\begin{corollary}[\textbf{FastMix: $\tau=1$}]
\label{corollary:fastmixing}
Under the conditions of
Theorem~\ref{theorem:DiMA:gradientnorm:main},
choose $r=\mathcal{O}\Big(\frac{1}{\sqrt{1-\lambda}}\Big)$
sufficiently large such that $\rho\le1/2$, and set
\[
h=\min\left\{
\frac{(NK)^{1/3}}{T^{2/3}},
\frac1{\sqrt T},
\frac1{\sqrt K}
\right\}.
\]
Then
\begin{align}
&\frac1T\sum_{i=0}^{T-1}
\mb{E}\Big[
\|\nabla P(\bx_{c,i})\|
+2L_f\delta^y_{c,i}
\Big] \notag\\
&\le
\mathcal{O}\left(
\frac{\kappa}{(NKT)^{1/3}}
+\frac{\kappa}{\sqrt T}
 +\frac{\kappa\sqrt K}{T}
\right).
\end{align}
\begin{enumerate}
\item
After the transient  period $T_{\rm tran} = \mathcal{O}((NK)^2)$, the leading term
$\mathcal{O}(\frac{\kappa}{(NKT)^{1/3}})$ dominates.
If $K\varepsilon\le\kappa$, choosing
$
N=\left\lceil\frac{\kappa\varepsilon^{-1}}{K}\right\rceil
$ gives an outer iteration count of $T= \mathcal{O}(\kappa^2\varepsilon^{-2})$, matching the order of the transient threshold $T_{\rm tran} = \mathcal{O}(\kappa^2\varepsilon^{-2})$. Accounting for the mixing rounds per outer iteration gives
\begin{align}
CC
&=\mathcal{O}\left(
\frac{\kappa^2\varepsilon^{-2}}{\sqrt{1-\lambda}}
\right), \quad 
OC 
=\mathcal{O}\left(
\kappa^3\varepsilon^{-3}
\right).
\end{align}

\item
Choosing $N=1$ gives
\begin{align}
CC
&= \mathcal{O}\left(
\frac1{\sqrt{1-\lambda}}
\left[
\frac{\kappa^3\varepsilon^{-3}}{K}
+\kappa^2\varepsilon^{-2}
+\kappa\sqrt K\,\varepsilon^{-1}
\right]
\right),\\
OC
&=\mathcal{O}\left(
\kappa^3\varepsilon^{-3}
+K\kappa^2\varepsilon^{-2}
+\kappa K^{3/2}\varepsilon^{-1}
\right).
\end{align}
When $K\varepsilon\le\kappa$, these simplify to
\begin{align}
CC=\mathcal{O}\left(
\frac{\kappa^3\varepsilon^{-3}}{K\sqrt{1-\lambda}}
\right),
\quad
OC=\mathcal{O}\left(
\kappa^3\varepsilon^{-3}
\right).
\end{align}
\end{enumerate}
\end{corollary}

\color{black}
We next discuss the implications of these bounds for the communication and oracle complexities of DiMA.

\textbf{Standard Gossip}.
Corollary~\ref{corollary:gossip} presents the complexity results for standard gossip. After the transient time $T_{tran} =\mathcal{O}(\frac{(NK)^2}{(1-\lambda)^3})$, the rate \eqref{equation:bound:gossip} is dominated by $\mathcal{O}(\frac{\kappa}{(NKT)^{1/3}})$. In this regime, the number of communication rounds needed to reach a target accuracy achieves linear speedup in terms of $N, K$.

When $K\varepsilon\le\kappa(1-\lambda)$, choosing $N=\lceil\frac{\kappa(1-\lambda)\varepsilon^{-1}}{K}\rceil$ balances different terms. 
This produces a communication complexity of $\mathcal{O}(\frac{\kappa^2\varepsilon^{-2}}{1-\lambda})$ and a total oracle complexity of $\mathcal{O}(\kappa^3\varepsilon^{-3})$. Increasing $N$ beyond this choice does not improve the order of the communication bound because the network-dependent term then dominates the required number of rounds.

With a single local step $N=1$, the communication complexity contains both the leading term $\mathcal{O}(\frac{\kappa^3\varepsilon^{-3}}{K})$ and the network-dependent term $\mathcal{O}(\frac{\kappa^2\varepsilon^{-2}}{1-\lambda})$. When $K\varepsilon\le\kappa(1-\lambda)$, the term  $\mathcal{O}(\frac{\kappa^3\varepsilon^{-3}}{K})$ dominates, giving a communication complexity of $\mathcal{O}(\frac{\kappa^3\varepsilon^{-3}}{K})$ and a total oracle complexity of $\mathcal{O}(\kappa^3\varepsilon^{-3})$, which gives linear communication speedup in $K$ and no spectral gap dependency. {\color{black}When $K\varepsilon>\kappa(1-\lambda)$, the network term exceeds the
first term. The term $\kappa\sqrt K/\varepsilon$ can also contribute.}

Compared with Dec-FedTrack \cite{ghiasvand2025robust}, DiMA with multiple local steps requires fewer local updates per round and improves both the communication and gradient oracle complexity bounds, thereby answering \textbf{Q1}. Among methods that use standard gossip, the single-step variant of DiMA improves the spectral-gap dependence of both bounds over DM-HSGD \cite{xian2021faster}, while the multiple-step variant improves the spectral-gap dependence of both bounds over DGDA-VR \cite{zhang2024jointly}.

\textbf{Accelerated Mixing}.
Accelerated mixing uses more gossip steps per outer iteration to reduce disagreement more rapidly. For each mixing operation,
it uses $\mathcal{O}\Big(\frac{1}{\sqrt{1-\lambda}}\Big)$ gossip steps.
For fixed $N$ and $K$, we can shorten the transient period from order $\frac{(NK)^2}{(1-\lambda)^3}$ to order $(NK)^2$.

With $N=\lceil\frac{\kappa}{K\varepsilon}\rceil$ and $K\varepsilon\le\kappa$, DiMA with accelerated mixing achieves a communication complexity of $\mathcal{O}(\frac{\kappa^2\varepsilon^{-2}}{\sqrt{1-\lambda}})$. This improves the spectral gap dependence over the standard-gossip results while retaining a total oracle complexity of $\mathcal{O}(\kappa^3\varepsilon^{-3})$.
The general communication bound for DREAM \cite{chen2023convergence} has the same dependence on $\kappa$, $\varepsilon$, and $1-\lambda$, with an additional $\log K$ factor. DREAM's refined bound removes this factor when the gradient refresh batch size of variance reduction is at least $K$. DiMA attains the same communication order using only batch size $\mathcal{O}(1)$.

With a single local step $N=1$, accelerated mixing can give a larger communication bound. When $K\varepsilon\le\kappa(1-\lambda)$, standard gossip already achieves a bound of $\mathcal{O}\left(\frac{\kappa^3\varepsilon^{-3}}{K}\right)$. In this regime, accelerated mixing does not improve the order of the outer iteration bound, so its extra gossip rounds increase the total communication bound.  Hence, accelerated mixing is useful when it saves enough outer iterations to offset the extra mixing rounds. Accelerated mixing also requires an estimate of the spectral gap to set the mixing parameters.

\textbf{Hyperparameter selection}.
In an unnormalized update, the distance an agent can move depends on both the step size and the norm of the stochastic update direction. In \textbf{DiMA}, normalization bounds that distance by the step size, which makes the control of local drift easier. For the model updates, our theorem requires
$\eta_x=\eta_y/(2\kappa)$.

Like Dec-FedTrack \cite{ghiasvand2025robust}, \textbf{DiMA} takes several local steps between communication rounds. With standard gossip, our communication bound is $\mathcal{O}\left(\frac{\kappa^2\varepsilon^{-2}}{1-\lambda}\right)$, while the bound for Dec-FedTrack is $\mathcal{O}\left(\frac{\kappa^3\varepsilon^{-2}}{(1-\lambda)^2}\right)$.
For these bounds, we choose $\mathcal{O}\left(\frac{\kappa(1-\lambda)\varepsilon^{-1}}{K}\right)$ local steps per outer iteration, compared with $\mathcal{O}\left(\varepsilon^{-2}+\frac{\kappa^2\varepsilon^{-2}}{K(1-\lambda)}\right)$ for Dec-FedTrack. Dec-FedTrack uses only stochastic gradients, while the improvement of \textbf{DiMA} comes from using additional HVPs to correct its momentum estimate. 

The convergence bound \eqref{equation:bound:gossip} explains why we do not keep increasing $N$. With standard gossip, more local steps reduce $\frac{\kappa}{(NKT)^{1/3}}$, but the network term $\frac{\kappa}{\sqrt{(1-\lambda)T}}$ is unchanged. When $K\varepsilon\le\kappa(1-\lambda)$, choosing $N=\left\lceil\frac{\kappa(1-\lambda)}{K\varepsilon}\right\rceil$ makes the two terms of the same order at the target accuracy.
Beyond this choice, the network term determines the order of the communication bound. Taking more local steps then increases the number of oracle evaluations per outer iteration without improving the order of the total communication rounds.

\vspace{-1em}
\color{black}

\color{black}

\subsection{Proof sketch}
Our proof aims to bound $\|\nabla P(\bx_{c,i})\|$ and $\delta^y_{c,i}  \triangleq \|\by_{c,i}-\by^o(\bx_{c,i})\|$. 
From Lemma \ref{lemma:DiMA},
the value function descent relation produces the desired negative gradient norm term, but it also contains a positive term proportional to $\delta^y_{c,i}$. 
From Lemma \ref{lemma:DiMA:envelope_gap}, the descent relation for $\delta^y_{c,i}$ provides a negative term to offset this effect.
To combine the two bounds, we introduce the following potential function
\color{black}
\begin{align}
&\bd{\Omega}_{i+1}\triangleq
\mb{E}
\Big[
\underbrace{P(\bx_{c,i+1})}_{\text{Value function}}+b_x\|\underbrace{\mU^x_{c,i} - \mU^x_{i}}_{\text{Tracking error}}\|+b_y\|\underbrace{\mU^y_{c,i} - \mU^y_{i}}_{\text{Tracking error}}\|\notag\\
& 
+
c_x\Big\|\underbrace{\mM^x_{i} - \mbox{\rm col}\Big\{\frac{1}{N}\sum_{n=0}^{N-1}\nabla^x_{k,i,n}\Big\}_{k=1}^{K}}_{\text{Gradient estimation error}: \, \widetilde{\mM}^x_{i}}\Big\|+ d_y\underbrace{\|\by_{c,i+1} - \by^o(\bx_{c,i+1})\|}_{y\text{-variable gap}: \, \delta^y_{c,i+1}}
\notag \\
 & + c_y
\Big\|\underbrace{\mM^y_{i} - \mbox{\rm col}\Big\{\frac{1}{N}\sum_{n=0}^{N-1}\nabla^y_{k,i,n}\Big\}_{k=1}^{K}}_{\text{Gradient estimation error}: \, \widetilde{\mM}^y_{i}} \Big\|\Big].
\end{align}
See Appendix~\ref{sec:nomenclature} for the definitions of the notation.
Taking the difference between two successive potential functions gives 
\begin{align}
&\bd{\Omega}_{i+1} - \bd{\Omega}_{i}=
\mb{E}\Big[
P(\bx_{c,i+1}) - P(\bx_{c,i})
+ b_x(\|\mU^x_{c,i} - \mU^x_{i}\|  \notag \\
&\quad - \|\mU^x_{c,{i-1}} - \mU^x_{i-1}\|)
+ b_y(\|\mU^y_{c,i} - \mU^y_{i}\|- \|\mU^y_{c,i-1} - \mU^y_{i-1}\|) \notag\\
&
\quad  + c_x(\|\widetilde{\mM}^x_{i}\| - \|\widetilde{\mM}^x_{i-1}\|)
+ c_y(\|\widetilde{\mM}^y_{i}
\|- \|\widetilde{\mM}^y_{i-1}\|) \notag\\
&\quad +d_y(\delta^y_{c,i+1} - \delta^y_{c,i})\Big]
\label{proof:sketch}
\end{align}
Lemma \ref{lemma:DiMA} gives the following descent for $P$:
\begin{align}
P(\bx_{c,i+1}) - P(\bx_{c,i})  &\le - \eta_x\|\nabla P(\bx_{c,i})\|+ 2L_f\eta_x\delta^y_{c,i}  +
\text{Others}.
\label{proof:sketch:value}
\end{align}
Lemma \ref{lemma:DiMA:envelope_gap} gives the following descent for the $y$-variable gap:
\begin{align}
&\delta^y_{c,i+1} - \delta^y_{c,i}
\le 
-\frac{1}{2}\delta^y_{c,i}+
\text{Others}. 
\label{proof:sketch:duality}
\end{align} Substituting \eqref{proof:sketch:value} and \eqref{proof:sketch:duality} into \eqref{proof:sketch} and choosing
$d_y=8L_f\eta_x$
gives the desired negative terms
\begin{align}
\mE[-2L_f\eta_x\delta^y_{c,i}-\eta_x\|\nabla P(\bx_{c,i})\|].
\label{proof:sketch:target}
\end{align} 
The terms grouped under “Others” in \eqref{proof:sketch:value} also enter the bound on $\bd{\Omega}_{i+1} - \bd{\Omega}_{i}$. Among them are the current momentum tracking errors $\mU^x_{c,i} - \mU^x_{i}, \mU^y_{c,i} - \mU^y_{i}$. Lemma~\ref{lemma:DiMA:consensus_Ux} bounds each tracking error in terms of the previous tracking error and networked momentum estimation errors $\widetilde{{\mM}}^x_{i-1}, \widetilde{{\mM}}^y_{i-1}$.  This dependence explains why both errors are involved in designing the potential function. We then substitute the momentum tracking error from Lemma~\ref{lemma:DiMA:consensus_Ux} and the networked momentum error from Lemma~\ref{lemma:DiMA:deviation:momentum_network} into \eqref{proof:sketch}, then choose $b_x,b_y,c_x,c_y$ so that the terms involving the previous tracking and estimation errors cancel.

Both recursions in Lemmas \ref{lemma:DiMA:consensus_Ux} and \ref{lemma:DiMA:deviation:momentum_network} also contain terms involving the squared model increments $\|\bx_{k,i}-\bx_{k,i-1}\|^2$ and $\|\by_{k,i}-\by_{k,i-1}\|^2$. These terms arise because the momentum correction uses the model displacement between successive communication stages. Lemma~\ref{lemma:temporal:increment} bounds them by treating the consensus and disagreement modes separately.

Finally, we sum the one-step inequality over $T$ iterations and the potential differences telescope. We use Lemma~\ref{lemma:DiMA:network_consensus} to bound model disagreement and Lemma~\ref{lemma:DiMA:deviation_ui} to bound the iteration-averaged centroid momentum error. Applying these bounds and dividing by $\eta_xT$ gives the stated convergence bound.
\color{black}

\section{Simulations}
\label{appendices:robust}
\begin{figure*}[!htbp]
\centering
\vspace{-2em}
\textbf{(a) Unconstrained formulation}\par\smallskip
\includegraphics[width=0.3\textwidth]{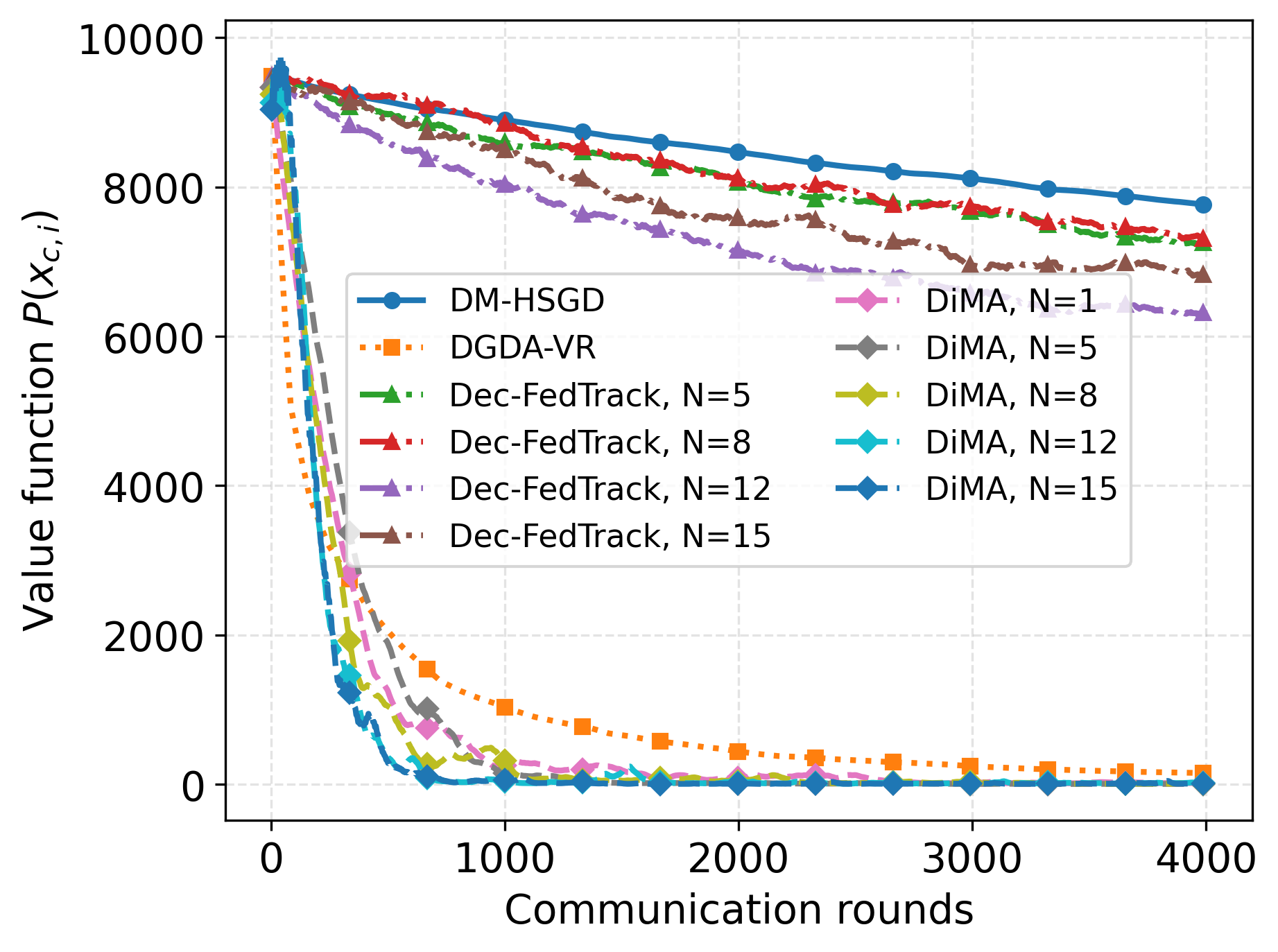}
\includegraphics[width=0.3\textwidth]{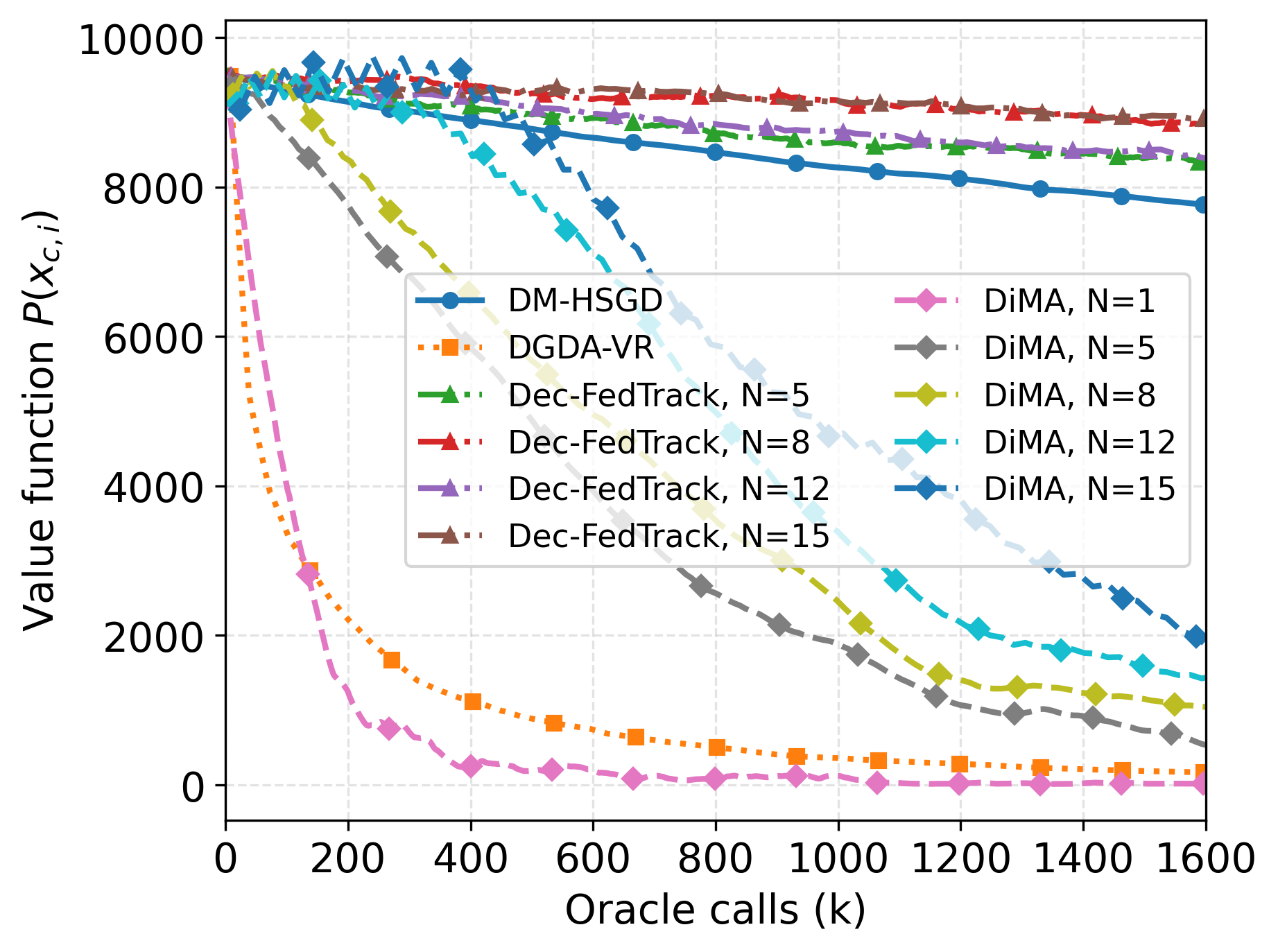}
\includegraphics[width=0.3\textwidth]{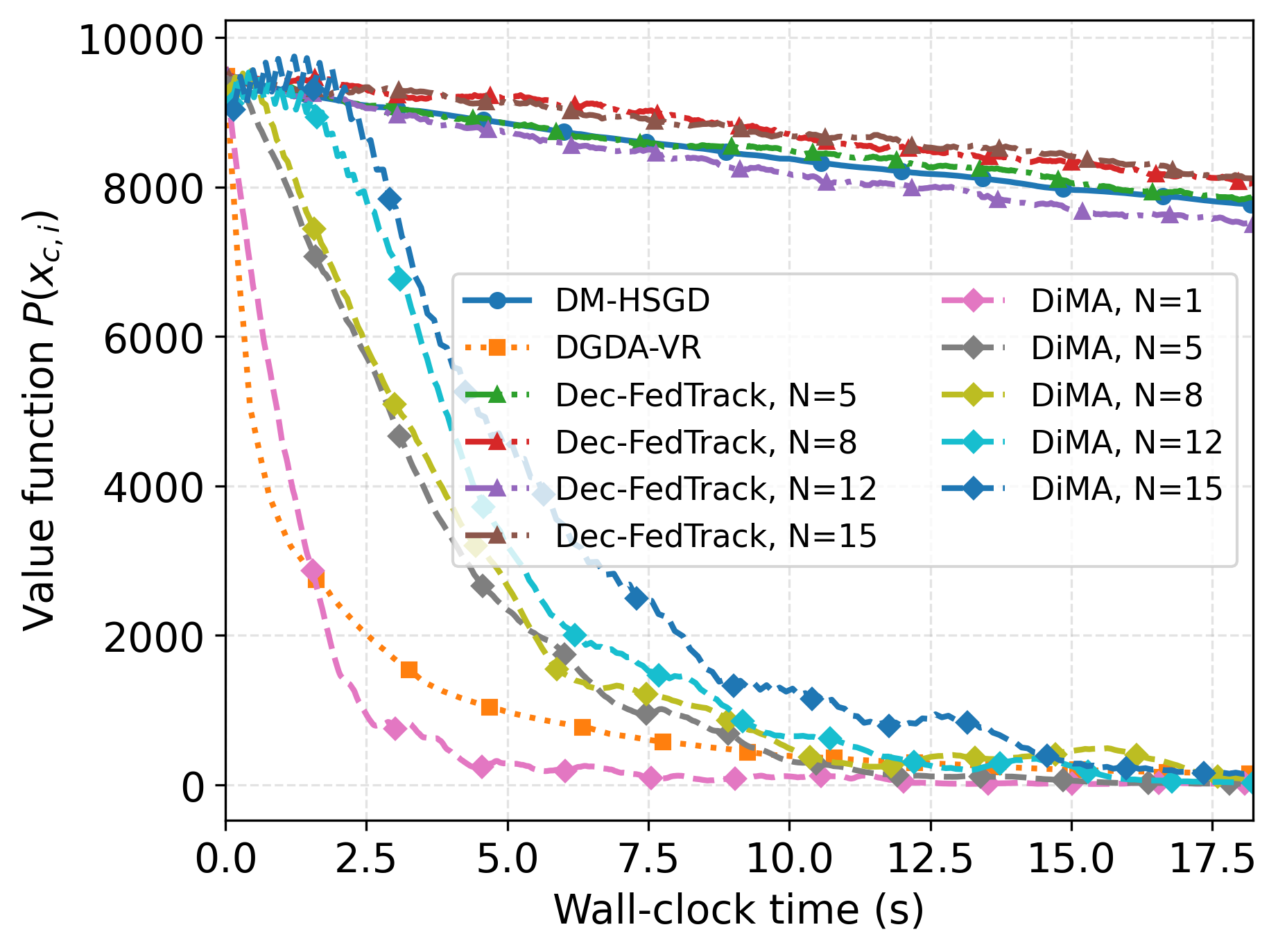}
\par
\includegraphics[width=0.3\textwidth]{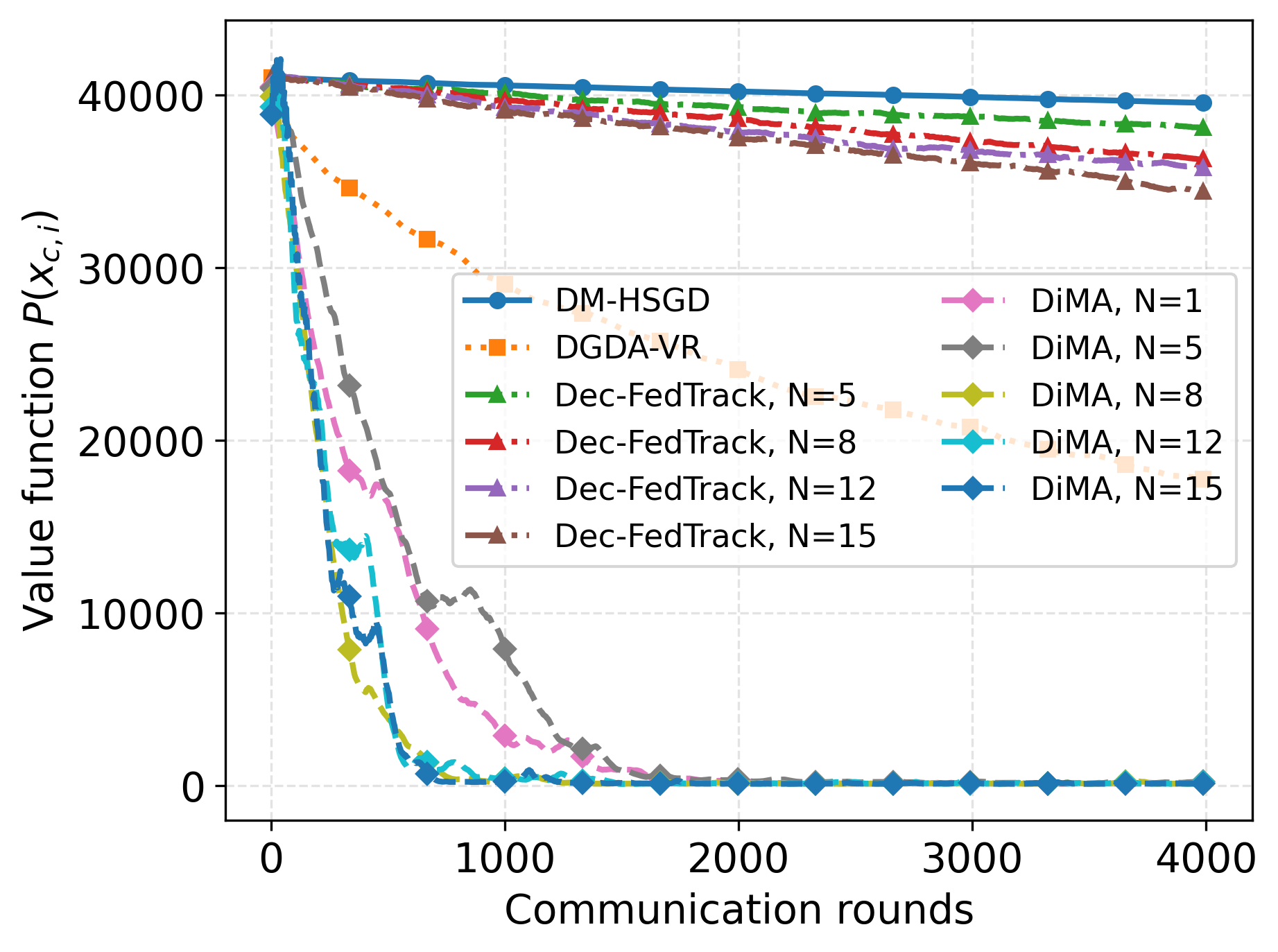}
\includegraphics[width=0.3\textwidth]{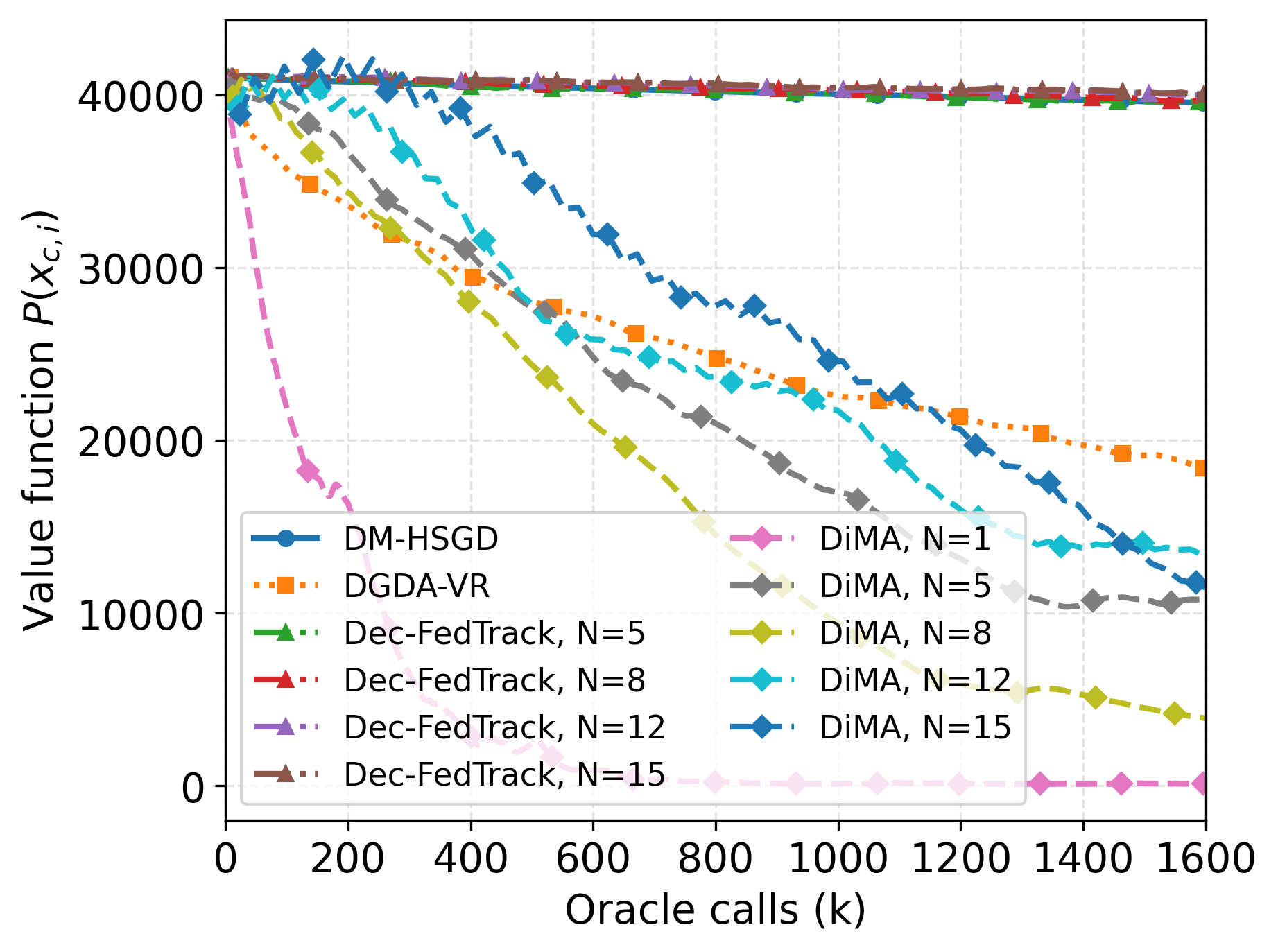}
\includegraphics[width=0.3\textwidth]{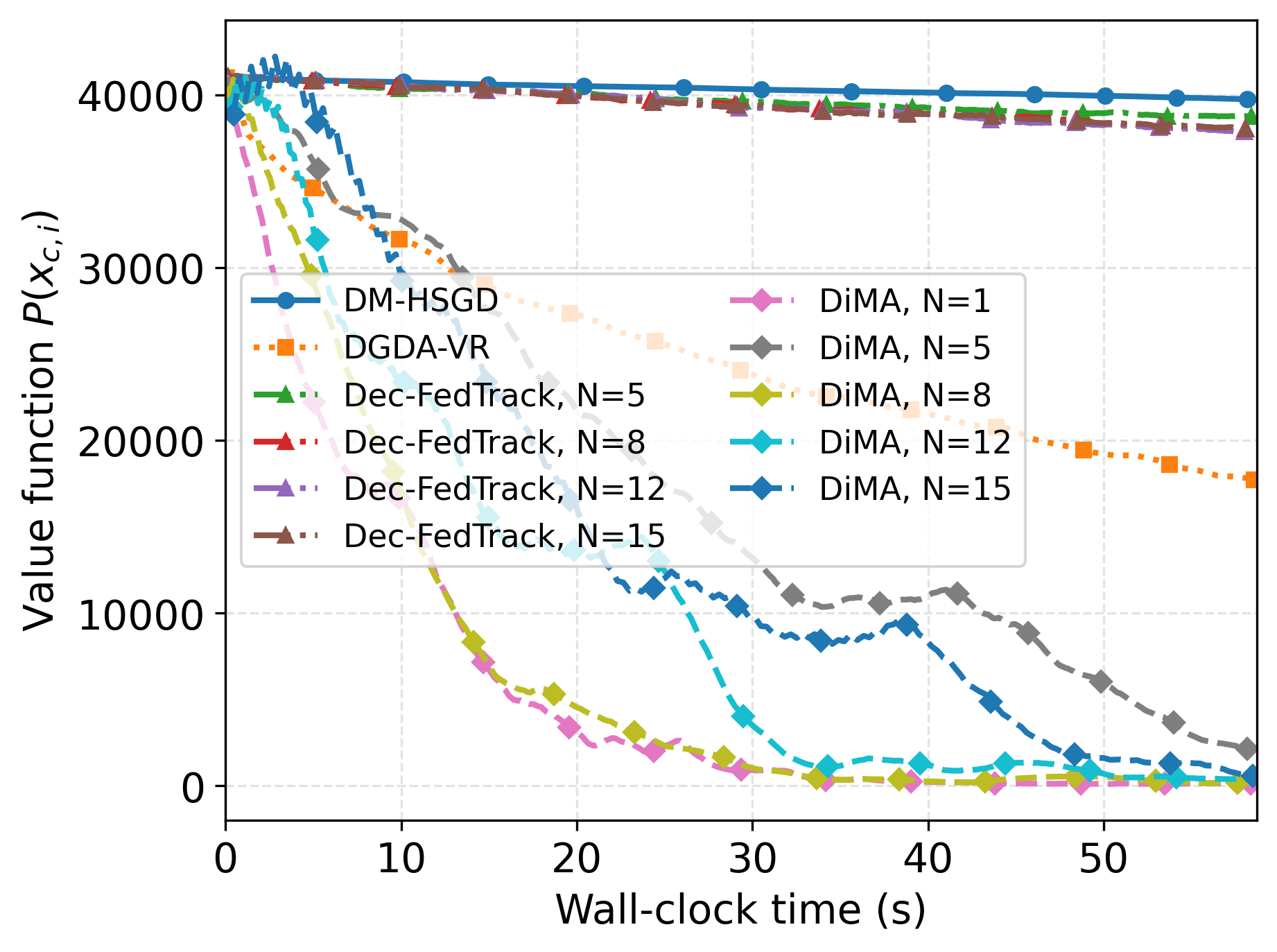}

\textbf{(b) Simplex-constrained formulation}\par\smallskip
\includegraphics[width=0.3\textwidth]{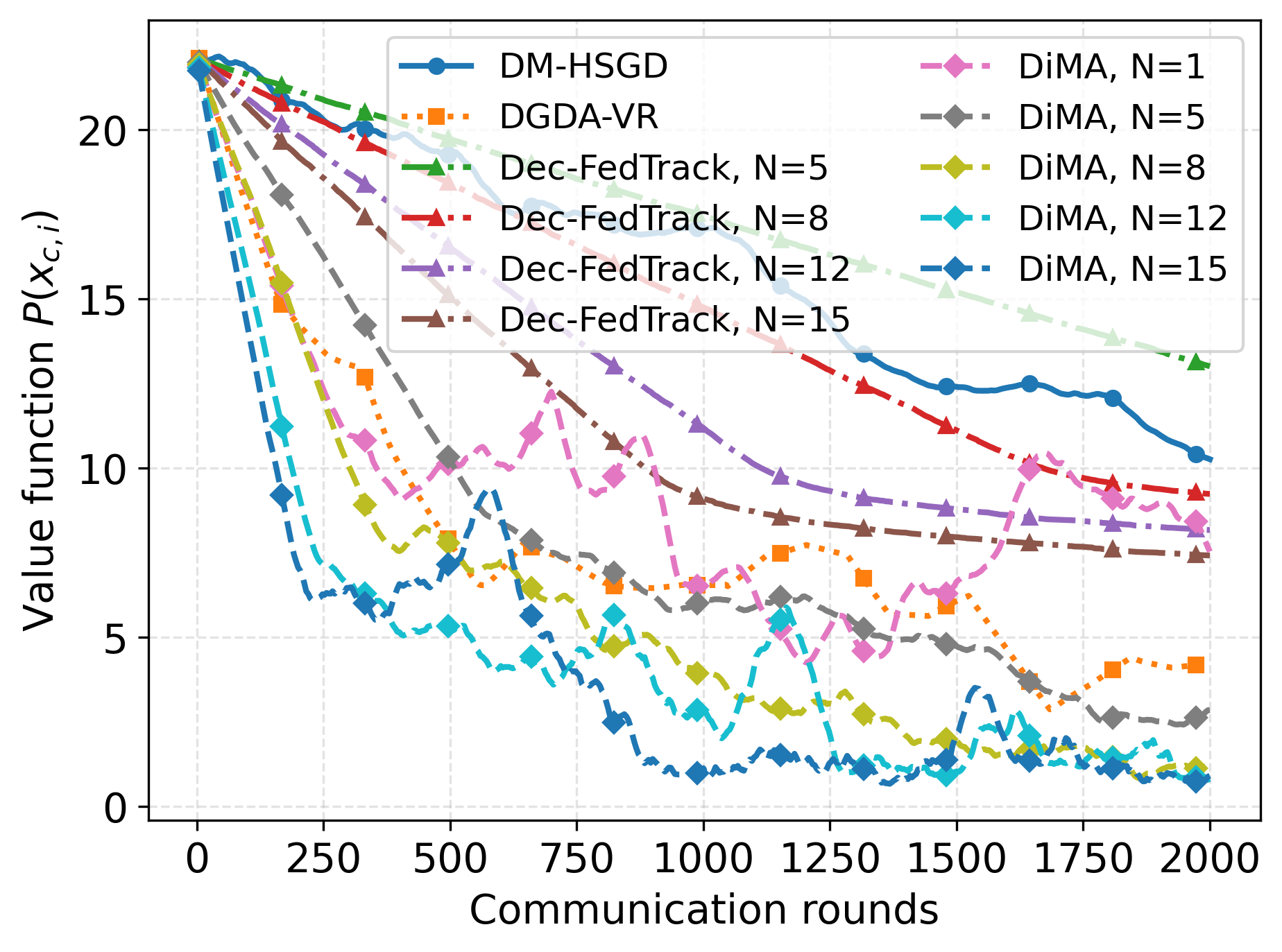}
\includegraphics[width=0.3\textwidth]{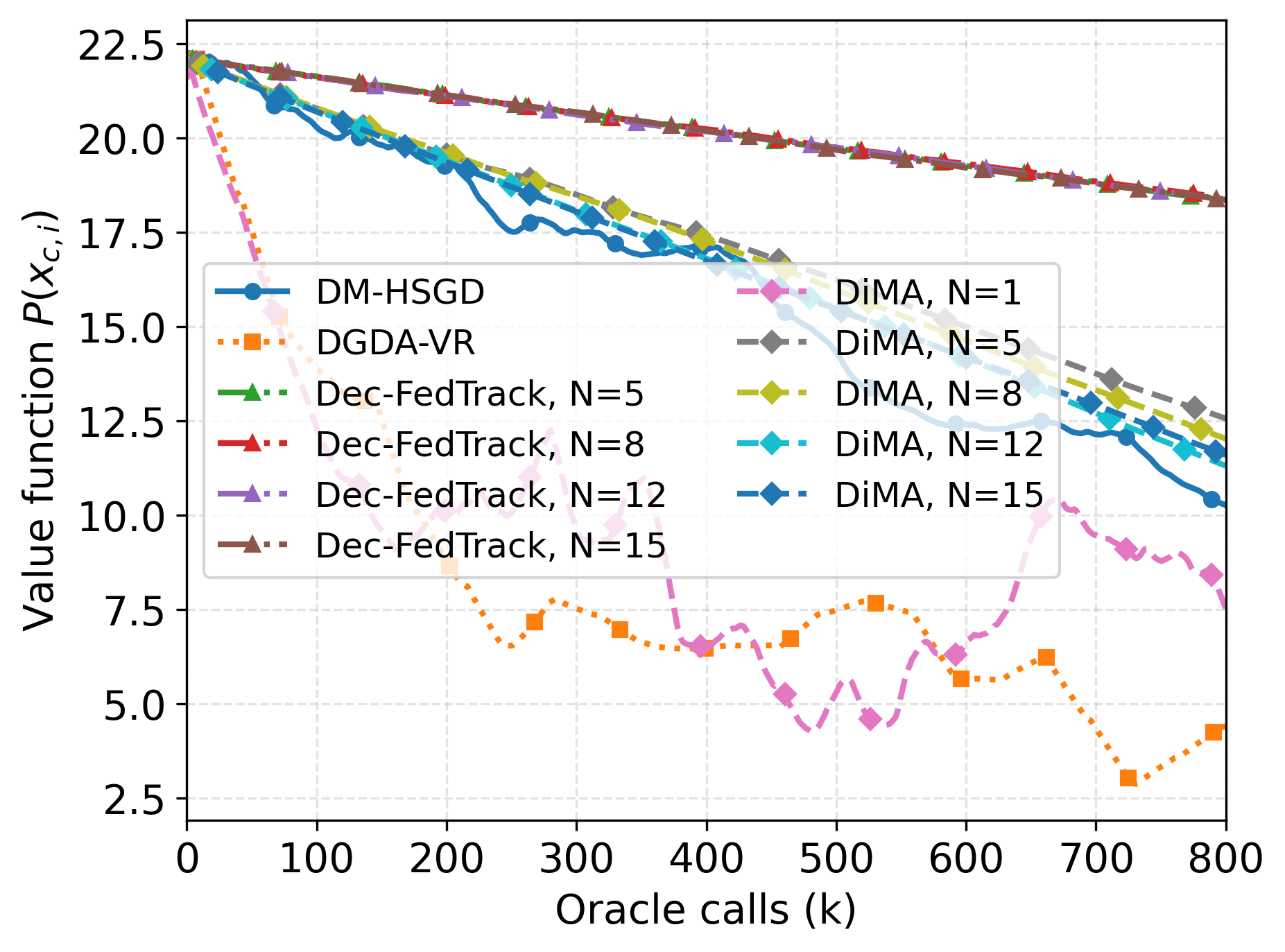}
\includegraphics[width=0.3\textwidth]{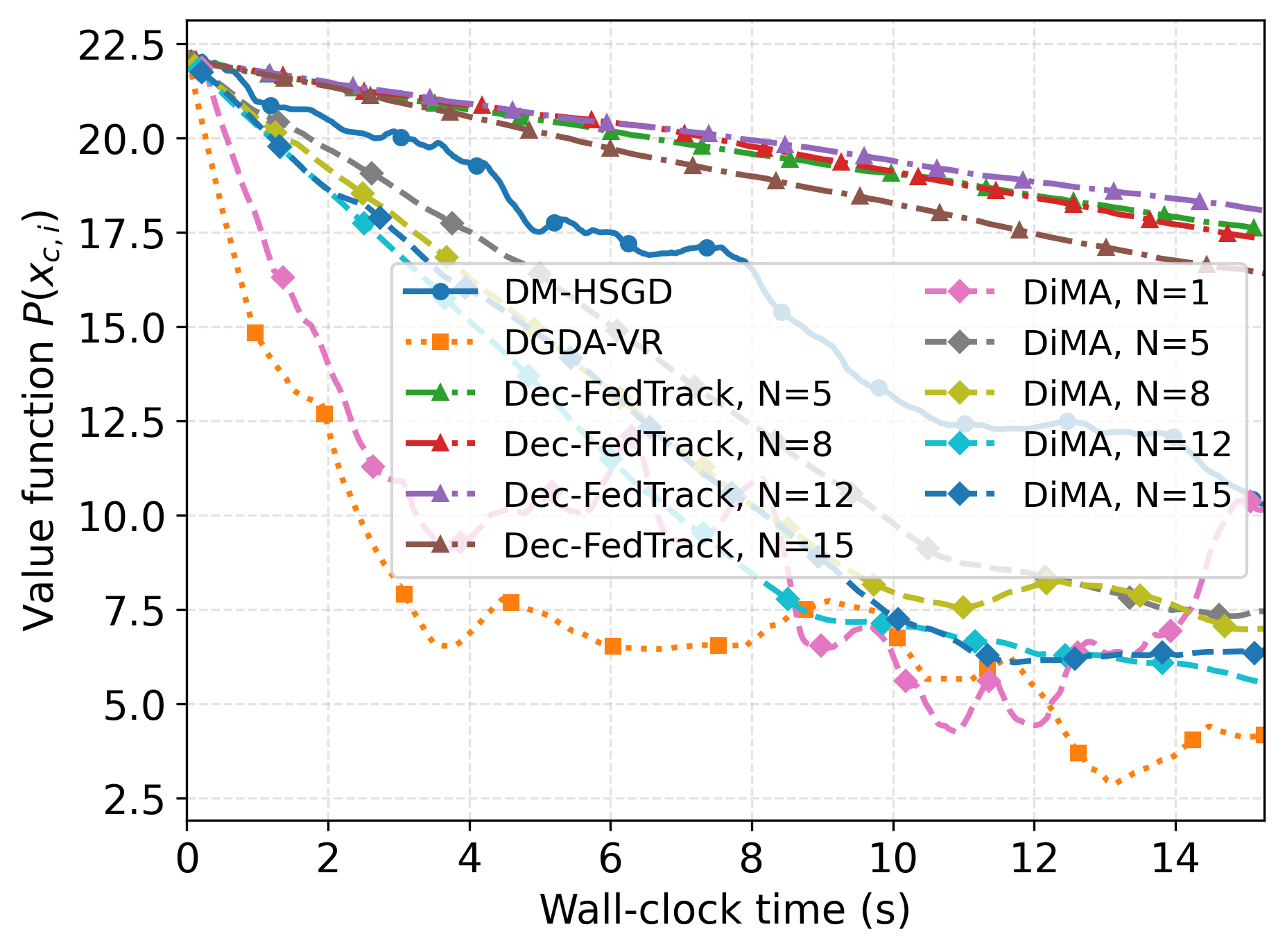}
\par
\includegraphics[width=0.3\textwidth]{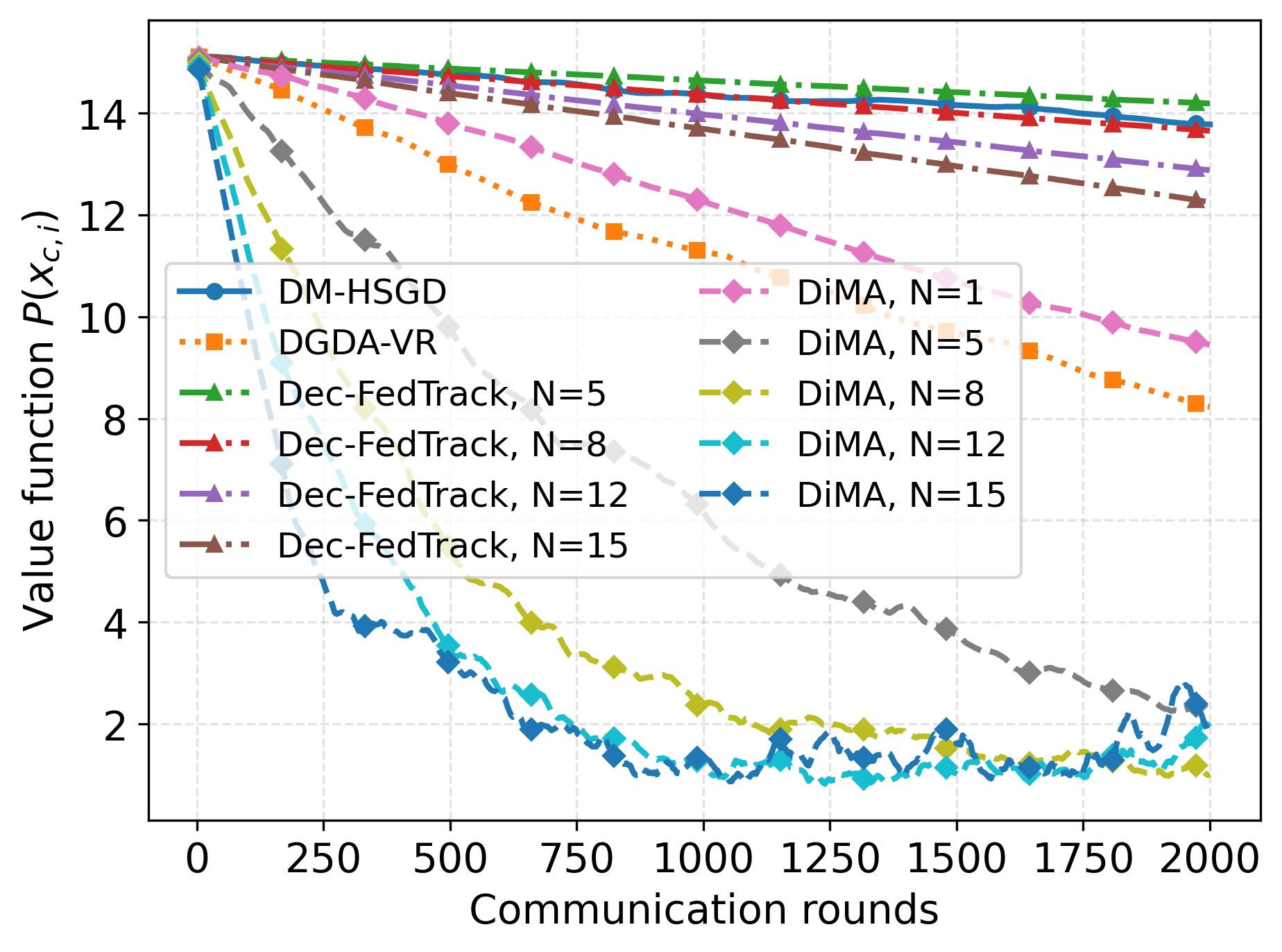}
\includegraphics[width=0.3\textwidth]{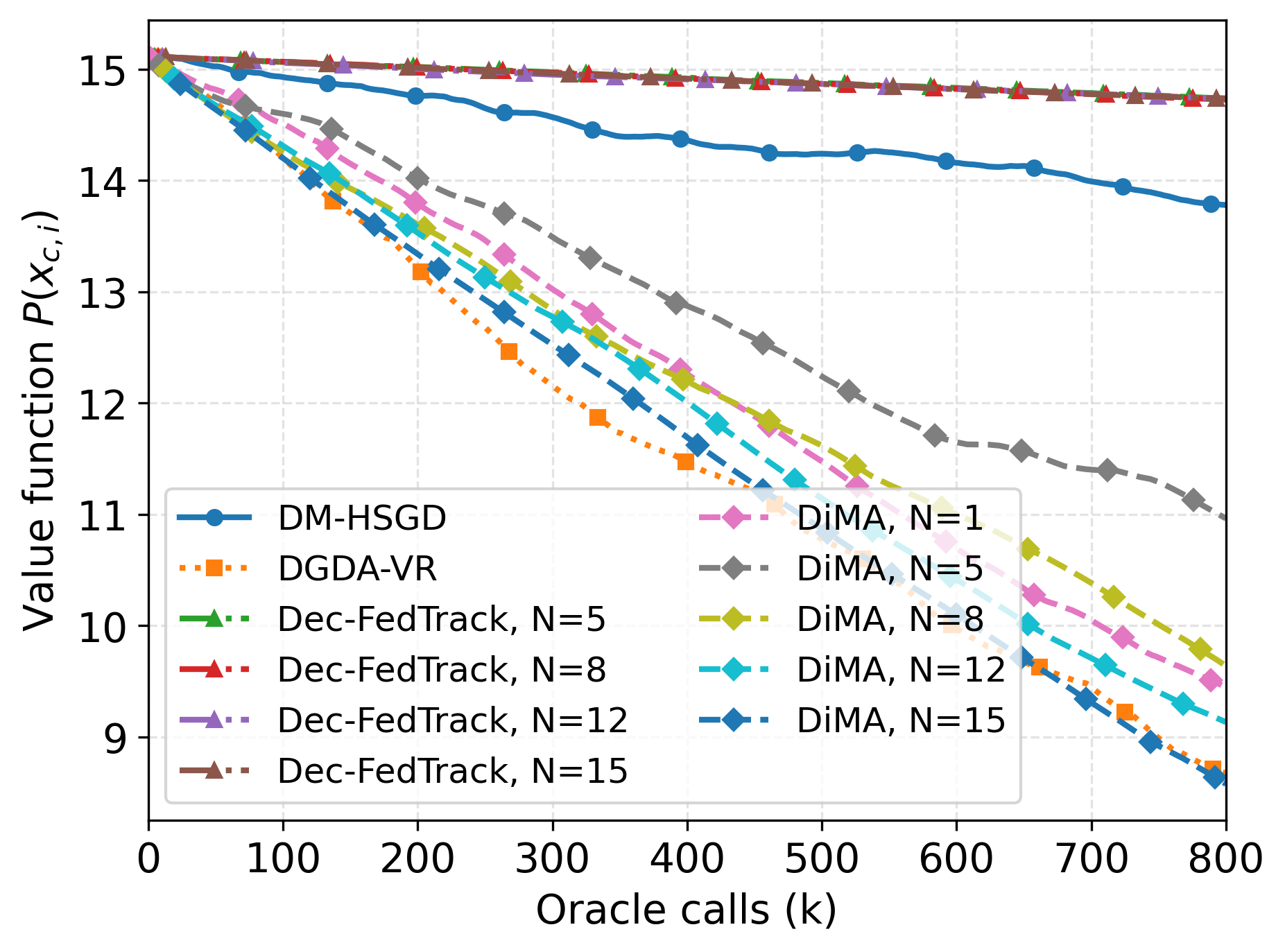}
\includegraphics[width=0.3\textwidth]{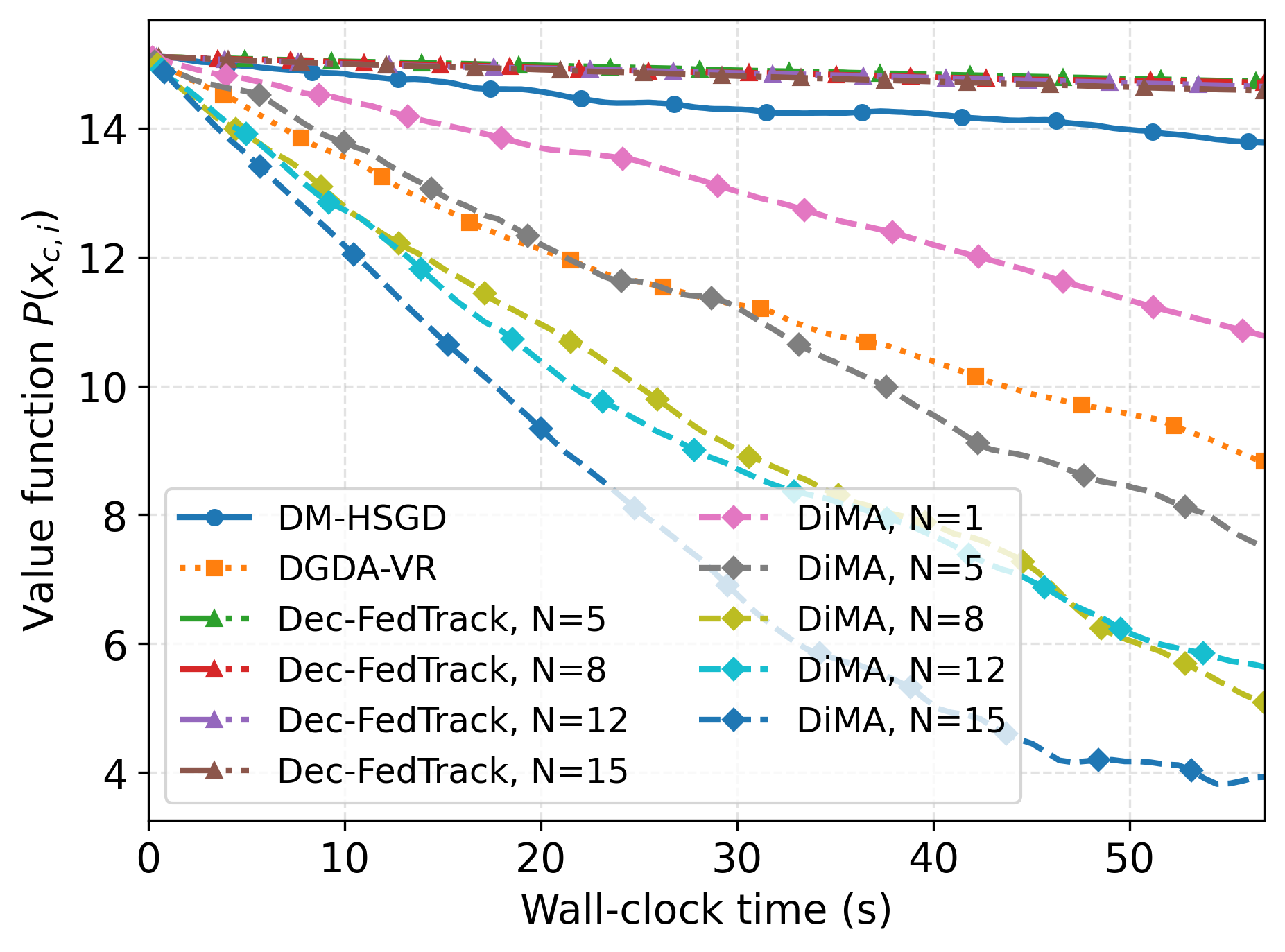}

\caption{\color{black} Comparison under (a) unconstrained and
(b) simplex-constrained formulations.
Within each formulation, the rows correspond to the
``mushrooms'' and ``a9a'' datasets, respectively.
The columns show the worst-case value $P(\bx_{c,i})$
at the network centroid versus communication rounds,
oracle calls, and wall-clock time, respectively.}
\label{figure:decentralized}
\vspace{-1.5em}
\end{figure*}

We now illustrate our algorithm \textbf{DiMA} using the robust logistic regression example considered by several works \cite{ xian2021faster, chen2024efficient,luo2020stochastic,yan2019stochastic}.
Suppose that the training data set is given by $\{(r_i, l_i)\}_{i=1}^{L}$, 
where $r_i \in \mb{R}^{d}$
is the feature vector and 
$l_i \in \{+1, -1\}$
is the binary label.
We evenly split the $L$ random samples into $K$ subsets, and each agent only has access to its local data. For convenience, we use $S_{1}, \cdots, S_{K}$
to denote the indices within these subsets, which
 correspond to a partition of the global index set $[1, \dots, L]$.
In this problem, the optimization objective is given as follows:
\begin{align}
\min _{x \in \mathbb{R}^d} \max _{y \in Y^L} J(x, y)
&=  \frac{1}{K} \sum_{k=1}^{K} J_{k}(x,y ) , \\
J_{k}(x, y)&= K\sum_{i \in S_{k}} y_i Q_i(x)-V(y)+g(x)\notag,
\end{align}
where
\begin{align}
Q_i(x) &=\log \left(1+\exp \left(-l_i r_i^T x\right)\right) \notag ,
g(x)=\rho_2 \sum_{i=1}^d \frac{\rho x_i^2}{1+\rho x_i^2} \notag ,\\
V(y)&=\frac{\rho_1}{2}\|y-\frac{\mathds{1}_L}{L}\|^2.
\end{align}
Here, $Q_i(x)$ is the loss incurred 
by the sample $(r_i, l_i)$, $g(x)$ 
is a nonconvex regularizer,
$V(y)$ is the divergence measure, 
$\mathds{1}_L$
is the $L$-dimensional one vector,
\textcolor{black}{and $Y^L =\mathbb{R}^L$ or $Y^L = \Delta_L$, where  $\Delta_L$ is a $L$-dimensional simplex defined as 
\begin{align}
\Delta_L &= \{y \in \mathbb{R}^L: 0 \leq y_i \leq 1, \sum_{i=1}^L y_i=1\}
\end{align}
We consider both the simplex-constrained formulation and an unconstrained variant with \(y\in\mathbb{R}^L\), retaining the penalty \(V(y)\) in both cases. Under the simplex constraint, \(y\) represents a probability distribution over the samples. 
Removing the simplex constraint gives the adversary more freedom to choose the sample weights. For each fixed 
$x$, this larger search space yields a worst-case regularized objective at least as large as in the constrained formulation under the same penalty parameters. For the constrained formulation, we project \(y\) onto \(\Delta_L\) after each ascent step in all compared algorithms. The projected \textbf{DiMA} is evaluated as an empirical extension of our algorithm, while the unconstrained version follows Algorithm~\ref{DiMA} directly.
We compare algorithms 
using the datasets from the
LIBSVM repository\footnote{\url{ https://www.csie.ntu.edu.tw/~cjlin/libsvmtools/datasets/}}.
}

\textbf{\color{black} Experimental setting and baselines.}
\textcolor{black}{We set $\rho_1=1.0$ for the constrained formulation and $\rho_1=100.0$, using a stronger penalty in the latter to prevent excessively large sample weights. Following \cite{xian2021faster}, we use $\rho_2=0.001$ and $\rho=10$ in both cases. All simulations are conducted on an iMac with a 4.2 GHz Intel i7 CPU and 16 GB RAM. We plot the worst-case value \(P(\bx_{c,i})=\max_y J(\bx_{c,i},y)\) at the network centroid against communication rounds, oracle calls, and wall-clock time.}
This quantity can be evaluated by plugging
$x$-variable into the objective  $J(x,y)$ and solving $P(x) = \max_y J(x,y)$ offline.

We first compare \textbf{DiMA} using standard gossip with the following baselines: the local GDA method Dec-FedTrack \cite{ghiasvand2025robust}, which has been shown to outperform DREAM \cite{chen2024efficient}; the momentum-based method DM-HSGD \cite{xian2021faster}; and the method DGDA-VR \cite{zhang2024jointly}. \textcolor{black}{We consider a ring network with \(K=20\) agents. The mixing matrix $A$  assigns weight \(2/3\) to each agent itself and \(1/6\) to each neighbor so that $\lambda \approx 0.9837$.}
For all algorithms, the mini-batch size is set to $20$.
For Dec-FedTrack and \textbf{DiMA}, we choose the number of local updates from $N \in \{5, 8, 12, 15\}$.
The local learning rates $\bar{\eta}_x, \bar{\eta}_y$ are tuned over
$\{1e^{-1}, 5e^{-2}, 1e^{-2}, 5e^{-3}, 1e^{-3}\}$, ensuring $\bar{\eta}_x \le 
\bar{\eta}_y$.
The global learning rates are set as $\eta_x = N \bar{\eta}_x, \eta_y = N\bar{\eta}_y $.
For \textbf{DiMA} and DM-HSGD, we tune the momentum parameter $\beta$
over $\{1e^{-1},1e^{-2},1e^{-3}\}$.
For DM-HSGD and DGDA-VR, we tune $\eta_x$ over 
$\{1e^{-1}, 5e^{-2}, 1e^{-2}, 5e^{-3}\}$ and $\eta_y$ over 
$\{1e^{0}, 5e^{-1},1e^{-1}, 5e^{-2}\}$.
For DGDA-VR, the large-batch computation period is set to $q= 20$, and the large-batch size is set to $S =100$. 

\color{black} 
\textbf{Experimental results.}
Figure~\ref{figure:decentralized} reports the simulation results. Under both formulations,  we can see that our algorithms can reduce the worst-case value effectively. By combining local updates with momentum acceleration, \textbf{DiMA} reduces the communication needed to reach a given objective value compared with the baselines. However, increasing the number of local steps \(N\) does not always improve oracle efficiency or wall-clock time since additional local computation and local drift can offset the communication savings.
 Moreover, \textbf{DiMA} with \(N=1\) requires less wall-clock time and fewer oracle calls in several settings.

\textbf{Additional ablation studies.}
We conduct additional ablation studies on the number of local steps \(N\), the network spectral gap, and accelerated mixing when \(\lambda\) is known to the designer. These results are provided in the supplementary material due to the page limit.

\color{black}
\vspace{-1em}
\section{Conclusion}
In this work, we proposed a communication-efficient decentralized minimax algorithm  \textbf{DiMA} by integrating accelerated momentum with local updates. To effectively combine these two schemes, we designed a novel structured local–global learning strategy for \textbf{DiMA}. 
\textcolor{black}{We established convergence guarantees for standard gossip and accelerated mixing, showing how the number of local steps affects their relative communication costs.}
Future work could explore extensions to broader classes of minimax problems, and adaptive choices of local update periods, under decentralized minimax settings.

\appendices
\section{Notation}
\label{sec:nomenclature}

\textcolor{black}{Throughout the proof, we write $\mathcal{N}(\bg)=\bg/\|\bg\|$ for $\bg \ne 0$ and set
$\mathcal{N}(0)=0$. With this convention, $\|\mathcal{N}(\bg)\|\le 1$ and
$\bg=\|\bg\|\mathcal{N}(\bg)$ for all $\bg$.}
We then use $\bz \triangleq [\bx; \by]$ to denote 
the concatenated optimization variable
 and define the true gradients.
\begin{subequations}
\begin{align}
&\nabla^z_{k,i}\triangleq[\nabla^x_{k,i};\nabla^y_{k,i}] \triangleq [\nabla_x J_k(\bz_{k,i});\nabla_y J_k(\bz_{k,i})] \in \mR^{M_1+M_2}, \\
&\nabla^z_{k,i,n}\triangleq
[\nabla^x_{k,i,n};\nabla^y_{k,i,n}]\notag\\
&\qquad   \triangleq [\nabla_x J_k(\bz_{k,i,n});\nabla_y J_k(\bz_{k,i,n})] \in \mR^{M_1+M_2}, 
\end{align}
where $k,i$, and $n$ denote the agent,  communication round, and local step index, respectively.
The superscripts $x,y$, and $z$ indicate the gradient associated with $x, y$ or $z= [x;y]$, respectively.

Let \(\bd{\mathcal{F}}_{i,n}\) denote the filtration before the agents draw their samples at local step \(n\) of round \(i\). Then, both \(\bz_{k,i,n}\) and \(\bz_{k,i-1,n}\) are \(\mathcal F_{i,n}\)-measurable. We can define
\begin{align}
&[h^x_{k,i,n}; h^y_{k,i,n}] \triangleq \mE[\bh^x_{k,i,n}; \bh^y_{k,i,n} \mid \bd{\mathcal{F}}_{i,n}] \in \mR^{M_1+M_2}, \\
&[h^x_{k,i}; h^y_{k,i}] \triangleq
\frac{1}{N}\sum_{n=0}^{N-1}[h^x_{k,i,n}; h^y_{k,i,n}] \in \mR^{M_1+M_2},\\
&[h^x_{c,i}; h^y_{c,i}] \triangleq \frac{1}{NK}\sum_{k=1}^{K}\sum_{n=0}^{N-1} [h^x_{k,i,n};h^y_{k,i,n}]\in \mR^{M_1+M_2}.
\end{align}
\end{subequations}
For the network-averaged quantities, i.e., the averaged quantities over $k \in [K]$, we denote
\begin{subequations}
\begin{align}
&[\bx_{c,i};\bg^x_{c,i};\bm^x_{c,i};\bh^x_{c,i};\bu^x_{c,i}] \notag \\
&\triangleq \frac{1}{K}\sum_{k=1}^{K}[\bx_{k,i};\bg^x_{k,i};\bm^x_{k,i};\bh^x_{k,i};\bu^x_{k,i}] \in \mR^{5M_1},  \\
&[\by_{c,i};\bg^y_{c,i};\bm^y_{c,i};\bh^y_{c,i};\bu^y_{c,i}] \notag \\
&\triangleq \frac{1}{K}\sum_{k=1}^{K}[\by_{k,i};\bg^y_{k,i};\bm^y_{k,i};\bh^y_{k,i};\bu^y_{k,i}] \in \mR^{5M_2}, \\
&\bz_{c,i} \triangleq [\bx_{c,i};\by_{c,i}] \triangleq  \frac{1}{K}\sum_{k=1}^{K}[\bx_{k,i};\by_{k,i}] \in \mR^{M_1+M_2}. 
\end{align}
\end{subequations}
Here, 
the subscript ``$c$" means centroid.

We denote \emph{network-extended} quantities, i.e.,  concatenated variables over $k \in [K]$, as:
\begin{subequations}
\begin{align}
[\mX_i;\mY_i ] &\triangleq [\mbox{\rm col}\{\bx_{k,i}\}_{k=1}^{K} ; \mbox{\rm col}\{\by_{k,i}\}_{k=1}^{K}]\in\mb{R}^{KM} ,
  \\
[\mG^x_{i};\mG^y_{i}] &\triangleq [\mbox{\rm col}
\{\bg^x_{k,i}\}_{k=1}^{K} ;\mbox{\rm col}\{\bg^y_{k,i}\}_{k=1}^{K}]\in\mb{R}^{KM},   \\
[\mM^x_{i}; \mM^y_{i}]&\triangleq [\mbox{\rm col} \{\bm^x_{k,i}\}_{k=1}^{K}; \mbox{\rm col} \{\bm^y_{k,i}\}^{K}_{k=1}]\in\mb{R}^{KM}, \\
[\mH^x_{i}; \mH^y_{i}] &\triangleq [\mbox{\rm col}
\{\bh^x_{k,i}\}_{k=1}^{K}; \mbox{\rm col}
\{\bh^y_{k,i}\}_{k=1}^{K}] \in\mb{R}^{KM},  \\
[\mU^x_{i};\mU^y_{i}] &\triangleq [\mbox{\rm col}
\{\bu^x_{k,i}\}_{k=1}^{K};\mbox{\rm col}
\{\bu^y_{k,i}\}_{k=1}^{K}] \in \mb{R}^{KM}.
\end{align}
\end{subequations}
The {\em network-extended averaged}  quantities
are defined as:
\begin{subequations}
\begin{align}
[{\mX}_{c,i}; {\mY}_{c,i}] &\triangleq 
[\mathds{1}_K \otimes {\bx}_{c,i}; \mathds{1}_K \otimes {\by}_{c,i}]
\in\mb{R}^{KM},  \\
[\mG^x_{c,i};\mG^y_{c,i}] &\triangleq 
[\mathds{1}_K \otimes\bg^x_{c,i};\mathds{1}_K \otimes\bg^y_{c,i}]
\in\mb{R}^{KM} ,
 \\
[\mM^x_{c, i};\mM^y_{c, i} ] &\triangleq 
[\mathds{1}_K \otimes\bm^x_{c,i}; \mathds{1}_K \otimes \bm^y_{c,i}]
\in\mb{R}^{KM} ,
\\
[\mH^x_{c,i};\mH^y_{c,i}] &\triangleq [\mathds{1}_K \otimes \bh^x_{c,i}; \mathds{1}_K \otimes\bh^y_{c,i}]
\in\mb{R}^{KM} ,
\\
[\mU^x_{c,i};\mU^y_{c,i}] &\triangleq 
[\mathds{1}_K \otimes\bu^x_{c,i};\mathds{1}_K \otimes\bu^y_{c,i}] 
\in\mb{R}^{KM} . 
\end{align}
\end{subequations}
For  \emph{network-extended normalized}   quantities and the associated averaged quantities, we denote 

{\small
\begin{subequations}
\begin{align}
[\overline{\mU}^x_{i}; \overline{\mU}^y_{i}]&\triangleq \Big[ \mbox{\rm col}\Big\{\mathcal{N}(\bu^x_{k,i})\Big\}_{k=1}^{K};
\mbox{\rm col}\Big\{ \mathcal{N}(\bu^y_{k,i})\Big\}_{k=1}^{K}\Big]\in\mb{R}^{KM} ,
\\
&[\overline{\mU}^x_{c,i};\overline{\mU}^y_{c,i}] \notag\\
\notag&\triangleq \Big[ \mathds{1}_K \otimes  \frac{1}{K} \sum_{r=1}^K\mathcal{N}(\bu^x_{r,i}); \mathds{1}_K \otimes\frac{1}{K} \sum_{r=1}^{K}\mathcal{N}(\bu^y_{r,i})
\Big]\in\mb{R}^{KM}. 
\end{align}
\end{subequations}
}
We further use regular font symbols $ \mnH^x_{i}, \mnH^y_{i}$ to denote the \emph{deterministic} realization of 
$ \mH^x_{i}, \mH^y_{i}$, i.e.,
\begin{subequations}
\begin{align}
\mnH^x_i
&\triangleq \operatorname{col}\Big\{
    \frac{1}{N}\sum_{n=0}^{N-1} h^x_{k,i,n}
    \Big\}_{k=1}^{K}
    \in \mR^{KM_1}, \\
\mnH^y_i
&\triangleq \operatorname{col}\Big\{
    \frac{1}{N}\sum_{n=0}^{N-1} h^y_{k,i,n}
    \Big\}_{k=1}^{K}
    \in \mR^{KM_2}.
\end{align}
\end{subequations}
For some other types of error terms,
we denote 
\begin{subequations}
\begin{align}
 \mathcal{E}_{i} &\triangleq \|\mX_{i}-\mX_{c,i}\|
+\|\mY_{i} - \mY_{c,i}\| \in \mR, \\
\Delta^y_{c, i} &\triangleq P(\bx_{c,i}) - J(\bx_{c,i},\by_{c,i}) \in \mR, \\
\delta^y_{c,i}
&\triangleq
\operatorname{dist}\!\left(
\by_{c,i},Y^o(\bx_{c,i})
\right)
\triangleq
\|\by_{c,i}-\by^o(\bx_{c,i})\| \in \mR.
\end{align}
\end{subequations}
where $\by^o(\bx_{c,i}) \in Y^o(\bx_{c,i})$ is the closest maximizer to $\by_{c,i}$.
Additionally, 
we introduce 
\begin{subequations}
\begin{align}
\mc{J} \triangleq& 
\frac{\mathds{1}_K \mathds{1}^\top_{K}}{K} \otimes I_{M'}, \quad 
\mc{A} \triangleq  A \otimes I_{M'} \in \mR^{KM'\times KM'},
\end{align}
\end{subequations}
where $M' = M_1 \text{ or } M_2$ is compatible to the dimension of $x$ or $y$.

\vspace{-1em}
\section{Basic Lemmas}
 We now provide some basic lemmas.
\begin{Lemma}[\textbf{Danskin-type lemma}
\cite{nouiehed2019solving}]
\label{lemma:Danskin2}
Under Assumptions~\ref{assumption:lipchitzgradient}
and~\ref{assumption:riskfunction}, define the nonempty
maximizer set
$
Y^o(x)=\operatorname*{argmax}_y J(x,y).
$
\textcolor{black}{
The value function \(P(x)=\max_y J(x,y)\) is
\(L\)-smooth with
$
L=L_f+L_f\kappa, \kappa=\frac{L_f}{\nu},
$
and satisfies
$
\nabla P(x)=\nabla_x J(x,y^o), \forall y^o\in Y^o(x).
$
Moreover,
\begin{align}
P(x)-J(x,y)
\geq\frac{\nu}{2}
\operatorname{dist}\bigl(y,Y^o(x)\bigr)^2,
\end{align}
where \(\operatorname{dist}(y,S)\triangleq\inf_{v\in S}\|y-v\|\) denotes the distance from \(y\) to the set \(S\).
}
\end{Lemma}

\begin{Lemma}[Lemma 5 \cite{cai2024accelerated}]
\label{lemma:meanvalue}
Let $J(z)$ be twice continuously differentiable with
$L_h$-Lipschitz Hessian. Then, for all $z_1,z_2$,
it holds that 
\begin{align}
&\Big\|
\nabla_z J (z_2)
- \nabla_z J (z_1)+ \nabla^2_{z}J(z_1)
(z_1-z_2) 
\Big\|
\le 
\frac{L_h}{2}\|z_1 -z_2\|^2.
\end{align}
\end{Lemma}

\color{black}
For our mixing schemes, we provide the following results.
\begin{Lemma}[\textbf{Mixing subroutine}]
\label{lem:mix}
Let Assumption~\ref{assumption:DiMA} hold.
Let
$
\mV^s=\operatorname{col}\{\bd v_k^{(s)}\}_{k=1}^K
$
denote the network quantities after the $s$-th mixing step,
with $\mV^{-1}=\mV^0=\operatorname{col}\{\bd v_k\}_{k=1}^K$. Denote the final output by
\begin{align}
\mV^R=\operatorname{Mix}(\mV;A,r,\tau).
\end{align}
where $R = r$ (FastMix) or $R=1$ (standard gossip).
Both mixing subroutines satisfy
\begin{align}
\mathcal{J}\mV^R&=\mathcal{J}\mV^0, \label{eq:subroutine:average}\\
\|\mV^R-\mathcal{J}\mV^0\|
&\leq\rho\|\mV^0-\mathcal{J}\mV^0\|,\label{eq:subroutine:contraction}
\end{align}
where 
\begin{align}
\rho=
\begin{cases}
\lambda, & \tau=0 \quad \text{(Diffusion)},\\
c_1(1-c_2\sqrt{1-\lambda})^r,
& \tau=1 \quad \text{(FastMix)},
\end{cases}
\label{eq:subroutine:rho}
\end{align}
with \(c_1=\sqrt{14}\) and \(c_2=1-1/\sqrt{2}\).
For FastMix, choosing
\begin{align}
r\geq
\left\lceil
\frac{\log(2c_1)}{c_2\sqrt{1-\lambda}}
\right\rceil
\label{Lemma:r:choice}
\end{align}
ensures \(\rho\leq1/2\) \cite{ye2023multi}.
\end{Lemma}
\begin{proof}
Since \(A\) is doubly stochastic,
$
\mathcal J(A\otimes I_{M'})=\mathcal J
$.
For standard diffusion, \(\boldsymbol v^R=(A\otimes I)\boldsymbol v^0\), the relation \eqref{eq:subroutine:average} holds directly. For FastMix,
we have $\forall s \ge 0$
\begin{align}
\mathcal{J}\mV^{s+1}
=(1+\eta_a)\mathcal{J}\mV^s
-\eta_a\mathcal{J}\mV^{s-1}.
\end{align}
Since \(\mV^{-1}=\mV^0\), an induction proves
\eqref{eq:subroutine:average}.

For Diffusion, \(R=1\) and \(\eta_a=0\), so
\begin{align}
\mV^R-\mathcal{J}\mV^0
=
\left(A\otimes I_{M'}-\mathcal{J}\right)(I-\mathcal{J})\mV^0.
\end{align}
Using
\(\|A-\mathds{1}_K\mathds{1}_K^{\top}/K\|=\lambda\)
gives \eqref{eq:subroutine:contraction}. For FastMix, the relations \eqref{eq:subroutine:contraction} and \eqref{eq:subroutine:rho} follow
from
\cite{ye2023multi}. Finally, \(1-u\leq e^{-u}\) implies
\begin{align}
\rho= c_1(1-c_2\sqrt{1-\lambda})^r \le c_1e^{-c_2r\sqrt{1-\lambda}}\leq\frac12,
\end{align}
where the last inequality follows from \eqref{Lemma:r:choice}.
\end{proof}
We next bound the network model increment. To this end,
 recall the model updating rule:
\begin{align}
\mX_{i+1}
&=\operatorname{Mix}(\mX_i+\bd{v}_i^x;A,r,\tau), \label{update:rule:mixX}\\
\mY_{i+1}
&=\operatorname{Mix}(\mY_i+\bd{v}_i^y;A,r,\tau),
\end{align}
where
\begin{align}
\bd{v}_i^x=-\eta_x{\rm col}
\Big\{\mathcal N(\bu^x_{k,i})\Big\}_{k=1}^K,
\bd{v}_i^y=\eta_y {\rm col}
\Big\{\mathcal N(\bu^y_{k,i})\Big\}_{k=1}^K,\notag
\end{align}
and $\mathcal N(v)=v/\|v\|$. We state the following results.

\begin{Lemma}[\textbf{Model increment}]
\label{lemma:temporal:increment}
Under Assumption~\ref{assumption:DiMA}, suppose that
$\bz_{k,0}=\bz_{k,-1}=\bz_0$ and
$\bz_{k,-1,n}=\bz_0$ for every agent $k$ and
$n=0,\ldots,N-1$.
We choose $\eta_x\le\eta_y$,
$\bar{\eta}_x\le\bar{\eta}_y$, and
use a fixed mixing subroutine for Algorithm \ref{DiMA}. If using FastMix, we further choose $r$ so that $\rho\le1/2$.
For $T\ge1$, we have
\begin{align}
\sum_{i=0}^{T-1}\sum_{k=1}^K
\|\bz_{k,i+1}-\bz_{k,i}\|^2
&\le2C_\tau KT\eta_y^2,
\label{eq:DiMA:temporal_increment}\\
\sum_{i=0}^{T-1}\sum_{k=1}^K
\|\bz_{k,i}-\bz_{k,i-1}\|^2
&\le2C_\tau KT\eta_y^2,
\label{eq:DiMA:temporal_increment_shifted}
\end{align}
where $C_0=1$ for standard gossip and $C_1=4$ for
FastMix.
Moreover, defining
\begin{align}
D_T\triangleq
\frac1{KNT}
\sum_{i=0}^{T-1}\sum_{k=1}^K\sum_{n=0}^{N-1}
\mb{E}\|\bz_{k,i,n}-\bz_{k,i-1,n}\|^2,
\label{equation:definition:DT}
\end{align}
we have
\begin{equation}
D_T\le6C_\tau\eta_y^2+4N^2\bar{\eta}_y^2.
\label{proof:DiMA:local_temporal_difference0}
\end{equation}
\end{Lemma}

\begin{proof}
We first bound the temporal increments of $\mX_i$.
Since $A$ is symmetric, we obtain the following eigendecomposition
\begin{align}
A=Q\Lambda Q^\top,
\quad
Q=[q_1,\ldots,q_K],
\quad
\Lambda=\operatorname{diag}\{\lambda_1,\ldots,\lambda_K\},
\end{align}
where $Q$ is orthogonal and $Aq_j=\lambda_jq_j$.
Assumption~\ref{assumption:DiMA} implies
\[
q_1=\frac{\mathds{1}_K}{\sqrt K},
\qquad
\lambda_1=1,
\qquad
0\le\lambda_j<1,\quad j\ge2.
\]
Define the transformed model and direction vectors by
\[
\widehat{\bx}_{j,i}
=(q_j^\top\otimes I_{M_1})\mX_i,
\qquad
\widehat{\bd{v}}^x_{j,i}
=(q_j^\top\otimes I_{M_1})\bd{v}_i^x.
\]
Identical initial models give
$\mX_0=\mathds{1}_K\otimes\bx_0$, and hence
\begin{align}
\widehat{\bx}_{j,0}
=(q_j^\top\mathds{1}_K)\bx_0=0,\qquad j\ge2.
\label{equation:initial:disagreemode}
\end{align}
For standard gossip, the mixing operator is $A$.
For FastMix, we define the following polynomial
\begin{align}
p_{-1}(t)&=p_0(t)=1,\label{equation:polynomial}\\
p_{s+1}(t)
&=(1+\eta_a)t\,p_s(t)-\eta_a p_{s-1}(t),
\quad s=0,\ldots,r-1, \notag
\end{align}
where $\eta_a$ is the fixed acceleration parameter.
After $r$ inner iterations, $\operatorname{Mix}(\mX;A,r,\tau)$ maps an input
$\mX\in\mathbb R^{KM_1}$ to
\begin{align}
(p_r(A)\otimes I_{M_1})\mX.
\end{align}
Since $A=Q\Lambda Q^\top$,
we can write
\begin{align}
p_r(A)=Qp_r(\Lambda)Q^\top.
\end{align}
Using the eigencomponents,
we can transform the updating rule \eqref{update:rule:mixX} for both mixing subroutines into each mode $j\ge 1$:
\begin{align}
\widehat{\bx}_{j,i+1}
=a_j\bigl(
\widehat{\bx}_{j,i}+\widehat{\bd{v}}^x_{j,i}
\bigr),
\quad
a_j=
\begin{cases}
\lambda_j, & \tau=0,\\
p_r(\lambda_j), & \tau=1.
\end{cases}
\label{equation:mode:recursion}
\end{align}
Below, we discuss two cases separately.

\emph{Case 1: Standard gossip ($\tau=0$).}
For $\lambda_j>0$, set
$c_j=(1-\lambda_j)/\lambda_j\ge0$.
Expanding the squared norms gives
\begin{align}
&\|\widehat{\bx}_{j,i+1}-\widehat{\bx}_{j,i}\|^2
+c_j\bigl(
\|\widehat{\bx}_{j,i+1}\|^2
-\|\widehat{\bx}_{j,i}\|^2
\bigr) \notag\\
&=
\lambda_j\|\widehat{\bd{v}}^x_{j,i}\|^2
-\frac{(1-\lambda_j)^2}{\lambda_j}
\|\widehat{\bx}_{j,i}\|^2
\le\|\widehat{\bd{v}}^x_{j,i}\|^2.
\end{align}
Summing over $i=0,\ldots,T-1$ yields
\begin{align}
&\sum_{i=0}^{T-1}
\|\widehat{\bx}_{j,i+1}-\widehat{\bx}_{j,i}\|^2
+c_j\|\widehat{\bx}_{j,T}\|^2\le
\sum_{i=0}^{T-1}\|\widehat{\bd{v}}^x_{j,i}\|^2
+c_j\|\widehat{\bx}_{j,0}\|^2.
\end{align}
The initial term vanishes because $c_1=0$ for the
consensus mode and $\widehat{\bx}_{j,0}=0$ for
$j\ge2$; see \eqref{equation:initial:disagreemode}. Dropping the nonnegative term $c_j\|\widehat{\bx}_{j,T}\|^2$ gives
\begin{align}
\sum_{i=0}^{T-1}
\|\widehat{\bx}_{j,i+1}-\widehat{\bx}_{j,i}\|^2
\le
\sum_{i=0}^{T-1}\|\widehat{\bd{v}}^x_{j,i}\|^2.
\label{bound:case:gossip}
\end{align}
If $\lambda_j=0$, the corresponding disagreement mode, i.e., 
$\widehat{\bx}_{j,i+1}$ with $j\ge 2$,
starts at zero and remains zero. The consensus mode
$\lambda_1=1$ satisfies \eqref{bound:case:gossip} with equality.

\emph{Case 2: FastMix ($\tau=1$).}
The polynomial recursion \eqref{equation:polynomial} gives $p_r(1)=1$, so $a_1=1$.
The consensus mode increment therefore satisfies
\begin{align}
\widehat{\bx}_{1,i+1}-\widehat{\bx}_{1,i}
=\widehat{\bd{v}}^x_{1,i}.
\label{consensus:fastmix:mode}
\end{align}
For $j\ge2$, let
$J_K=\mathds{1}_K\mathds{1}_K^\top/K$.
Since $J_Kq_j=0$, $p_r(A)q_j=a_jq_j$, and $\|q_j\|=1$,
we have
\begin{align}
|a_j|
=\|(p_r(A)-J_K)q_j\|
\le\rho\|q_j-J_Kq_j\|
=\rho\le\frac12,
\end{align}
where the first inequality follows
by applying \eqref{eq:subroutine:contraction} with \(M'=1\) and initial input \(q_j\).
Recursively applying the mode recursion \eqref{equation:mode:recursion} and using
$\widehat{\bx}_{j,0}=0$ gives
\begin{align}
\widehat{\bx}_{j,i+1}
=\sum_{\ell=0}^{i}
a_j^{\,i-\ell+1}\widehat{\bd{v}}^x_{j,\ell}.
\end{align}
For $a_j\ne0$, we can derive
\begin{align}
\|\widehat{\bx}_{j,i+1}\|^2
&\le
\left(
\sum_{\ell=0}^{i}|a_j|^{i-\ell+1}
\|\widehat{\bd{v}}^x_{j,\ell}\|
\right)^2 \notag\\
&\le
\left(\sum_{\ell=0}^{i}|a_j|^{i-\ell+1}\right)
\left(\sum_{\ell=0}^{i}|a_j|^{i-\ell+1}
\|\widehat{\bd{v}}^x_{j,\ell}\|^2\right) \notag\\
&\le
\frac{|a_j|}{1-|a_j|}
\sum_{\ell=0}^{i}|a_j|^{i-\ell+1}
\|\widehat{\bd{v}}^x_{j,\ell}\|^2,
\end{align}
where the first inequality follows from 
triangle inequality and the second inequality follows from the Cauchy-Schwarz inequality.
Summing over $i$ and collecting the terms for each
fixed $\ell$ gives
\begin{align}
\sum_{i=0}^{T-1}\|\widehat{\bx}_{j,i+1}\|^2
&\le
\frac{|a_j|}{1-|a_j|}
\sum_{\ell=0}^{T-1}
\|\widehat{\bd{v}}^x_{j,\ell}\|^2
\sum_{i=\ell}^{T-1}|a_j|^{i-\ell+1} \notag\\
&\le
\left(\frac{|a_j|}{1-|a_j|}\right)^2
\sum_{\ell=0}^{T-1}
\|\widehat{\bd{v}}^x_{j,\ell}\|^2 \le
\sum_{\ell=0}^{T-1}
\|\widehat{\bd{v}}^x_{j,\ell}\|^2,
\end{align}
The last inequality uses $|a_j|\le1/2$.
For $a_j=0$, the same bound holds because the model
mode remains zero. Since $\widehat{\bx}_{j,0}=0$,
\begin{align}
&\sum_{i=0}^{T-1}
\|\widehat{\bx}_{j,i+1}-\widehat{\bx}_{j,i}\|^2\le
2\sum_{i=0}^{T-1}
\bigl(
\|\widehat{\bx}_{j,i+1}\|^2
+\|\widehat{\bx}_{j,i}\|^2
\bigr) \notag\\
&\le
4\sum_{i=0}^{T-1}\|\widehat{\bx}_{j,i+1}\|^2
\le
4\sum_{i=0}^{T-1}\|\widehat{\bd{v}}^x_{j,i}\|^2.
\end{align}
Considering \eqref{consensus:fastmix:mode},
 the above inequality holds $\forall j\ge 1$.
 
For either mixing subroutine, summing over $j$ and using
the orthogonal invariance under the Euclidean norm  with $Q^\top\otimes I_{M_1}$  gives

\begin{align}
\sum_{i=0}^{T-1}\|\mX_{i+1}-\mX_i\|^2
\le C_\tau\sum_{i=0}^{T-1}\|\bd{v}_i^x\|^2
\le C_\tau KT\eta_x^2.
\end{align}
where $C_0=1$ for standard gossip and $C_1=4$ for
FastMix.
The same argument for $\mY_i$ gives
\begin{align}
\sum_{i=0}^{T-1}\|\mY_{i+1}-\mY_i\|^2
\le C_\tau KT\eta_y^2.
\end{align}
Adding these bounds and using $\eta_x\le\eta_y$
proves \eqref{eq:DiMA:temporal_increment}.

Since $\bz_{k,-1}=\bz_{k,0}$,
\begin{align}
\sum_{i=0}^{T-1}\sum_{k=1}^K
\|\bz_{k,i}-\bz_{k,i-1}\|^2
&=
\sum_{i=0}^{T-2}\sum_{k=1}^K
\|\bz_{k,i+1}-\bz_{k,i}\|^2 \notag\\
&\le2C_\tau KT\eta_y^2,
\end{align}
which proves \eqref{eq:DiMA:temporal_increment_shifted}.
Moreover, we have
\begin{align}
\bz_{k,i,n}-\bz_{k,i-1,n}
&=(\bz_{k,i,n}-\bz_{k,i})
+(\bz_{k,i}-\bz_{k,i-1}) \notag\\
&\quad+(\bz_{k,i-1}-\bz_{k,i-1,n}).
\end{align}
The normalized local updates ensure
\begin{align}
\|\bz_{k,i,n}-\bz_{k,i}\|^2
\le n^2(\bar{\eta}_x^2+\bar{\eta}_y^2).
\end{align}
The same bound holds for $\bz_{k,i-1}-\bz_{k,i-1,n}$,
including $i=0$ because of identical initialization.
Therefore,
\begin{align}
&\|\bz_{k,i,n}-\bz_{k,i-1,n}\|^2 \le
3\|\bz_{k,i}-\bz_{k,i-1}\|^2
+6n^2(\bar{\eta}_x^2+\bar{\eta}_y^2). \notag
\end{align}
Taking expectations, averaging over $i,k,n$, and
using \eqref{eq:DiMA:temporal_increment_shifted} gives
\begin{align}
&D_T
\le
6C_\tau\eta_y^2
+\frac6N\sum_{n=0}^{N-1}
n^2(\bar{\eta}_x^2+\bar{\eta}_y^2) \\
&=
6C_\tau\eta_y^2
+(N-1)(2N-1)(\bar{\eta}_x^2+\bar{\eta}_y^2)
\le
6C_\tau\eta_y^2+4N^2\bar{\eta}_y^2, \notag
\end{align}
which proves
\eqref{proof:DiMA:local_temporal_difference0}.
\end{proof}

\color{black}
\begin{Lemma}[\textbf{Descent relation}]
\label{lemma:DiMA}
Under Assumptions~\ref{assumption:lipchitzgradient}, \ref{assumption:boundvariance}, \ref{assumption:riskfunction}, and~\ref{assumption:DiMA}, suppose $0<\bar\eta_x\le\bar\eta_y$. Then
\begin{align}
&P(\bx_{c,i+1}) \le P(\bx_{c,i}) - \eta_x\|\nabla P(\bx_{c,i})\|+ 2L_f\eta_x
\delta^y_{c,i}+ \notag\\ 
& \frac{2\eta_xL_f\mathcal{E}_i}{\sqrt{K}}+ \frac{L\eta^2_x}{2} 
 +
2\eta_x
\Big\|\frac{1}{NK}\sum_{k=1}^{K}\sum_{n=0}^{N-1}
\nabla^x_{k,i,n}
- \bu^x_{c,i}\Big\|
 \notag\\
 &\quad 
 +
\frac{\eta_x}{\sqrt{K}}\|\mU^x_{c,i} - \mU^x_{i}\|+2NL_f\eta_x \bar{\eta}_y. 
\end{align}
\end{Lemma}
\textcolor{black}{
The above lemma establishes the error terms that need to be controlled to ensure descent of the value function.}

\begin{Lemma}[\textbf{Consensus error}]
\label{lemma:DiMA:network_consensus}
With the initialization in Algorithm~\ref{DiMA}, $\eta_x\le\eta_y$,
and a fixed mixing operator satisfying \eqref{eq:subroutine:contraction}
with $\rho<1$, the model disagreement satisfies, for all $i\ge-1$,
\begin{align}
\mathcal{E}_{i}\le \frac{4\sqrt{K}\eta_y}{1-\rho}, 
\end{align}
where $\rho$ is a constant; see Lemma \ref{lem:mix}.
\end{Lemma}
\textcolor{black}{
The above lemma provides a unified consensus error bound for \textbf{DiMA} with local updates under either standard diffusion or accelerated mixing.
Below, we bound the error of averaged momentum $\bm^w_{c,i}$ using $\bu^w_{c,i} = \bm^w_{c,i}, \forall w\in\{x, y\}$. We introduce 
\begin{align}
 \widetilde{\bu}^w_{c,i}
&\triangleq
\bu^w_{c,i}
-\frac1{KN}\sum_{k=1}^{K}\sum_{n=0}^{N-1}
\nabla^w_{k,i,n}   
\end{align}
}

\begin{Lemma}[\textbf{Averaged momentum error}]
\label{lemma:DiMA:deviation_ui}
Under Assumptions~\ref{assumption:lipchitzgradient}--\ref{assumption:DiMA}, use the initialization of Algorithm~\ref{DiMA} and the mixing conditions of Lemma~\ref{lemma:temporal:increment}. If $\eta_x\le\eta_y$, $\bar\eta_x\le\bar\eta_y$, and $0<\beta\le1$, then for each integer $T\ge1$,
\color{black}
\begin{align}
&\frac{1}{T}\sum_{i=0}^{T-1}\mb{E}\Big\| \widetilde{\bu}^w_{c,i}\Big\|  \le   Q_1,
\end{align}
where $w = x \text{ or y}$ and
\begin{align}
Q_1 &\triangleq
\frac{\sigma}{\sqrt K\beta T}
+\frac{L_h(3C_{\tau}\eta_y^2
+2N^2\bar{\eta}_y^2)}{\beta} \notag\\
&\quad 
+
\frac{\sigma_h
\sqrt{6C_{\tau}\eta_y^2+4N^2\bar{\eta}_y^2}}
{\sqrt{NK\beta}}
+\frac{\sigma\sqrt{\beta}}{\sqrt{NK}},
\end{align}
where $C_{\tau} =1$ for standard gossip and $C_{\tau} =4$  for FastMix.
\end{Lemma}

\textcolor{black}{Below, we derive the momentum tracker and update error. To this end, we define 
\begin{align}
\widetilde{\mM}^w_{i} 
&\triangleq\mM^w_{i} - \mbox{\rm col}\Big\{\frac{1}{N}\sum_{n=0}^{N-1}\nabla^w_{k,i,n}\Big\}_{k=1}^{K},
\end{align}
where $w = x \text{ or } y$}.
\begin{Lemma}[\textbf{Momentum tracking error}]
\label{lemma:DiMA:consensus_Ux}
Under Assumptions \ref{assumption:lipchitzgradient}, \ref{assumption:boundvariance}, \ref{assumption:riskfunction}, \ref{assumption:DiMA},   by choosing parameters $\eta_x \le \eta_y, \bar{\eta}_x \le \bar{\eta}_y, 0<\beta \le 1$, we can bound the momentum tracking error as follows: \color{black}
\begin{align}
&\mb{E}\|\mU^w_{c,i}  - \mU^w_{i}\| \le
\rho\mb{E}\|\mU^w_{i-1}-\mU^w_{c,i-1}\|
+\rho\beta\mb{E}\|\widetilde{\mM}^w_{i-1}\|
\notag\\
&\quad+
\frac{\rho KL_h}{2}d_i^2
+\rho\left(
L_f\sqrt K
+\sigma_h\sqrt{\frac KN}
\right)d_i +\rho\beta\sigma\sqrt{\frac KN}. \notag
\end{align}
where $w = x \text{ or } y$,
and 
\begin{align}
d_i^2
&\triangleq
\frac{1}{KN}
\sum_{k=1}^{K}\sum_{n=0}^{N-1}
\mb{E}\|\bz_{k,i,n}-\bz_{k,i-1,n}\|^2,
\quad
d_i\triangleq\sqrt{d_i^2}.
\notag
\end{align}
\end{Lemma}

\begin{Lemma}[\textbf{Networked momentum error}]
\label{lemma:DiMA:deviation:momentum_network}
Under Assumptions~\ref{assumption:lipchitzgradient},
\ref{assumption:boundvariance},
\ref{assumption:riskfunction}, and
\ref{assumption:DiMA}, choose
$\eta_x\le\eta_y$, $\bar{\eta}_x\le\bar{\eta}_y$,
and $0<\beta\le1$. \color{black}
Then, for $w\in\{x,y\}$,
\begin{align}
\mb{E}\|\widetilde{\mM}^w_i\|
&\le
(1-\beta)\mb{E}\|\widetilde{\mM}^w_{i-1}\|
+\frac{KL_h}{2}d_i^2 +
\sigma_h\sqrt{\frac KN}\,d_i
\notag\\
&\quad 
+\beta\sigma\sqrt{\frac KN}.
\end{align}
\end{Lemma}

\begin{Lemma}[\textbf{Duality gap}]
\label{lemma:DiMA:envelope_gap}
Under Assumptions~\ref{assumption:lipchitzgradient}, \ref{assumption:boundvariance}, \ref{assumption:riskfunction}, and~\ref{assumption:DiMA}, suppose $0<\eta_x\le\eta_y/(2\kappa)$ and $0<\bar\eta_x\le\bar\eta_y$. Then
\begin{align}
&\delta^y_{c,i+1}
\le 
\frac{1}{2}\delta^y_{c,i}
+ \frac{(\Delta^y_{c,i} - \Delta^y_{c,i+1})}{\eta_y \nu}+
\frac{2\kappa}{\sqrt{K}}\mathcal{E}_i
+2\kappa N \bar{\eta}_y+5\kappa \eta_y
\notag \\
&\quad+
\frac{2}{\nu}\Big\|\bu^y_{c,i} -\frac{1}{KN}
\sum_{k=1}^{K}\sum_{n=0}^{N-1} \nabla^y_{k,i,n}\Big\| +\frac{1}{\sqrt{K}\nu}
\|\mU^y_i - \mU^y_{c,i}\|. 
\end{align}
\end{Lemma}
\vspace{-2em}
\section{Proof of Theorem \ref{theorem:DiMA:gradientnorm:main}}
\label{appendix:DiMA}

\color{black}
\begin{proof}
For $w\in\{x,y\}$, recall the notation
\begin{align}
\widetilde{\bu}^w_{c,i}
&\triangleq
\bu^w_{c,i}
-\frac1{KN}\sum_{k=1}^{K}\sum_{n=0}^{N-1}
\nabla^w_{k,i,n},\\
\widetilde{\mM}^w_i
&\triangleq
\mM^w_i-
\operatorname{col}\Big\{
\frac1N\sum_{n=0}^{N-1}\nabla^w_{k,i,n}
\Big\}_{k=1}^{K}.
\end{align}
We use the following potential function
\begin{align}
\bd{\Omega}_{i+1}&
\triangleq\mb{E}\Big[
P(\bx_{c,i+1})
+b_x\|\mU^x_i-\mU^x_{c,i}\|
+b_y\|\mU^y_i-\mU^y_{c,i}\| \notag\\
&+c_x\|\widetilde{\mM}^x_i\|
+c_y\|\widetilde{\mM}^y_i\|
+d_y\delta^y_{c,i+1}
\Big].
\label{proof:DiMA:new_potential}
\end{align}
Note that step sizes and coefficients are constant
over the $T$ iterations.
Combining Lemmas~\ref{lemma:DiMA}
and~\ref{lemma:DiMA:envelope_gap}, and choosing
$d_y=8L_f\eta_x$, gives
\begin{align}
&\bd{\Omega}_{i+1}-\bd{\Omega}_i \notag\\
&\le
-\eta_x\mb{E}\Big[
\|\nabla P(\bx_{c,i})\|
+2L_f\delta^y_{c,i}
\Big]
+\frac{(2L_f+16\kappa L_f)\eta_x}{\sqrt K}
\mb{E}\mathcal E_i \notag\\
&\quad+
2\eta_x\mb{E}\|\widetilde{\bu}^x_{c,i}\|
+16\kappa\eta_x
\mb{E}\|\widetilde{\bu}^y_{c,i}\|
+\frac{L\eta_x^2}{2}
+2NL_f\eta_x\bar\eta_y \notag\\
&\quad+
\left(b_x+\frac{\eta_x}{\sqrt K}\right)
\mb{E}\|\mU^x_i-\mU^x_{c,i}\|
-b_x\mb{E}\|\mU^x_{i-1}-\mU^x_{c,i-1}\| \notag\\
&\quad+
\left(b_y+\frac{8\kappa\eta_x}{\sqrt K}\right)
\mb{E}\|\mU^y_i-\mU^y_{c,i}\|
-b_y\mb{E}\|\mU^y_{i-1}-\mU^y_{c,i-1}\| \notag\\
&\quad+
c_x\left(
\mb{E}\|\widetilde{\mM}^x_i\|
-\mb{E}\|\widetilde{\mM}^x_{i-1}\|
\right)
+c_y\left(
\mb{E}\|\widetilde{\mM}^y_i\|
-\mb{E}\|\widetilde{\mM}^y_{i-1}\|
\right) \notag\\
&\quad+
\frac{8\kappa\eta_x}{\eta_y}
\mb{E}\left[
\Delta^y_{c,i}-\Delta^y_{c,i+1}
\right]
+8L_f\eta_x
\left(5\kappa\eta_y+2\kappa N\bar\eta_y\right).
\label{proof:DiMA:new_potential_increment}
\end{align}
We next apply
Lemmas~\ref{lemma:DiMA:consensus_Ux}
and~\ref{lemma:DiMA:deviation:momentum_network}.
Choose
\begin{align}
b_x=c_x
&=\frac{\rho\eta_x}{(1-\rho)\sqrt K},
&
b_y=c_y
&=\frac{8\kappa\rho\eta_x}{(1-\rho)\sqrt K}.
\label{proof:DiMA:new_coefficients}
\end{align}
These choices satisfy
\begin{align}
\rho\left(b_x+\frac{\eta_x}{\sqrt K}\right)-b_x
&=0, \
\rho\left(b_y+\frac{8\kappa\eta_x}{\sqrt K}\right)-b_y
=0,
\end{align}
and
\begin{align}
\rho\beta\left(b_x+\frac{\eta_x}{\sqrt K}\right)
-\beta c_x&=0,\
\rho\beta\left(b_y+\frac{8\kappa\eta_x}{\sqrt K}\right)
-\beta c_y=0.
\end{align}
Thus, the terms involving the tracker and
networked momentum errors are eliminated.

The remaining contributions from the Lemmas~\ref{lemma:DiMA:consensus_Ux}
and~\ref{lemma:DiMA:deviation:momentum_network} are
\begin{align}
&(b_x+b_y)\Big[
KL_h d_i^2
+\left(L_f\sqrt K
+2\sigma_h\sqrt{\frac KN}\right)d_i
+2\beta\sigma\sqrt{\frac KN}
\Big] \notag\\
&=
\frac{\rho(1+8\kappa)\eta_x}{1-\rho}
\Big[
L_h\sqrt K\,d_i^2
+\left(L_f+\frac{2\sigma_h}{\sqrt N}\right)d_i
+\frac{2\beta\sigma}{\sqrt N}
\Big]. \notag
\end{align}
Therefore,
\begin{align}
&\bd{\Omega}_{i+1}-\bd{\Omega}_i \label{proof:DiMA:new_descent}\\
&\le
-\eta_x\mb{E}\Big[
\|\nabla P(\bx_{c,i})\|
+2L_f\delta^y_{c,i}
\Big]
+\frac{(2L_f+16\kappa L_f)\eta_x}{\sqrt K}
\mb{E}\mathcal E_i \notag\\
&\quad+
2\eta_x\mb{E}\|\widetilde{\bu}^x_{c,i}\|
+16\kappa\eta_x
\mb{E}\|\widetilde{\bu}^y_{c,i}\|
+\frac{L\eta_x^2}{2}
+2NL_f\eta_x\bar\eta_y \notag\\
&\quad+
\frac{\rho(1+8\kappa)\eta_x}{1-\rho}
\Big[
L_h\sqrt K\,d_i^2
+\left(L_f+\frac{2\sigma_h}{\sqrt N}\right)d_i
+\frac{2\beta\sigma}{\sqrt N}
\Big] \notag\\
&\quad+
\frac{8\kappa\eta_x}{\eta_y}
\mb{E}\left[
\Delta^y_{c,i}-\Delta^y_{c,i+1}
\right]
+8L_f\eta_x
\left(5\kappa\eta_y+2\kappa N\bar\eta_y\right). \notag
\end{align}
Summing over $i=0,\ldots,T-1$, dividing by $\eta_x T$,
and using $\bd{\Omega}_T\ge P^\star$ and
$\Delta^y_{c,T}\ge0$, we obtain
\begin{align}
&\frac1T\sum_{i=0}^{T-1}
\mb{E}\Big[
\|\nabla P(\bx_{c,i})\|
+2L_f\delta^y_{c,i}
\Big] \notag\\
&\le
\frac{\bd{\Omega}_0-P^\star}{\eta_xT}
+\frac{8\kappa\Delta^y_{c,0}}{\eta_yT}
+\frac{2L_f+16\kappa L_f}{T\sqrt K}
\sum_{i=0}^{T-1}\mb{E}\mathcal E_i \notag\\
&\quad+
\frac1T\sum_{i=0}^{T-1}
\left(
2\mb{E}\|\widetilde{\bu}^x_{c,i}\|
+16\kappa\mb{E}\|\widetilde{\bu}^y_{c,i}\|
\right)
+\frac{L\eta_x}{2}
+2NL_f\bar\eta_y \notag\\
&\quad+
\frac{\rho(1+8\kappa)}{1-\rho}
\Big[
L_h\sqrt K\,D_T
+\left(L_f+\frac{2\sigma_h}{\sqrt N}\right)\sqrt{D_T}
+\frac{2\beta\sigma}{\sqrt N}
\Big] \notag\\
&\quad+
8L_f\left(5\kappa\eta_y+2\kappa N\bar\eta_y\right),
\label{proof:DiMA:new_summed_descent}
\end{align}
where we used
$
\frac1T\sum_{i=0}^{T-1}d_i^2=D_T,
\frac1T\sum_{i=0}^{T-1}d_i\le\sqrt{D_T}.
$
We next use
Lemma~\ref{lemma:DiMA:deviation_ui} to bound each averaged
centroid momentum error by $Q_1$.
Moreover, identical initial models and
Lemma~\ref{lemma:DiMA:network_consensus} give
$
\frac1{T\sqrt K}
\sum_{i=0}^{T-1}\mb{E}\mathcal E_i
\le\frac{4\eta_y}{1-\rho}.
$
Therefore,
\begin{align}
&\frac1T\sum_{i=0}^{T-1}
\mb{E}\Big[
\|\nabla P(\bx_{c,i})\|
+2L_f\delta^y_{c,i}
\Big] \label{proof:DiMA:new_master_explicit}\\
&\le
\frac{\bd{\Omega}_0-P^\star}{\eta_xT}
+\frac{8\kappa\Delta^y_{c,0}}{\eta_yT}
+\frac{4(2L_f+16\kappa L_f)\eta_y}{1-\rho}
+(2+16\kappa)Q_1 \notag\\
&\quad+
\frac{L\eta_x}{2}
+2NL_f\bar\eta_y
+8L_f\left(5\kappa\eta_y+2\kappa N\bar\eta_y\right)
\notag\\
&\quad+
\frac{\rho(1+8\kappa)}{1-\rho}
\Big[
L_h\sqrt K\,D_T
+\left(L_f+\frac{2\sigma_h}{\sqrt N}\right)\sqrt{D_T}
+\frac{2\beta\sigma}{\sqrt N}
\Big].
\notag
\end{align}

To bound the initial potential, define
\begin{align}
H_{0,w}
\triangleq
\sqrt{
\frac1K\sum_{k=1}^{K}
\|\nabla_wJ_k(\bz_0)-\nabla_wJ(\bz_0)\|^2
}, \, w\in\{x,y\}. \notag
\end{align}
Initialize each local momentum with one unbiased
stochastic gradient and set its tracker equal to that
momentum, i.e., 
$\bu^x_{k,-1}={\bd{m}}^x_{k, -1}=\nabla_x Q_k(\bz_{k,-1};\bxi_{k,-1}), \bu^y_{k,-1}={\bd{m}}^y_{k,-1} =\nabla_y Q_k(\bz_{k,-1};\bxi_{k,-1}), \forall k \in [K]$. Then, we have
\begin{align}
\mb{E}\|\widetilde{\mM}^w_{-1}\|
&\le\sqrt K\,\sigma,\\
\mb{E}\|\mU^w_{-1}-\mU^w_{c,-1}\|
&\le\sqrt K\,(H_{0,w}+\sigma).
\end{align}
Using \eqref{proof:DiMA:new_coefficients}, we obtain
\begin{align}
\frac{\bd{\Omega}_0-P^\star}{\eta_xT}
&\le
\frac{P(\bx_{c,0})-P^\star}{\eta_xT}
+\frac{B_0}{T},
\end{align}
where
\begin{align}
B_0\triangleq
8L_f\delta^y_{c,0}
+\frac{\rho}{1-\rho}
\Big[
H_{0,x}+8\kappa H_{0,y}
+2(1+8\kappa)\sigma
\Big].
\label{proof:DiMA:new_initialization}
\end{align}
We now choose
\begin{align}
\eta_y=\beta=h,\quad
\eta_x=\frac{h}{2\kappa},\quad
\bar\eta_y=\frac{h}{NL_f},\quad
\bar\eta_x=\frac{\eta_x}{NL_f},
\notag
\end{align}
where $0<h\le1$.
By Lemma \ref{lemma:temporal:increment},
we have
\begin{align}
D_T
&\le6C_\tau\eta_y^2+4N^2\bar\eta_y^2
=\left(6C_\tau+\frac4{L_f^2}\right)h^2.
\end{align}
Substituting $\beta=\eta_y=h$ and
$\bar\eta_y=h/(NL_f)$ into the definition of $Q_1$ gives
\begin{align}
Q_1
&=
\frac{\sigma}{\sqrt K hT}
+L_h\left(3C_\tau+\frac{2}{L_f^2}\right)h
+
\Bigg(
\sigma_h\sqrt{6C_\tau+\frac{4}{L_f^2}}\notag\\
&\quad+\sigma
\Bigg)\sqrt{\frac{h}{NK}}
\le
\frac{\sigma}{\sqrt K hT}
+C\left(h+\sqrt{\frac{h}{NK}}\right),
\end{align}
where, using $C_\tau\le4$, we define
\begin{align}
C\triangleq
\max\left\{
L_h\left(12+\frac{2}{L_f^2}\right),
\;
\sigma_h\sqrt{24+\frac{4}{L_f^2}}+\sigma
\right\}.
\end{align}
This constant is independent of $K,N,T,\rho,\kappa$
and applies to both mixing strategies.
Substituting these estimates into
\eqref{proof:DiMA:new_master_explicit} yields
\begin{align}
&\frac1T\sum_{i=0}^{T-1}
\mb{E}\Big[
\|\nabla P(\bx_{c,i})\|
+2L_f\delta^y_{c,i}
\Big] \\
&\le
\frac{\kappa A_0}{hT}
+C'\kappa\left[ \frac{1}{hT}+
\sqrt{\frac{h}{NK}}
+\frac{h}{1-\rho}
+\frac{\sqrt K h^2}{1-\rho}
\right]
+\frac{B_0}{T}. \notag
\end{align}
where 
\begin{align}
&A_0\triangleq
2\big(P(\bx_{c,0})-P^\star\big)
+8\Delta^y_{c,0} \\
&C' \triangleq
\max \Big\{
18\sigma,
18C,113L_f+18+18C+9(L_f
\notag \\
&
+2\sigma_h)\sqrt D+18\sigma,
9L_hD
\Big\}, \quad
D\triangleq24+4/L_f^2
\end{align}
Here we used $\kappa\ge1$, $N,K\ge1$, $\rho<1$,
and $L\le2L_f\kappa$, and  $C'$ depends only on the smoothness factor and noise
constants and is independent of $K,N,T,\rho,\kappa$.
\end{proof}
\color{black}

\vspace{-1em}
\bibliographystyle{IEEEtran}
\bibliography{example_paper}

\clearpage 

\newcounter{sequation}
\newcommand{\svar}{%
    \stepcounter{sequation}%
    s\thesequation%
}
\renewcommand\theequation{s\arabic{equation}}
\section*{\bf {\Large Supplementary Material}}
\section*{\bf ``Communication-Efficient Decentralized Stochastic Minimax
Optimization'' by H. Cai, S. A. Alghunaim and A. H. Sayed}

\subsection{Proof of Lemma \ref{lemma:DiMA}}
\begin{proof}
By Lemma~\ref{lem:mix}, the mixing subroutine preserves the
network average. 
Hence,
\begin{align}
    \bx_{c,i+1} = \bx_{c,i} -  \frac{\eta_x}{K}\sum_{k=1}^K \mathcal{N}(\bu^x_{k,i}).
\end{align}
Using the $L$-smooth property of
$P(\cdot)$, we get
\begin{align}
\label{proof:DiMA:boundP}
&P(\bx_{c,i+1})\le P(\bx_{c,i})
    -
    \Big\langle
    \nabla P(\bx_{c,i}),
    \frac{\eta_x}{K} 
\sum_{k=1}^K\mathcal{N}(\bu^x_{k,i})
    \Big\rangle
+\frac{L\eta^2_x}{2}
    \notag \\
    &\le
    P(\bx_{c,i})
    - 
    \eta_x\Big\langle 
    \nabla P(\bx_{c,i})-\bu^x_{c,i},\frac{1}{K}
\sum_{k=1}^K\mathcal{N}(\bu^x_{k,i})
    \Big\rangle
     \notag\\
     &\quad 
  - \eta_x\Big\langle  \bu^x_{c,i},\frac{1}{K}
\sum_{k=1}^K\mathcal{N}(\bu^x_{k,i})
    \Big\rangle + \frac{L\eta^2_x}{2}
     \notag \\
&\overset{(a)}{\le} 
P(\bx_{c,i})
+\eta_x\|\nabla P(\bx_{c,i})-\bu^x_{c,i}\|
-\eta_{x} \Big \langle 
\bu^x_{c,i}, \mathcal{N}(\bu^x_{c,i})
\Big\rangle \notag \\
&\quad 
-\eta_x\Big\langle 
\bu^x_{c,i}, \frac{1}{K}\sum_{k=1}^K\Big(
\mathcal{N}(\bu^x_{k,i})-\mathcal{N}(\bu^x_{c,i}) \Big)
\Big\rangle + \frac{L\eta^2_x}{2}
\notag \\
&\overset{(b)}{\le}  P(\bx_{c,i})
+\eta_x\|\nabla P(\bx_{c,i}) - \bu^x_{c,i}\|
-\eta_x\|\bu^x_{c,i}\|
\notag \\
&\quad 
+
\frac{\eta_x}{K}\sum_{k=1}^{K}|\|\bu^x_{c,i}\| - \|\bu^x_{k,i}\||+\frac{L\eta^2_x}{2}
\notag \\
    & 
    \overset{(c)}{\le}
    P(\bx_{c,i}) 
    +
    2\eta_x\|\nabla P(\bx_{c,i}) - \bu^x_{c,i}\|
    -\eta_x \|\nabla P(\bx_{c,i})\|
    \notag\\
    &\quad 
    + \frac{\eta_x}{K}\sum_{k=1}^{K}
    \|\bu^x_{c,i} - \bu^x_{k,i}\| + \frac{L\eta^2_x}{2}
    \notag \\
    &\le
    P(\bx_{c,i})
    -\eta_x\|\nabla P(\bx_{c,i})\| 
    +2\eta_x\|\nabla P(\bx_{c,i}) -\nabla_x J(\bz_{c,i})\| \notag\\
    &\quad 
  + 2\eta_x\|\nabla_x J(\bz_{c,i})-\bu^x_{c,i}\|
 +\frac{\eta_x}{K}
 \sum_{k=1}^{K}
 \|\bu^x_{c,i} - \bu^x_{k,i}\|+\frac{L\eta^2_x}{2}  \notag \\
&\overset{(d)}{\le} 
P(\bx_{c,i}) - \eta_x\|\nabla P(\bx_{c,i})\|
+ 2L_f\eta_x
\delta^y_{c,i}
\notag\\
&\quad +2\eta_x
\|\nabla_x J(\bz_{c,i}) -\bu^x_{c,i} \|
+\frac{\eta_x}{\sqrt{K}}\|\mU^x_{c,i} - \mU^x_{i}\| + \frac{L\eta^2_x}{2} \notag \\
&\overset{(e)}{\le}
P(\bx_{c,i}) -\eta_x\|\nabla P(\bx_{c,i})\|+ 2L_f\eta_x
\delta^y_{c,i}+ \frac{2\eta_xL_f\mathcal{E}_i}{\sqrt{K}}
\notag \\
&\quad 
 +
2\eta_x
\Big\|\frac{1}{K}\sum_{k=1}^{K}
\nabla^x_{k,i} 
- \bu^x_{c,i}\Big\|
 +
\frac{\eta_x}{\sqrt{K}}\|\mU^x_{c,i} - \mU^x_{i}\|+ \frac{L\eta^2_x}{2}, \notag
\end{align}
Where $(a)$ follows from the Cauchy--Schwarz inequality and
$\big\|K^{-1}\sum_{k=1}^K \mathcal{N}(\bu^x_{k,i})\big\|\le 1$.
For step $(b)$, we use
$\langle\bu^x_{c,i},\mathcal{N}(\bu^x_{c,i})\rangle=\|\bu^x_{c,i}\|$
and $\bu^x_{c,i}=K^{-1}\sum_{k=1}^K\bu^x_{k,i}$ to obtain
\begin{align}
&-\eta_x\Big\langle
\bu^x_{c,i},
\frac{1}{K}\sum_{k=1}^K \mathcal{N}(\bu^x_{k,i})-\mathcal{N}(\bu^x_{c,i})
\Big\rangle \notag\\
&\quad =
-\frac{\eta_x}{K}\sum_{k=1}^K
\Big\langle \mathcal{N}(\bu^x_{c,i}),
\|\bu^x_{c,i}\|\mathcal{N}(\bu^x_{k,i})-\bu^x_{k,i}
\Big\rangle \notag\\
&\quad = -\frac{\eta_x}{K}\sum_{k=1}^K \Big\langle\mathcal{N}(\bu^x_{c,i}),
\bigl(\|\bu^x_{c,i}\|-\|\bu^x_{k,i}\|\bigr)
\mathcal N(\bu^x_{k,i}) \Big\rangle  \notag\\
&\quad \le
\frac{\eta_x}{K}\sum_{k=1}^K
\Big|\|\bu^x_{c,i}\|-\|\bu^x_{k,i}\|\Big|.
\end{align}
The last inequality follows from Cauchy--Schwarz and
$\|\mathcal{N}(\bg)\|\le 1$ for a given $\bg$.
Step $(c)$ follows from 
$-\|b\| \le \|a-b\| - \|a\|$ and $|\|a\| - \|b\|| \le \|a-b\|$, step
$(d)$ follows from the Cauchy-Schwarz inequality, i.e.,
\begin{align}
\label{proof:DiMA:cauchyschwarz}
&\sum_{k=1}^{K}
\|\bu^x_{c,i} - \bu^x_{k,i}\| = \sum_{k=1}^{K} 1 \cdot
\|\bu^x_{c,i} - \bu^x_{k,i}\|  \\
&= 
\mathds{1}^\top_{K}
\mbox{\rm col}\{\|\bu^x_{c,i} - \bu^x_{k,i}\|\}_{k=1}^{K} \notag \\
&\le 
\sqrt{1^2 + \cdots + 1^2}\sqrt{\sum_{k=1}^{K}\|\bu^x_{c,i}-\bu^x_{k,i}\|^2}
\le 
\sqrt{K}\|\mU^x_{c,i}-\mU^x_{i}\|,\notag
\end{align}
and step $(e)$ is due to 
\begin{align}
&\|\nabla_x J(\bz_{c,i}) - \bu^x_{c,i}\| \le 
\Big\|\nabla_x J(\bz_{c,i}) 
- \frac{1}{K}\sum_{k=1}^{K}
\nabla^x_{k,i}\Big\|
\notag\\
&\quad 
+ 
\Big\|\frac{1}{K}\sum_{k=1}^{K}
\nabla^x_{k,i}
- \bu^x_{c,i}\Big\| \le 
\frac{1}{K}
\sum_{k=1}^{K} 
\| 
\nabla_x J_k(\bz_{c,i})
-\nabla^x_{k,i}\|
 \notag \\
&+ 
\Big\|\frac{1}{K}\sum_{k=1}^{K}
\nabla^x_{k,i}
- \bu^x_{c,i}\Big\|\le
\frac{L_f\mathcal{E}_{i}}{\sqrt{K}}
+
\Big\|\frac{1}{K}\sum_{k=1}^{K}
\nabla^x_{k,i}
- \bu^x_{c,i}\Big\|,
\label{eq:deviation:ux_centroid}
\end{align}
where the last inequality follows from $L_f$-smooth assumption and an argument similar to \eqref{proof:DiMA:cauchyschwarz}.
Finally, 
we use 
\begin{align}
&\Big\|\frac{1}{K}
\sum_{k=1}^{K} \nabla^x_{k,i} - \bu^x_{c,i}\Big\|\le 
\Big\|\frac{1}{K}
\sum_{k=1}^{K} \nabla^x_{k,i} - \frac{1}{KN}\sum_{k=1}^K\sum_{n=0}^{N-1} \nabla^x_{k,i,n}\Big\| \notag \\
&\quad 
+ \Big\|\bu^x_{c,i}-\frac{1}{KN}\sum_{k=1}^K\sum_{n=0}^{N-1} \nabla^x_{k,i,n}
\Big\| \notag\\
&\overset{(a)}{\le} 
\frac{1}{KN}\sum_{k=1}^{K}
\sum_{n=0}^{N-1}
\|[\nabla^x_{k,i};\nabla^y_{k,i}] - [\nabla^x_{k,i,n};\nabla^y_{k,i,n}]\|
\notag\\
&\quad 
+ \Big\|\bu^x_{c,i}-\frac{1}{KN}\sum_{k=1}^K\sum_{n=0}^{N-1} \nabla^x_{k,i,n}
\Big\| \notag\\
&\overset{(b)}{\le} 
\frac{L_f}{KN}
\sum_{k=1}^{K}
\sum_{n=0}^{N-1}
\|\bz_{k,i} - \bz_{k,i,n}\|
+ \Big\|\bu^x_{c,i}-\frac{1}{KN}\sum_{k=1}^K\sum_{n=0}^{N-1} \nabla^x_{k,i,n}
\Big\| \notag\\
&\overset{(c)}{\le} 
\frac{L_f}{KN}
\sum_{k=1}^{K}
\sum_{n=0}^{N-1}
2n\bar{\eta}_y
+ \Big\|\bu^x_{c,i}-\frac{1}{KN}\sum_{k=1}^K\sum_{n=0}^{N-1} \nabla^x_{k,i,n}
\Big\| \notag \\
&\overset{(d)}{\le} NL_f\bar{\eta}_y+ \Big\|\bu^x_{c,i}-\frac{1}{KN}\sum_{k=1}^K\sum_{n=0}^{N-1} \nabla^x_{k,i,n}
\Big\|,
\end{align}
where $(a)$ follows from Jensen's inequality and $\|u\| \le \|[u;v]\|, \forall v$,
$(b)$ follows from the $L_f$-smooth assumption, $(c)$ follows from the local normalized learning rule and the choice $\bar{\eta}_x \le \bar{\eta}_y$, 
$(d)$ follows from 
$\sum_{n=0}^{N-1}2n \le N^2$.
Incorporating the above results into \eqref{proof:DiMA:boundP}, we complete the proof.
\end{proof}

\subsection{Proof of Lemma \ref{lemma:DiMA:network_consensus}}
\begin{proof}
Setting identical values for 
$\bx_{k,0},\bx_{k,-1}  \in \mR^{M_1}, \forall k \in [K]$
and $\by_{k,0}, \by_{k,-1}  \in \mR^{M_2}, \forall k \in [K]$, the initial consensus error vanishes, we have  $\mathcal{E}_0= \mathcal{E}_{-1} = 0$.
\color{black}
By Lemma~\ref{lem:mix}, the mixing subroutine preserves the
network average; we get 
\begin{align}
    \mX_{c,i+1}  = \mc{J} (\mX_{i} - \eta_x \overline{\mU}^x_{i}).
\label{proof:DiMA:NetworkXcentroid}
\end{align}
Invoking Lemma \ref{lem:mix} again, 
we can derive
\begin{align}
\label{proof:consenusX}
&\|\mX_{c,i+1}-\mX_{i+1}\| \notag\\
&=
\left\|
\mc{J}(\mX_i-\eta_x\overline{\mU}_i^x)
-\operatorname{Mix}
(\mX_i-\eta_x\overline{\mU}_i^x;A,r,\tau)
\right\| \notag\\
&\overset{(a)}{\leq}
\rho\left\|
(I-\mc{J})(\mX_i-\eta_x\overline{\mU}_i^x)
\right\| \notag\\
&=
\rho\left\|
(\mX_i-\mX_{c,i})
-\eta_x(\overline{\mU}_i^x-\overline{\mU}_{c,i}^x)
\right\| \notag\\
&\overset{(b)}{\leq}
\rho\|\mX_i-\mX_{c,i}\|
+\eta_x\rho
\|\overline{\mU}_i^x-\overline{\mU}_{c,i}^x\|,
\end{align}
where $(a)$ follows from \eqref{eq:subroutine:contraction}, and
$(b)$ follows from Assumption \ref{assumption:DiMA}
and the triangle inequality.
\color{black}
For $\|\overline{\mU}^x_{i}  - \overline{\mU}^x_{c,i}\|$, we have
\begin{align}
&\|\overline{\mU}^x_{i} - \overline{\mU}^x_{c,i}\|^2 \notag \\
&=
\Bigg\| 
\begin{bmatrix}
    \mathcal{N}(\bu^x_{1,i})- \frac{1}{K}\sum_{k=1}^{K} \mathcal{N}(\bu^x_{k,i}) \notag \\
    \vdots \\
     \mathcal{N}(\bu^x_{k,i}) - \frac{1}{K}\sum_{k=1}^{K} \mathcal{N}(\bu^x_{k,i})
\end{bmatrix}\Bigg\|^2 \notag \\
&=
\sum_{j=1}^{K}
\Big\| \frac{1}{K}\sum_{k=1}^{K} \Big(\mathcal{N}(\bu^x_{j,i}) - \mathcal{N}(\bu^x_{k,i})\Big)\Big\|^2 \notag \\
&\overset{(a)}{\le}
\sum_{j=1}^{K}
\frac{1}{K}\sum_{k=1}^{K}\Big\|\mathcal{N}(\bu^x_{j,i}) - \mathcal{N}(\bu^x_{k,i})\Big\|^2\overset{(b)}{\le}\sum_{j=1}^{K}
\frac{1}{K}\sum_{k=1}^{K}4 \le 4K,
\end{align}
where $(a)$ and $(b)$ follow from Jensen's inequality.
Therefore
\begin{align}
\|\overline{\mU}^x_{i} - \overline{\mU}^x_{c,i}\| \le 2\sqrt{K}.
\end{align}
The similar results hold for  $\|\mY_{c,i+1} - \mY_{i+1}\|$.
Setting $\eta_x \le \eta_y$, 
we have 
\begin{align}
\mathcal{E}_{i+1} \le \rho 
\mathcal{E}_{i}
+4\sqrt{K}\rho\eta_y.
\end{align}
Iterating the above inequality to $i=0$
and using $\mathcal{E}_0 =0$, 
 we get 
 \begin{align}
     \mathcal{E}_{i+1} \le 4\sqrt{K} \eta_y \sum_{\tau=1}^{i+1}\rho^{\tau}
     \le \frac{4\sqrt{K}\eta_y}{1-\rho} \quad (\rho <1).
 \end{align}
\end{proof}
\subsection{Proof of Lemma \ref{lemma:DiMA:deviation_ui}}
\begin{proof}
By Lemma \ref{lem:mix}, the averaged momentum recursion is
\begin{align}
    \bu^x_{c,i} 
    =  \bu^x_{c,i-1}+\bm^x_{c,i}  - \bm^x_{c,i-1}.
\end{align}
Initializing the momentum and tracker identically therefore
ensures $\bu^x_{c,i}=\bm^x_{c,i}$ for all $i\ge-1$.
Averaging the local momentum updates gives
\begin{align}
\bm^x_{c,i}
=(1-\beta)(\bm^x_{c,i-1}+\bh^x_{c,i})
+\beta\bg^x_{c,i}.
\label{proof:DiMA:momentum_centroid}
\end{align}
Subtracting the average true gradient yields
\begin{align}
&\bu^x_{c,i} -\frac{1}{KN}\sum_{k=1}^{K}\sum_{n=0}^{N-1}
\nabla^x_{k,i,n} =\bm^x_{c,i} -\frac{1}{KN}\sum_{k=1}^{K}\sum_{n=0}^{N-1}
\nabla^x_{k,i,n} \notag \\
&=(1-\beta)\Big(\bm^x_{c,i-1}- \frac{1}{KN}\sum_{k=1}^{K}\sum_{n=0}^{N-1} \nabla^x_{k,i-1,n}\Big)
\\
&\quad 
+(1-\beta)
\Big(h^x_{c,i} - \frac{1}{KN}\sum_{k=1}^{K}\sum_{n=0}^{N-1}
(\nabla^x_{k,i,n}-
\nabla^x_{k,i-1,n})\Big) \notag\\
&\quad 
+ (1-\beta)(\bh^x_{c,i} - h^x_{c,i})
+\beta \Big(\bg^x_{c,i} - \frac{1}{KN}\sum_{k=1}^{K}\sum_{n=0}^{N-1}\nabla^x_{k,i,n}\Big). \notag
\end{align}
Repeatedly applying the recursion and using the triangle inequality gives
\begin{align}
&\mb{E}\Big\| \bu^x_{c,i} -\frac{1}{KN}\sum_{k=1}^{K}\sum_{n=0}^{N-1}
\nabla^x_{k,i,n}\Big\|  \\
&\le 
(1-\beta)^{i+1}
\mb{E}\Big\| 
\bm^x_{c,-1} - \frac{1}{NK}
\sum_{k=1}^{K} \sum_{n=0}^{N-1}\nabla^x_{k,-1,n})
\Big\|  \notag\\
&\quad 
+ \sum_{\tau=0}^{i} (1-\beta)^{i-\tau+1}\mb{E}\Big\|
h^x_{c,\tau} -
\frac{1}{NK}
\sum_{k=1}^{K}
\sum_{n=0}^{N-1} (\nabla^x_{k,\tau,n} 
\notag \\
&\quad - 
\nabla^x_{k,\tau-1,n}
)\Big\|
+
\mb{E}
\Big\| \sum_{\tau=0}^{i}
(1-\beta)^{i-\tau+1}
(\bh^x_{c,\tau} - h^x_{c,\tau})\Big\|
\notag \\
&\quad+\mb{E}
\Big\| \beta \sum_{\tau=0}^{i}
(1-\beta)^{i-\tau}
\Big(
\bg^x_{c,\tau} - \frac{1}{KN}\sum_{k=1}^{K}\sum_{n=0}^{N-1}
\nabla^x_{k,\tau,n}
\Big)
\Big\|. \notag
\end{align}
For simplicity, let $\bm^x_{k,-1} = 
\nabla_x Q_k(\bz_{k,-1};\bd{\xi}_{k,-1}), \forall k \in [K]$ and let 
$\bz_{k,-1}, \bz_{k,-1,n}, \forall k \in [K], \forall n \in [N]$ be initialized with the same value, it follows that 
\begin{align}
&\mb{E}\Big\| 
\bm^x_{c,-1} - \frac{1}{NK}
\sum_{k=1}^{K} \sum_{n=0}^{N-1}\nabla^x_{k,-1,n}
\Big\|\le \mb{E}\Big\| 
\bm^x_{c,-1} - \frac{1}{K}
\sum_{k=1}^{K} \nabla^x_{k,-1}
\Big\|  \notag \\
&\overset{(a)}{\le} 
\sqrt{ 
\mb{E}\Big\| 
\frac{1}{K}\sum_{k=1}^{K}
(\nabla_x Q_k(\bz_{k,-1};\bd{\xi}_{k,-1}) - \nabla^x_{k,-1})
\Big\|^2
} \overset{(b)}{\le} \frac{\sigma}{\sqrt{K}},
\end{align}
where $(a)$ follows from 
the Jensen's inequality for concave function, and $(b)$ follows from
the fact that $\bd{\xi}_{k, -1}$ are sampled i.i.d. over $k \in [K]$.

\color{black}
We proceed to bound the 
second term. 
The triangle inequality and
Lemma~\ref{lemma:meanvalue} give
\begin{align}
&\mb{E}\Big\|
h^x_{c,\tau}
-\frac1{KN}\sum_{k=1}^{K}\sum_{n=0}^{N-1}
(\nabla^x_{k,\tau,n}-\nabla^x_{k,\tau-1,n})
\Big\| \notag\\
&\le
\frac1{KN}\sum_{k=1}^{K}\sum_{n=0}^{N-1}
\mb{E}\Big\|
[h^x_{k,\tau,n};h^y_{k,\tau,n}]
-\nabla^z_{k,\tau,n}
+\nabla^z_{k,\tau-1,n}
\Big\| \notag\\
&\le
\frac{L_h}{2KN}\sum_{k=1}^{K}\sum_{n=0}^{N-1}
\mb{E}\|\bz_{k,\tau,n}-\bz_{k,\tau-1,n}\|^2.
\end{align}
Summing over $i$ and collecting the terms for each
fixed $\tau$ therefore gives
\begin{align}
&\frac1T\sum_{i=0}^{T-1}\sum_{\tau=0}^{i}
(1-\beta)^{i-\tau+1}
\mb{E}\Big\|
h^x_{c,\tau}
-\frac1{KN}\sum_{k=1}^{K}\sum_{n=0}^{N-1}
(\nabla^x_{k,\tau,n}\notag\\
&\quad -\nabla^x_{k,\tau-1,n})
\Big\| \notag\\
&\le
\frac{L_h}{2KNT}
\sum_{\tau=0}^{T-1}\sum_{k=1}^{K}\sum_{n=0}^{N-1}
\mb{E}\|\bz_{k,\tau,n}-\bz_{k,\tau-1,n}\|^2
\sum_{i=\tau}^{T-1}(1-\beta)^{i-\tau+1}
\notag\\
&\le\frac{L_hD_T}{2\beta},
\label{proof:consensusuxci}
\end{align}
since the geometric sum is at most $1/\beta$.
Moreover,
Assumption~\ref{assumption:boundvariance} gives
\[
\mb{E}\|\bh^x_{k,\tau,n}-h^x_{k,\tau,n}\|^2
\le
\sigma_h^2\mb{E}
\|\bz_{k,\tau,n}-\bz_{k,\tau-1,n}\|^2.
\]
Thus,
\begin{align}
&\mb{E}\Big\|
\sum_{\tau=0}^{i}(1-\beta)^{i-\tau+1}
(\bh^x_{c,\tau}-h^x_{c,\tau})
\Big\| \notag\\
&\le
\sqrt{
\frac{\sigma_h^2}{N^2K^2}
\sum_{\tau=0}^{i}(1-\beta)^{2(i-\tau+1)}
\sum_{k=1}^{K}\sum_{n=0}^{N-1}
\mb{E}\|\bz_{k,\tau,n}-\bz_{k,\tau-1,n}\|^2
},
\end{align}
where the cross-term is zero
because of Assumption~\ref{assumption:boundvariance} and i.i.d. samples.
Applying Jensen's inequality to the average over $i$
and reordering the finite sums yields
\begin{align}
&\frac1T\sum_{i=0}^{T-1}
\mb{E}\Big\|
\sum_{\tau=0}^{i}(1-\beta)^{i-\tau+1}
(\bh^x_{c,\tau}-h^x_{c,\tau})
\Big\| \notag\\
&\le
\Bigg(
\frac{\sigma_h^2}{N^2K^2T}
\sum_{\tau=0}^{T-1}\sum_{k=1}^{K}\sum_{n=0}^{N-1}
\mb{E}\|\bz_{k,\tau,n}-\bz_{k,\tau-1,n}\|^2 \notag\\
&\quad
\sum_{i=\tau}^{T-1}(1-\beta)^{2(i-\tau+1)}
\Bigg)^{1/2}
\le
\frac{\sigma_h\sqrt{D_T}}
{\sqrt{NK\beta}},
\end{align}
where
\[
\sum_{i=\tau}^{T-1}(1-\beta)^{2(i-\tau+1)}
\le
\frac1{1-(1-\beta)^2}
\le\frac1{{}\beta} \quad (\beta \le 1).
\]
and $D_T$ is defined in \eqref{equation:definition:DT}.
\color{black}
For the last term,
we can bound it through
\begin{align}
&\mb{E}
\Big\|\beta \sum_{\tau=0}^{i}
(1-\beta)^{i-\tau}
\Big(
\bg^x_{c,\tau} - \frac{1}{KN}\sum_{k=1}^{K}\sum_{n=0}^{N-1}\nabla^x_{k,\tau,n}
\Big)
\Big\| \notag \\
& \overset{(a)}{\le}  \beta
\sqrt{\mb{E}
\Big\| 
\frac{1}{NK}
\sum_{k=1}^{K}
\sum_{n=0}^{N-1}
\sum_{\tau=0}^{i}
(1-\beta)^{i-\tau}
\Big(
\bg^x_{k,\tau,n} -\nabla^x_{k,\tau,n} \Big)
\Big\|^2} \notag \\
& \overset{(b)}{\le}  \beta
\sqrt{\frac{1}{N^2K^2}\sum_{k=1}^{K}
\sum_{n=0}^{N-1}
\sum_{\tau=0}^{i}
(1-\beta)^{2(i-\tau)}
\mb{E}\Big\| 
\bg^x_{k,\tau,n} -\nabla^x_{k,\tau,n} 
\Big\|^2} \notag \\
& \overset{(c)}{\le} 
\beta \sqrt{\frac{\sigma^2}{NK\beta}} \le \frac{\sigma \sqrt{\beta}}{\sqrt{NK}},
\end{align}
where $(a)$ follows from $
\bg^x_{c,\tau} = \frac{1}{KN}\sum_{k=1}^{K}\sum_{n=0}^{N-1}
\bg^x_{k,\tau,n}$,
 $(b)$ follows from the fact $\bxi_{k,\tau,n}$ are i.i.d.
and independent of $\bz_{k,\tau, n}$, and $(c)$ follows from $\sum_{\tau=0}^{i}(1-\beta)^{2(i-\tau)}  \le \frac{1}{1-(1-\beta)^2} \le \frac{1}{\beta}$ for $\beta \le 1$.
\color{black}
Combining the bounds and using \eqref{proof:DiMA:local_temporal_difference0}, we obtain
\begin{align}
&\frac1T\sum_{i=0}^{T-1}
\mb{E}\Big\|
\bu^x_{c,i}
-\frac1{KN}\sum_{k=1}^{K}\sum_{n=0}^{N-1}
\nabla^x_{k,i,n}
\Big\| \notag\\
&\le
\frac{\sigma}{\sqrt K\beta T}
+\frac{L_hD_T}{2\beta}
+\frac{\sigma_h\sqrt{D_T}}
{\sqrt{NK\beta}}
+\frac{\sigma\sqrt{\beta}}{\sqrt{NK}} \notag\\
&\le
\frac{\sigma}{\sqrt K\beta T}
+\frac{L_h(3C_{\tau}\eta_y^2
+2N^2\bar{\eta}_y^2)}{\beta} \notag\\
&\quad 
+
\frac{\sigma_h
\sqrt{6C_{\tau}\eta_y^2+4N^2\bar{\eta}_y^2}}
{\sqrt{NK\beta}}
+\frac{\sigma\sqrt{\beta}}{\sqrt{NK}}.
\end{align}
The $y$-component satisfies the same bound by an
identical argument.
\end{proof}

\subsection{Proof of Lemma \ref{lemma:DiMA:consensus_Ux}}
\color{black}
\begin{proof}
We first prove the result for the $x$-component.
By Lemma~\ref{lem:mix},
\begin{align}
\mU^x_{c,i}
=\mathcal{J}\Big(
\mU^x_{i-1}+\mM^x_i-\mM^x_{i-1}
\Big).
\label{proof:DiMA:tracking_centoid}
\end{align}
Applying the contraction property in the same lemma gives
\begin{align}
\label{proof:DiMA:tracking_recursion}
&\|\mU^x_i-\mU^x_{c,i}\| =
\Big\|
\operatorname{Mix}\Big(
\mU^x_{i-1}+\mM^x_i-\mM^x_{i-1};A,r,\tau
\Big)\notag\\
&
-\mathcal{J}\Big(
\mU^x_{i-1}+\mM^x_i-\mM^x_{i-1}
\Big)
\Big\| \notag\\
&\overset{(a)}{\leq}
\rho\Big\|
(I-\mathcal{J})\Big(
\mU^x_{i-1}+\mM^x_i-\mM^x_{i-1}
\Big)
\Big\| \notag\\
&\overset{(b)}{\leq}
\rho\|\mU^x_{i-1}-\mU^x_{c,i-1}\|
+\rho\|(I-\mathcal{J})(\mM^x_i-\mM^x_{i-1})\| \notag\\
&\overset{(c)}{\leq}
\rho\|\mU^x_{i-1}-\mU^x_{c,i-1}\|
+\rho\|\mM^x_i-\mM^x_{i-1}\|,
\end{align}
where $(a)$ follows from Lemma~\ref{lem:mix}, $(b)$ follows
from the triangle inequality and
\((I-\mathcal{J})\mU^x_{i-1}
=\mU^x_{i-1}-\mU^x_{c,i-1}\), and $(c)$ uses
\(\|I-\mathcal{J}\|\leq1\).

For each outer iteration $i$, define
\begin{align}
d_i^2
&\triangleq
\frac{1}{KN}
\sum_{k=1}^{K}\sum_{n=0}^{N-1}
\mb{E}\|\bz_{k,i,n}-\bz_{k,i-1,n}\|^2,
\qquad
d_i\triangleq\sqrt{d_i^2}.
\label{proof:DiMA:tracking:temporal_difference}
\end{align}
Recall also that
\begin{align}
\widetilde{\mM}^x_{i-1}
\triangleq
\mM^x_{i-1}
-\operatorname{col}\Big\{
\frac1N\sum_{n=0}^{N-1}\nabla^x_{k,i-1,n}
\Big\}_{k=1}^{K}.
\end{align}
The momentum update can be written as
\begin{align}
\mM^x_i-\mM^x_{i-1}
&=-\beta\widetilde{\mM}^x_{i-1}
+(1-\beta)\bd{\Gamma}_2
+(1-\beta)\bd{\Gamma}_3 \notag\\
&\quad+\beta\bd{\Gamma}_4+\bd{\Gamma}_5,
\label{proof:DiMA:notation:Mxi}
\end{align}
where
\begin{align}
\bd{\Gamma}_2
&\triangleq
\mnH^x_i-
\operatorname{col}\Big\{
\frac1N\sum_{n=0}^{N-1}
(\nabla^x_{k,i,n}-\nabla^x_{k,i-1,n})
\Big\}_{k=1}^{K},\\
\bd{\Gamma}_3
&\triangleq \mH^x_i-\mnH^x_i,\\
\bd{\Gamma}_4
&\triangleq
\mG^x_i-
\operatorname{col}\Big\{
\frac1N\sum_{n=0}^{N-1}\nabla^x_{k,i,n}
\Big\}_{k=1}^{K},\\
\bd{\Gamma}_5
&\triangleq
\operatorname{col}\Big\{
\frac1N\sum_{n=0}^{N-1}
(\nabla^x_{k,i,n}-\nabla^x_{k,i-1,n})
\Big\}_{k=1}^{K}.
\end{align}
We bound these four terms separately.

First, the triangle inequality and
Lemma~\ref{lemma:meanvalue} yield
\begin{align}
\mb{E}\|\bd{\Gamma}_2\|
&\le
\frac1N\sum_{k=1}^{K}\sum_{n=0}^{N-1}
\mb{E}\Big\|
[h^x_{k,i,n};h^y_{k,i,n}]
-\nabla^z_{k,i,n}
+\nabla^z_{k,i-1,n}
\Big\| \notag\\
&\le
\frac{L_h}{2N}
\sum_{k=1}^{K}\sum_{n=0}^{N-1}
\mb{E}\|\bz_{k,i,n}-\bz_{k,i-1,n}\|^2 =\frac{KL_h}{2}d_i^2.
\label{eq:deviation:taylor}
\end{align}
For the stochastic terms, each evaluation point and HVP
direction is determined before its fresh sample is drawn.
Consequently, drawing i.i.d. samples and using Assumption~\ref{assumption:boundvariance}
gives
\begin{align}
&\mb{E}\|\bd{\Gamma}_3\|
\le
\sqrt{\mb{E}\|\mH^x_i-\mnH^x_i\|^2}
\notag\\
&\le
\left(
\frac1{N^2}
\sum_{k=1}^{K}\sum_{n=0}^{N-1}
\mb{E}\|\bh^x_{k,i,n}-h^x_{k,i,n}\|^2
\right)^{1/2}
\notag\\
&\le
\left(
\frac{\sigma_h^2}{N^2}
\sum_{k=1}^{K}\sum_{n=0}^{N-1}
\mb{E}\|\bz_{k,i,n}-\bz_{k,i-1,n}\|^2
\right)^{1/2}
 =\sigma_h\sqrt{\frac KN}\,d_i,
\label{eq:deviation:Hxi}
\end{align}
and
\begin{align}
\mb{E}\|\bd{\Gamma}_4\|
&\le
\left(
\frac1{N^2}
\sum_{k=1}^{K}\sum_{n=0}^{N-1}
\mb{E}\|\bg^x_{k,i,n}-\nabla^x_{k,i,n}\|^2
\right)^{1/2} \le \sigma\sqrt{\frac KN}.
\label{proof:DiMA:deviation:localgradient}
\end{align}
Finally, Jensen's inequality and $L_f$-Lipschitz Assumption give
\begin{align}
\mb{E}\|\bd{\Gamma}_5\|
&\le
\left(
\frac1N
\sum_{k=1}^{K}\sum_{n=0}^{N-1}
\mb{E}\|
\nabla^x_{k,i,n}-\nabla^x_{k,i-1,n}
\|^2
\right)^{1/2}
\notag\\
&\le
L_f\left(
\frac1N
\sum_{k=1}^{K}\sum_{n=0}^{N-1}
\mb{E}\|\bz_{k,i,n}-\bz_{k,i-1,n}\|^2
\right)^{1/2}
\notag\\
&=L_f\sqrt K\,d_i.
\end{align}

Combining these bounds with
\eqref{proof:DiMA:notation:Mxi} yields
\begin{align}
&\mb{E}\|\mM^x_i-\mM^x_{i-1}\| \notag\\
&\le
\beta\mb{E}\|\widetilde{\mM}^x_{i-1}\|
+\frac{(1-\beta)KL_h}{2}d_i^2
\notag\\
&\quad+
\left(
L_f\sqrt K
+(1-\beta)\sigma_h\sqrt{\frac KN}
\right)d_i
+\beta\sigma\sqrt{\frac KN}.
\end{align}
Since \(0<\beta\le1\), we have
$
0\le1-\beta\le1.
$
Substituting into
\eqref{proof:DiMA:tracking_recursion}, we obtain
\begin{align}
&\mb{E}\|\mU^x_i-\mU^x_{c,i}\| \notag\\
&\le
\rho\mb{E}\|\mU^x_{i-1}-\mU^x_{c,i-1}\|
+\rho\beta\mb{E}\|\widetilde{\mM}^x_{i-1}\|+
\frac{\rho KL_h}{2}d_i^2
\notag\\
&\quad
+\rho\left(
L_f\sqrt K
+\sigma_h\sqrt{\frac KN}
\right)d_i +\rho\beta\sigma\sqrt{\frac KN}.
\end{align}
The $y$-component follows by the same argument.
\end{proof}
\color{black}
\subsection{Proof of Lemma \ref{lemma:DiMA:deviation:momentum_network}}
\begin{proof}
The proof is similar to Lemma \ref{lemma:DiMA:consensus_Ux}.
\end{proof}

\subsection{Proof of Lemma \ref{lemma:DiMA:envelope_gap}}
\begin{proof}
We first recall the definition
$
\Delta^y_{c, i} 
\triangleq P(\bx_{c,i}) - J(\bx_{c,i}, \by_{c,i})$.
By Lemma \ref{lem:mix}, averaging the recursion for 
$\by_{k,i+1} $ over $k \in [K]$
gives
\begin{align}
\by_{c,i+1} = \by_{c,i}+\frac{\eta_y}{K}\sum_{k=1}^{K}
\mathcal{N}(\bu^y_{k,i}).
\end{align}
Using the $L_f$-smooth property for $-J(\bx_{c,i+1}, \cdot)$, we get
\begin{align}
&-J(\bx_{c,i+1}, \by_{c,i+1})  \\
&\le
-J(\bx_{c,i+1},\by_{c,i})
-\Big\langle
\nabla_y J(\bx_{c,i+1}, \by_{c,i}),
\frac{\eta_y}{K}\sum_{k=1}^{K}
\mathcal{N}(\bu^y_{k,i})
\Big\rangle
\notag\\
&\quad 
+\frac{L_f\eta^2_y}{2} \notag \\
&\overset{(a)}{=}  
-J(\bx_{c,i+1},\by_{c,i})
-\Big\langle \nabla_y J(\bx_{c,i+1},\by_{c,i}) -\bu^y_{c,i},  
\notag\\
&\quad
\frac{\eta_y}{K}\sum_{k=1}^{K}
\mathcal{N}(\bu^y_{k,i})\Big\rangle- \Big\langle \bu^y_{c,i} ,  \frac{\eta_y}{K}\sum_{k=1}^{K}
\mathcal{N}(\bu^y_{k,i})\Big\rangle   + \frac{L_f \eta^2_y}{2}
\notag \\
&\overset{(b)}{\le} 
-J(\bx_{c,i+1}, \by_{c,i})
+\eta_y\|\nabla_y J(\bx_{c,i+1},\by_{c,i}) -\bu^y_{c,i}\|
\notag\\
&\quad-\eta_{y} \Big \langle 
\bu^y_{c,i}, \frac{\bu^y_{c,i}}{\|\bu^y_{c,i}\|}
\Big\rangle -\eta_y\langle 
\bu^y_{c,i}, \frac{1}{K}\sum_{k=1}^K\Big(
\mathcal{N}(\bu^y_{k,i})\notag \\
&\quad -\frac{\bu^y_{c,i}}{\|\bu^y_{c,i}\|} \Big)
\Big\rangle
 + \frac{L_f\eta^2_y}{2} \notag\\
 &\overset{(c)}{\le} 
 -J(\bx_{c,i+1}, \by_{c,i})
+\eta_y\|\nabla_y J(\bx_{c,i+1},\by_{c,i}) -\bu^y_{c,i}\| \notag\\
&\quad 
-\eta_y\|\bu^y_{c,i}\|+
\frac{\eta_y}{K}\sum_{k=1}^{K}|\|\bu^y_{c,i}\| - \|\bu^y_{k,i}\||+\frac{L_f\eta^2_y}{2} \notag\\
 &\overset{(d)}{\le} 
 -J(\bx_{c,i+1}, \by_{c,i})
+\eta_y\|\nabla_y J(\bx_{c,i+1},\by_{c,i}) -\bu^y_{c,i}\| \notag\\
&\quad 
-\eta_y\|\bu^y_{c,i}\|+
\frac{\eta_y}{K}\sum_{k=1}^{K}\|\bu^y_{c,i} - \bu^y_{k,i}\|+\frac{L_f\eta^2_y}{2} \notag \\
&\overset{(f)}{\le} 
-J(\bx_{c,i+1}, \by_{c,i})
+2\eta_y\|\nabla_y J(\bx_{c,i+1},\by_{c,i}) -\bu^y_{c,i}\| \notag\\
&\quad 
-\eta_y\|\nabla_y J(\bx_{c,i+1},\by_{c,i})\|+
\frac{\eta_y}{K}\sum_{k=1}^{K}\|\bu^y_{c,i} - \bu^y_{k,i}\|+\frac{L_f\eta^2_y}{2}, \notag\\
&\overset{(g)}{\le} 
-J(\bx_{c,i+1},\by_{c,i})
+2\eta_y\eta_xL_f+2\eta_y\|\nabla_y J(\bx_{c,i},\by_{c,i}) - \bu^y_{c,i}\| \notag\\
&\quad 
-\nu\eta_y\|\by^o(\bx_{c,i+1};\by_{c,i})-\by_{c,i}\|
+\frac{\eta_y}{K}\sum_{k=1}^{K}\|\bu^y_{c,i} - \bu^y_{k,i}\|
+ \frac{L_f\eta^2_y}{2} \notag \\
&\overset{(h)}{\le}
-J(\bx_{c,i+1},\by_{c,i})
+2\eta_y\eta_xL_f+2\eta_y\|\nabla_y J(\bx_{c,i},\by_{c,i}) - \bu^y_{c,i}\| \notag\\
&\quad 
-\nu\eta_y\|\by^o(\bx_{c,i+1};\by_{c,i})-\by_{c,i}\|
+\frac{\eta_y}{\sqrt{K}}
\|\mU^y_{i} - \mU^y_{c,i}\|
+ \frac{L_f\eta^2_y}{2}
\end{align}
where $(a)-(d)$ follow from 
\eqref{proof:DiMA:boundP}, $(f)$ follows from 
$-\|\bu^y_{c,i}\| \le \| \nabla_y J(\bx_{c,i+1},\by_{c,i}) -\bu^y_{c,i}\|  -\|\nabla_y J(\bx_{c,i+1},\by_{c,i})\|$, $(g)$ follows from Lemma \ref{lemma:Danskin2} 
\begin{align}
\|\nabla_y J(\bx_{c,i+1},\by_{c,i})\|^2 &\ge \nu^2\|\by^o(\bx_{c,i+1};\by_{c,i})-\by_{c,i}\|^2, 
\end{align}
where
\begin{align}
\by^o(\bx_{c,i+1};\by_{c,i})
&\in\operatorname*{argmin}_{v\in Y^o(\bx_{c,i+1})}
\|\by_{c,i}-v\|,
\end{align}
and
$(h)$ is proved similar to \eqref{proof:DiMA:cauchyschwarz}.
Adding $P(\bx_{c,i+1})$ on both sides, we can derive that 
\begin{align}
&\Delta^y_{c,i+1} \notag
 \\
&\le \Delta^y_{c,i}
+P(\bx_{c,i+1})-P(\bx_{c,i})
+J(\bx_{c,i},\by_{c,i})-J(\bx_{c,i+1},\by_{c,i})
\notag\\
&\quad 
+2\eta_y\|\nabla_y J(\bx_{c,i},\by_{c,i}) - \bu^y_{c,i}\|
-\nu\eta_y\|\by^o(\bx_{c,i+1};\by_{c,i})\notag\\
&\quad -\by_{c,i}\|
+\frac{\eta_y}{\sqrt{K}}
\|\mU^y_i-\mU^y_{c,i}\|
+2\eta_y\eta_xL_f+\frac{L_f\eta^2_y}{2}.
\end{align}
Because $P(\cdot)$ is $L$-smooth and 
$-J(\cdot, \by_{c,i})$ is $L_f$-smooth, we have 
\begin{align}
&P(\bx_{c,i+1}) - P(\bx_{c,i})
+J(\bx_{c,i},\by_{c,i})-J(\bx_{c,i+1},\by_{c,i})
\notag\\
&\le \langle \nabla P(\bx_{c,i})-\nabla_x J(\bx_{c,i},\by_{c,i}),\bx_{c,i+1} -\bx_{c,i}\rangle
+\frac{(L+L_f)\eta^2_x}{2} \notag\\
&\le 
\eta_xL_f\|\by_{c,i} - \by^o(\bx_{c,i})\|
+\frac{(L+L_f)\eta^2_x}{2}.
\end{align}
Combining the results, we get
\begin{align}
&\Delta^y_{c,i+1} \\
&\le
\Delta^y_{c,i}+\eta_xL_f
\delta^y_{c,i}
-\nu\eta_y\|\by_{c,i}-\by^o(\bx_{c,i+1};\by_{c,i})\| \notag \\
&\quad 
+2\eta_y
\|\nabla_y J(\bx_{c,i},\by_{c,i}) -\bu^y_{c,i}\| + \frac{\eta_y}{\sqrt{K}} \|\mU^y_{i} - \mU^y_{c,i}\|\notag\\
&\quad 
+ \Big( 
\frac{L_f \eta^2_y}{2} + \frac{L_f \eta^2_x}{2} + \frac{L\eta^2_x}{2} + 2L_f \eta_x\eta_y
\Big).
\end{align}
Moving $\|\by_{c,i} - \by^o(\bx_{c,i+1};\by_{c,i})\|$ to the left hand side, we get 
\begin{align}
&\|\by_{c,i} - \by^o(\bx_{c,i+1};\by_{c,i})\|
\notag \\
&\le 
\frac{(\Delta^y_{c,i} - \Delta^y_{c,i+1})}{\eta_y \nu}
+ \frac{\kappa\eta_x }{\eta_y}
\delta^y_{c,i}+\frac{2}{\nu}\|\nabla_y J(\bx_{c,i}, \by_{c,i})- \bu^y_{c,i}\|
\notag\\
&\quad 
+ \frac{1}{\sqrt{K}\nu}\|\mU^y_{i}-\mU^y_{c,i}\|
+ 
\Big( 
\frac{\kappa\eta_y}{2} + \frac{\kappa\eta^2_x}{2\eta_y}
+ \frac{L\eta^2_x}{2\nu \eta_y} + 2\kappa \eta_x
\Big),
\end{align}
where $\kappa \triangleq L_f/\nu$.
By definition of $\by^o(\bx_{c,i+1})$
\[
\by^o(\bx_{c,i+1})
\in
\operatorname*{argmin}_{v\in Y^o(\bx_{c,i+1})}
\|\by_{c,i+1}-v\|,
\]
we proceed to derive that 
\begin{align}
&\|\by_{c,i+1} - \by^o(\bx_{c,i+1})\|  \le \|\by_{c,i+1} - \by^o(\bx_{c,i+1};\by_{c,i})\| \notag\\
&=\Big\|
\by_{c,i}+\frac{\eta_y}{K}\sum_{k=1}^{K}
\mathcal{N}(\bu^y_{k,i})
- \by^o(\bx_{c,i+1}; \by_{c,i})\Big\| \notag \\
&\le 
\|\by_{c,i} - \by^o(\bx_{c,i+1};\by_{c,i})\|+ \eta_y
\notag \\
& \le 
\frac{(\Delta^y_{c,i} - \Delta^y_{c,i+1})}{\eta_y \nu}
+ \frac{\kappa\eta_x}{\eta_y}
\delta^y_{c,i}+ 
\frac{2}{\nu}
\|\nabla_y J(\bx_{c,i}, \by_{c,i}) - \bu^y_{c,i}\|
\notag\\
&\quad 
+ \frac{1}{\sqrt{K}\nu} 
\|\mU^y_i - \mU^y_{c,i}\|
+\Big( 
\frac{\kappa\eta_y}{2} + \frac{\kappa\eta^2_x}{2\eta_y}
+ \frac{L\eta^2_x}{2\nu \eta_y} + 2\kappa \eta_x + \eta_y
\Big). \notag
\end{align}
We next choose
$ \frac{\kappa \eta_x}{\eta_y} \le \frac{1}{2} \Longrightarrow \eta_x \le \frac{\eta_y}{2\kappa}$ to obtain
\begin{align}
&\delta^y_{c,i+1} \notag \\
&\le 
\frac{1}{2}\delta^y_{c,i}
+ \frac{(\Delta^y_{c,i} - \Delta^y_{c,i+1})}{\eta_y \nu}+
\frac{2}{\nu}
\|\nabla_y J(\bx_{c,i}, \by_{c,i}) - \bu^y_{c,i}\|
\notag \\
&\quad  +\frac{1}{\sqrt{K}\nu}
\|\mU^y_i - \mU^y_{c,i}\|
 +\Big( 
\frac{\kappa\eta_y}{2} + \frac{\eta_x}{4}
 + 3\kappa \eta_x + \eta_y
\Big). \notag \\
\end{align}
where the last inequality follows from 
\begin{align}
\frac{L\eta^2_x}{2\nu\eta_y}\le 
\frac{(L_f+L_f\kappa)\eta^2_x}{2\nu\eta_y} \le \frac{\kappa^2 \eta^2_x}{\eta_y} \le \frac{\kappa \eta_x}{2} \le \kappa \eta_x.
\end{align}
Using $\kappa \geq1$ and $\eta_x \le \eta_y$, we have 
\begin{align}
\frac{\kappa\eta_y}{2} + \frac{\eta_x}{4}
 + 3\kappa \eta_x + \eta_y\le 5\kappa \eta_y.
 \end{align}
Next, we can derive that 
\begin{align}
&\|\nabla_y J(\bx_{c,i},\by_{c,i}) - \bu^y_{c,i}\| \notag\\
&\le
\Big\|\nabla_y J(\bx_{c,i},\by_{c,i}) - \frac{1}{K}\sum_{k=1}^{K}
\nabla^y_{k,i}\Big\|+\Big\|\frac{1}{K}
\sum_{k=1}^{K}\nabla^y_{k,i}
\notag\\
&\quad  - \frac{1}{KN}
\sum_{k=1}^{K}\sum_{n=0}^{N-1} \nabla^y_{k,i,n}\Big\|
+\Big\|\bu^y_{c,i} -\frac{1}{KN}
\sum_{k=1}^{K}\sum_{n=0}^{N-1} \nabla^y_{k,i,n}\Big\| \notag \\
&\overset{(a)}{\le}
\frac{L_f}{\sqrt{K}}
\mathcal{E}_i + \frac{1}{NK}
\sum_{k=1}^{K}\sum_{n=0}^{N-1}
\|\nabla^y_{k,i} - \nabla^y_{k,i,n}\|
\notag\\
&\quad 
+\Big\|\bu^y_{c,i} -\frac{1}{KN}
\sum_{k=1}^{K}\sum_{n=0}^{N-1} \nabla^y_{k,i,n}\Big\| \notag\\
&\overset{(b)}{\le} 
\frac{L_f}{\sqrt{K}}
\mathcal{E}_i+\frac{1}{NK}
\sum_{k=1}^{K}\sum_{n=0}^{N-1}2L_fn \bar{\eta}_y
\notag\\
&\quad 
+\Big\|\bu^y_{c,i} -\frac{1}{KN}
\sum_{k=1}^{K}\sum_{n=0}^{N-1} \nabla^y_{k,i,n}\Big\| \\
&\le 
\frac{L_f}{\sqrt{K}}
\mathcal{E}_i+ L_f N\bar{\eta}_y +\Big\|\bu^y_{c,i} -\frac{1}{KN}
\sum_{k=1}^{K}\sum_{n=0}^{N-1} \nabla^y_{k,i,n}\Big\|. \notag
\end{align}
where $(a)$ Jensen's inequality and the $L_f$-smooth assumption, $(b)$ follows from the $L_f$-smooth assumption and the local normalized learning rule.

Putting all the results together, we obtain
\begin{align}
&\delta^y_{c,i+1} \le 
\frac{1}{2}\delta^y_{c,i}
+ \frac{(\Delta^y_{c,i} - \Delta^y_{c,i+1})}{\eta_y \nu}+
\frac{2\kappa}{\sqrt{K}}\mathcal{E}_i+2\kappa N \bar{\eta}_y+5\kappa \eta_y
\notag \\
&\quad
+
\frac{2}{\nu}\Big\|\bu^y_{c,i} -\frac{1}{KN}
\sum_{k=1}^{K}\sum_{n=0}^{N-1} \nabla^y_{k,i,n}\Big\| +\frac{1}{\sqrt{K}\nu}
\|\mU^y_i - \mU^y_{c,i}\|. \notag \\
\end{align}
\end{proof}
\color{black}
\section{Additional experiments}
\begin{figure*}[t]
\vspace{-2em}
\centering
\textbf{(a) $\lambda=0.9$}\par\smallskip
\includegraphics[width=0.32\textwidth]{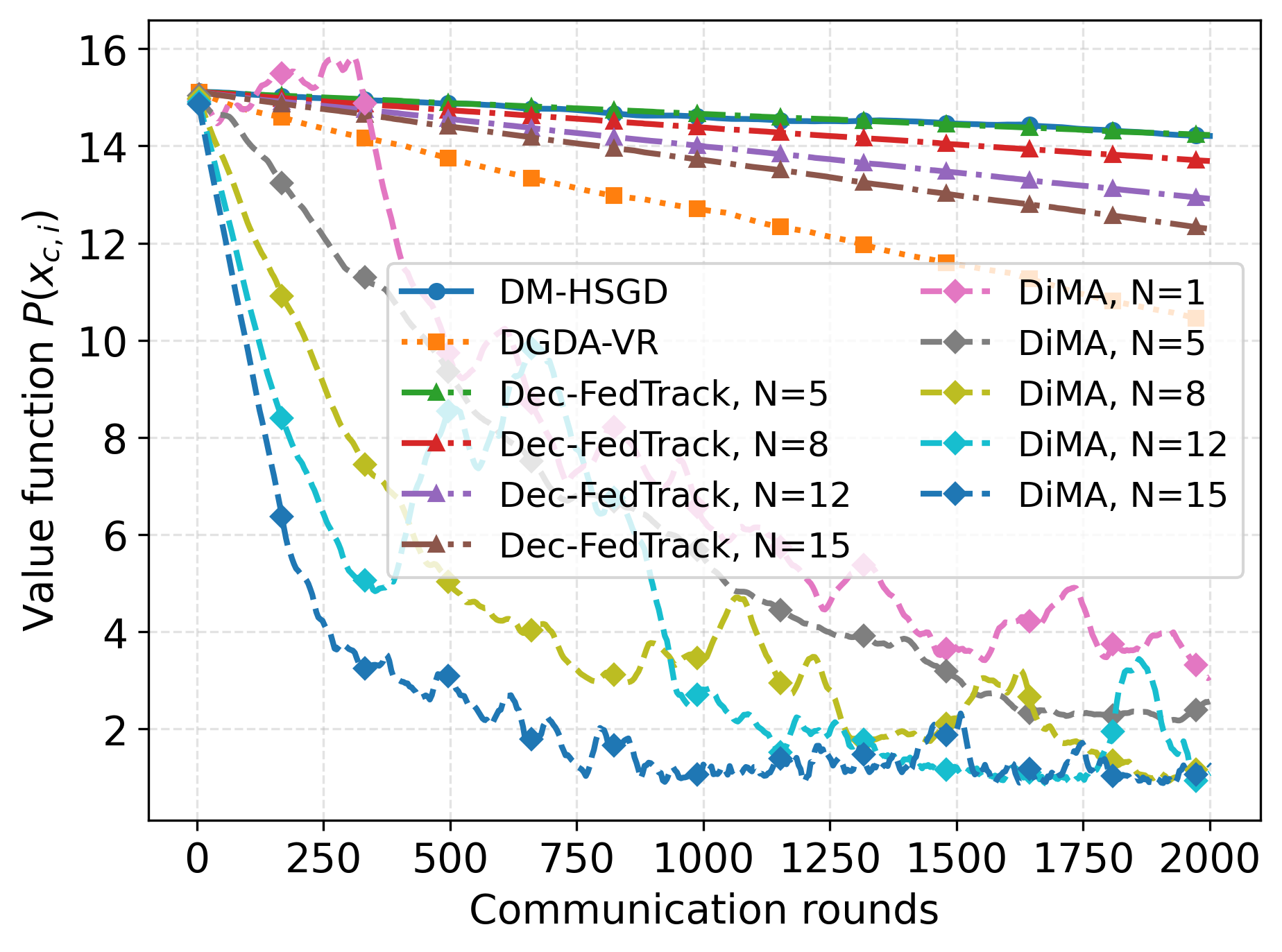}
\includegraphics[width=0.32\textwidth]{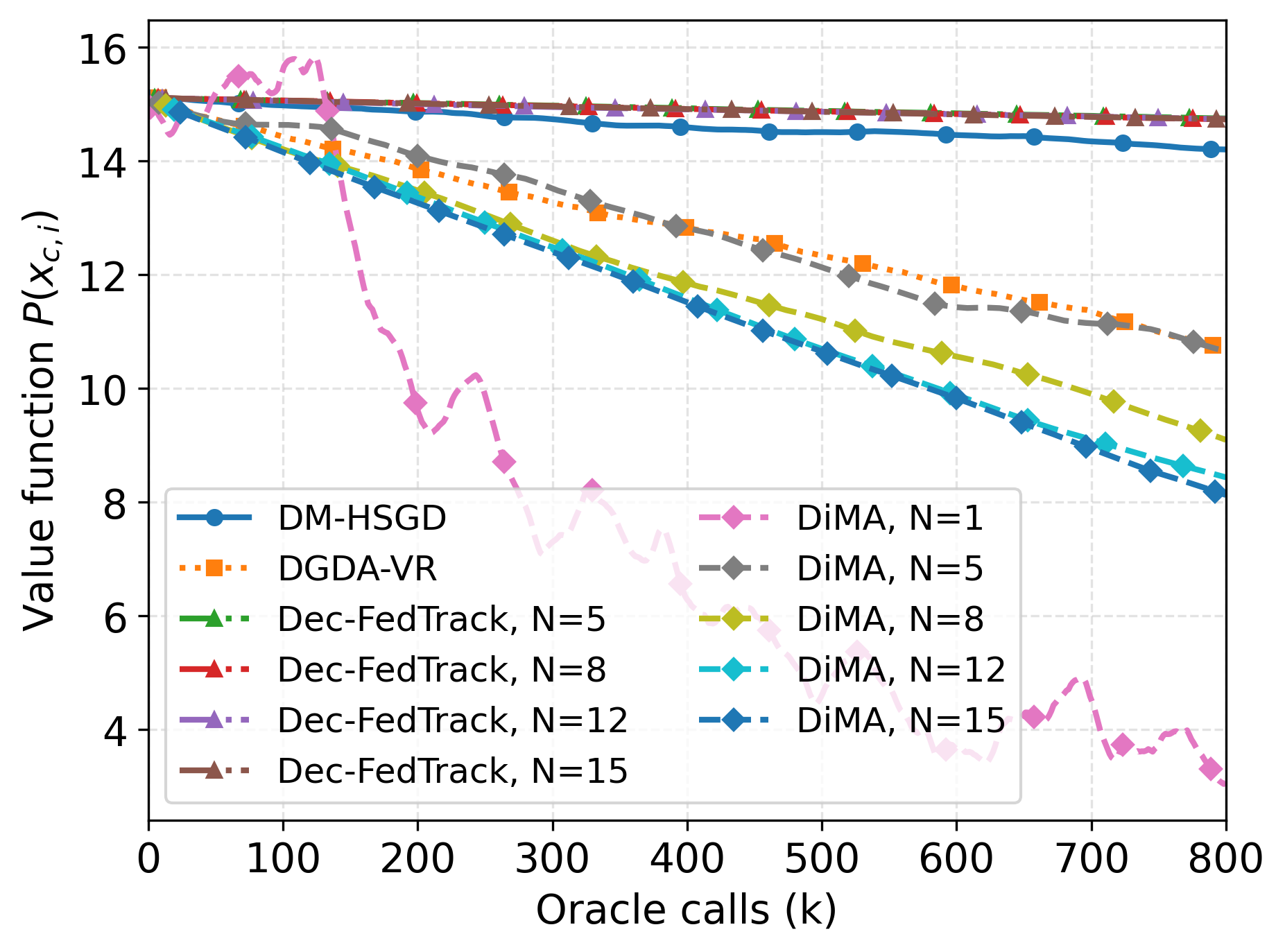}
\includegraphics[width=0.32\textwidth]{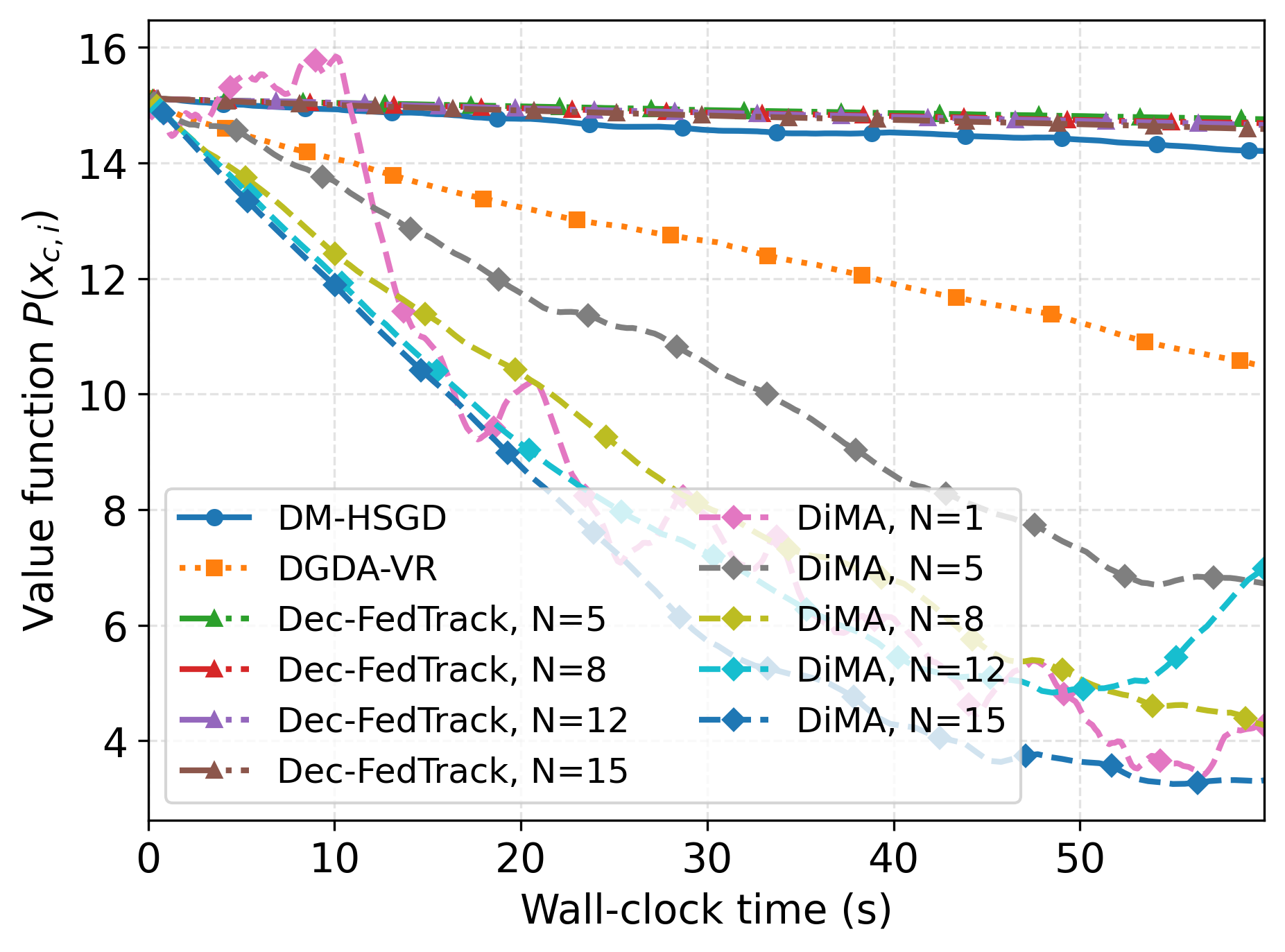}

\medskip
\textbf{(b) $\lambda=0.9999$}\par\smallskip
\includegraphics[width=0.32\textwidth]{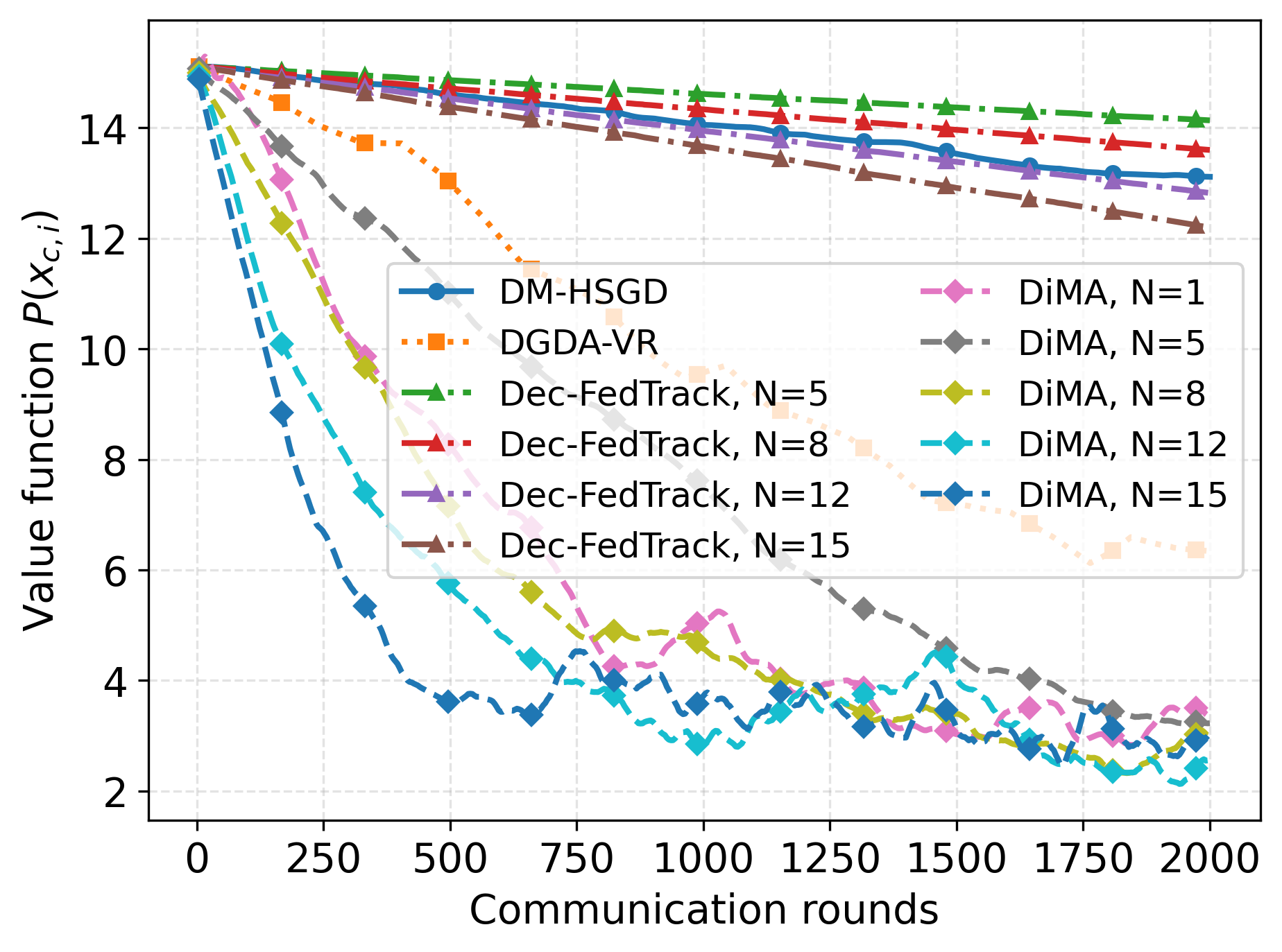}
\includegraphics[width=0.32\textwidth]{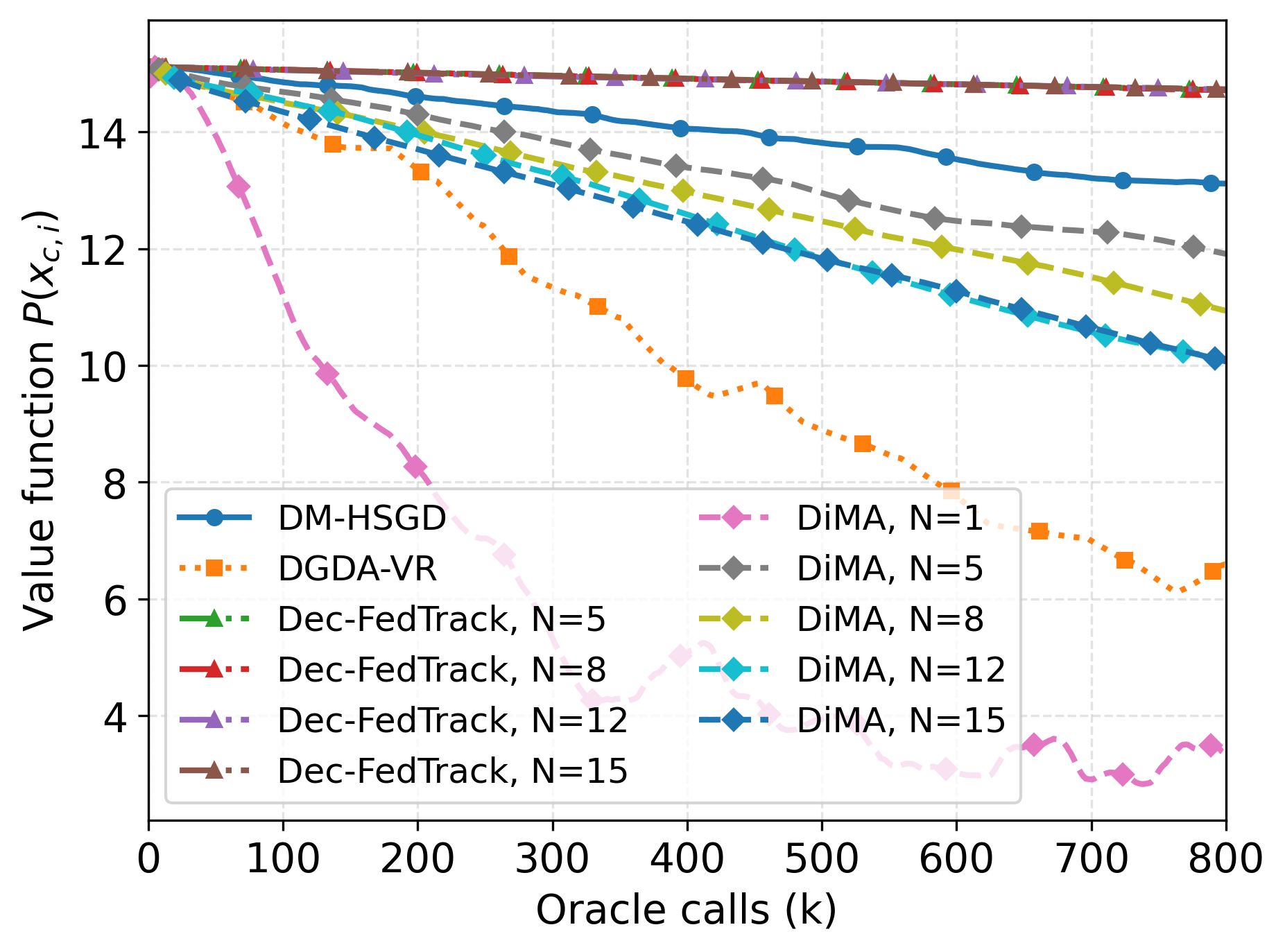}
\includegraphics[width=0.32\textwidth]{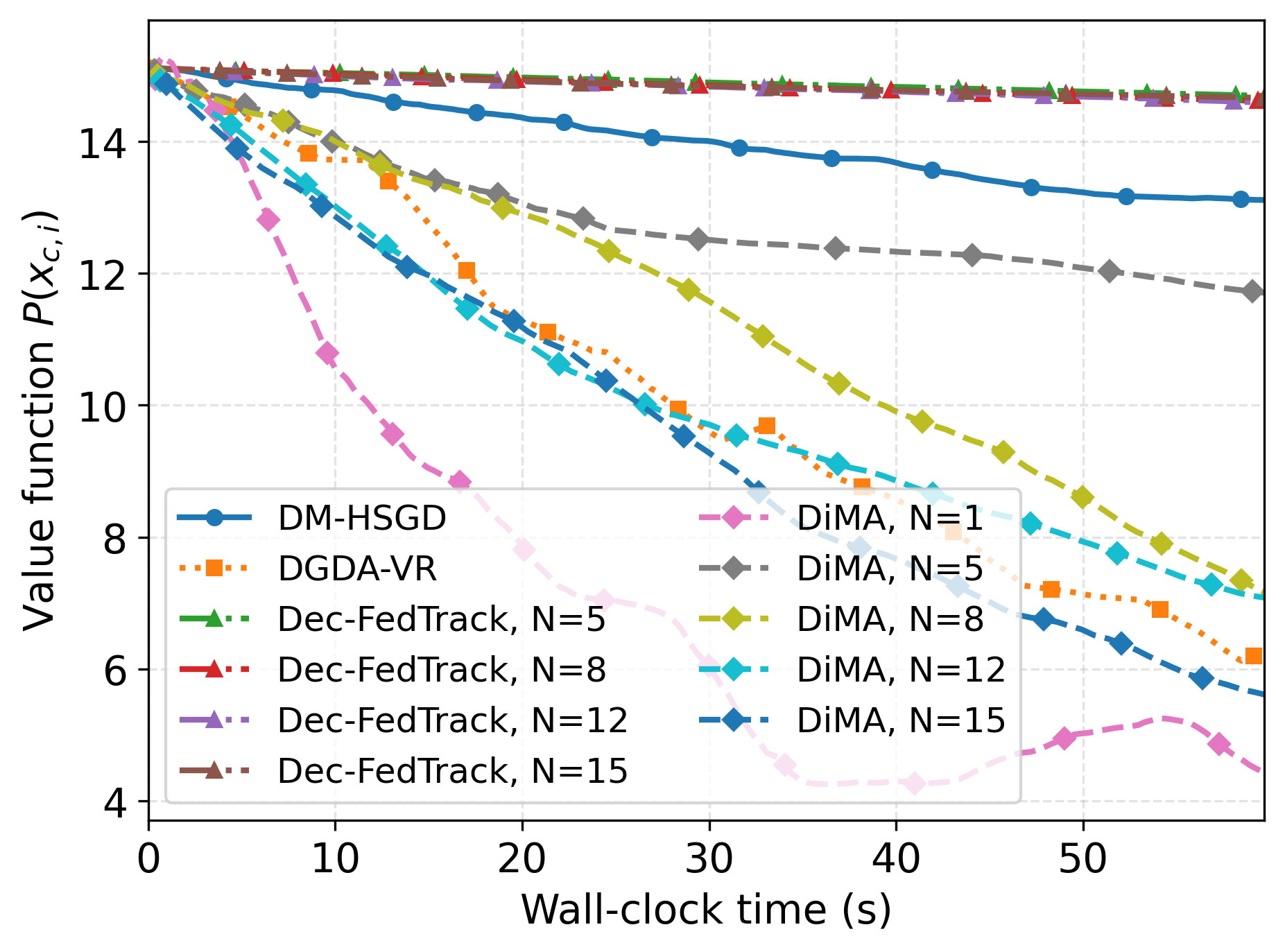}

\caption{Effect of network spectral gap on the ``a9a'' dataset.
The rows correspond to $\lambda=0.9$ and $\lambda=0.9999$.
The columns show the worst-case value $P(\bx_{c,i})$
versus communication rounds, oracle calls, and wall-clock time,
respectively.}
\label{figure:network:gap}
\vspace{-2em}
\end{figure*}

\begin{figure}[htbp]
\centering
\includegraphics[width=0.8\linewidth]{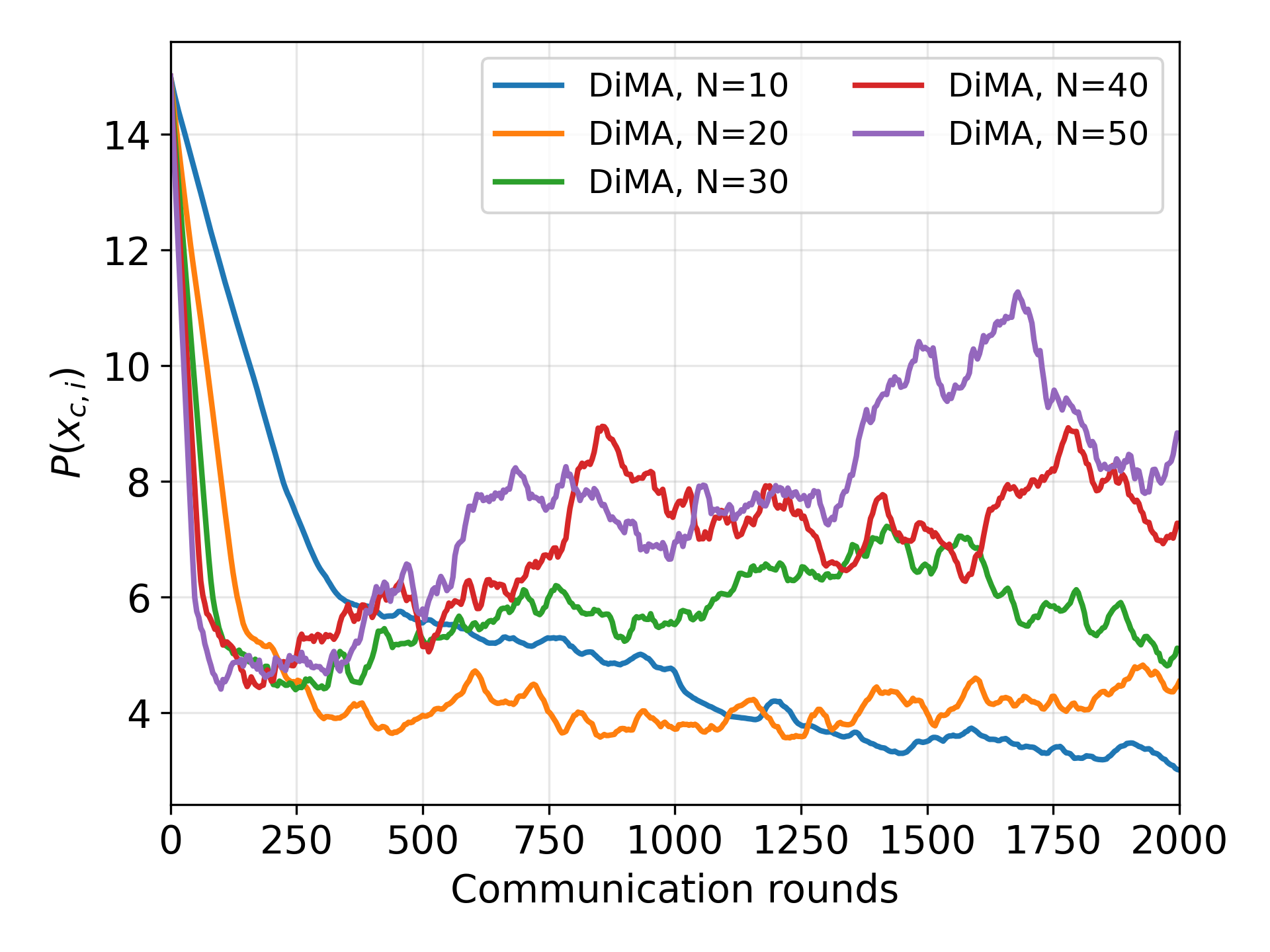}
\caption{Effect of the number of local steps $N$ on the dataset ``a9a".
We plot $P(\bx_{c,i})$ versus communication rounds.}
\label{figure:ablation:local}
\end{figure}
\begin{figure}[htbp]
\centering
\includegraphics[width=0.8\linewidth]{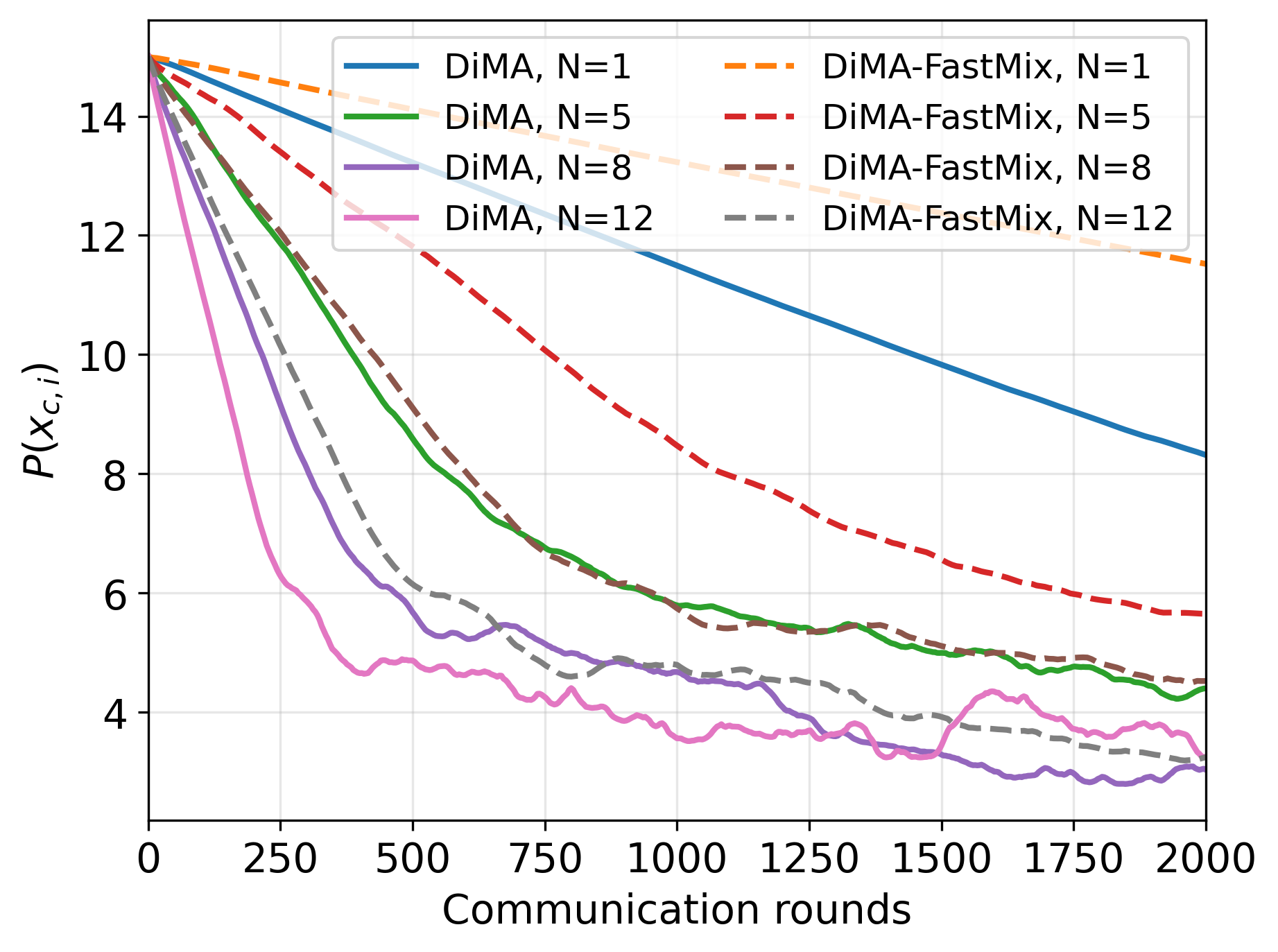}
\caption{Effect of accelerated mixing on the dataset ``a9a". We plot
$P(\bx_{c,i})$ versus communication rounds.
All gossip rounds are included in the communication count.}
\label{figure:ablation:fastmix}
\end{figure}
\textbf{Network spectral gap.}
We use the ``a9a'' dataset and consider the simplex-constrained problem. To study the effect of the spectral gap, we keep the ring topology fixed and adjust the mixing weights to obtain \(\lambda=0.9\) and \(\lambda=0.9999\), leaving all other settings unchanged. As shown in Fig.~\ref{figure:network:gap}, the smaller spectral gap has little effect on the performance of the proposed algorithms.
This may be due to the relative homogeneity of the training data and the use of momentum tracking in \textbf{DiMA}.

\textbf{Number of local steps.}
Using the same network factor \(\lambda=0.9999\), and similar experimental setup, we vary \(N\) over \(\{10,20,30,40,50\}\). Fig.~\ref{figure:ablation:local} shows that taking more local steps initially reduces the number of communication rounds needed to reach a given accuracy. Beyond a certain point, however, the local iterates may drift further away, offsetting the benefit of less frequent communication.

\textbf{Accelerated mixing.}
On the ring graph with \(\lambda=0.9999\), following the same experimental setup, we compare standard gossip with accelerated mixing, assuming \(\lambda\) is known and counting every gossip round. For simplicity, we choose $r=2$ for accelerated mixing. As shown in Fig.~\ref{figure:ablation:fastmix}, accelerated mixing does not always reduce the total number of communication rounds needed to reach a given accuracy. Accelerated mixing reduces disagreement by using more gossip rounds per outer iteration. It therefore saves communication only if it reduces the number of outer iterations enough to offset those extra rounds. For \(N=1\), standard gossip uses fewer communication rounds than accelerated mixing. This observation is consistent with our theoretical bounds: with $N=1$, the communication bound for standard gossip is better than that of accelerated mixing.

\ifCLASSOPTIONcaptionsoff
  \newpage
\fi

\end{document}